\documentclass[reqno,11pt]{amsart}
\usepackage{amsmath}
\usepackage[all]{xy}

\usepackage{amsthm,amsfonts, amssymb, amscd}
\usepackage{euscript}

\usepackage{tikz}
\usepackage{tikzit}
\usetikzlibrary{matrix,arrows}
\usetikzlibrary{positioning}

\usepackage[normalem]{ulem}

\usepackage{hyperref}
\hypersetup{%
    linktoc=page
}

\usepackage{todonotes}

\let\oldtocsection=\tocsection
\let\oldtocsubsection=\tocsubsection
\renewcommand{\tocsection}[2]{\hspace{0em}\oldtocsection{#1}{#2}}
\renewcommand{\tocsubsection}[2]{\hspace{1em}\oldtocsubsection{#1}{#2}}

\newtheorem{Thm}{Theorem}[section]
\newtheorem{Lem}[Thm]{Lemma}
\newtheorem{Cor}[Thm]{Corollary}
\newtheorem{Prop}[Thm]{Proposition}

\theoremstyle{remark}
\newtheorem{Rem}[Thm]{Remark}
\theoremstyle{remark}

\theoremstyle{definition}

\numberwithin{equation}{section}

\newcommand{\R}{\mathbb{ R}}           
\newcommand{\C}{\mathbb{C}}           
\newcommand{\Z}{\mathbb{ Z}}           

\newcommand{\ga}{\alpha}
\newcommand{\gb}{\beta}
\newcommand{\gd}{\delta}
\newcommand{\gD}{\Delta}

\newcommand{\gre}{\varepsilon}
\renewcommand{\gg}{\gamma}

\newcommand{\gl}{\lambda}

\newcommand{\cc}{\mathcal{C}}
 
 \newcommand{\ch}{\mathcal{H}}

 \newcommand{\cp}{\mathcal{P}}
 \newcommand{\cs}{\mathcal{S}}

 \newcommand{\cg}{\mathcal{G}}

 \newcommand{\cb}{\mathcal{B}}

\renewcommand{\tilde}{\widetilde}

\renewcommand{\bar}[1]{\overline{#1}}

\begin{document}
\parskip=4pt
\baselineskip=14pt

\title[Singular loci in type $\tilde A_2$]{Singular loci of Schubert varieties and the Lookup Conjecture in type $\tilde A_{2}$}

\author{Brian D.~Boe}
\address{
Department of Mathematics,
University of Georgia,
Athens, GA 30602
}
\email{brian@math.uga.edu}

\author{William Graham}
\address{
Department of Mathematics,
University of Georgia,
Athens, GA 30602
}
\email{wag@uga.edu}

\date{\today}

\begin{abstract}
We describe the loci of non-rationally smooth (nrs) points and of singular points for any non-spiral
Schubert variety of $\tilde{A}_2$ in terms of the geometry of the (affine) Weyl group action on the
plane $\R^2$.   Together with the results of Graham and Li for spiral elements, this allows us to
explicitly identify the maximal singular and nrs points in any Schubert variety of type $\tilde{A}_2$.
Comparable results are not known for any other infinite-dimensional Kac-Moody flag variety
(except for type $\tilde{A}_1$, where every Schubert variety is rationally smooth). 
As a consequence, we deduce that if $x$ is a point in a non-spiral Schubert
variety $X_w$, then $x$ is nrs in $X_w$ if and only if there are more than
$\dim X_w$ curves in $X_w$ through $x$ which are stable under the action of a maximal torus,
as is true for Schubert varieties in (finite) type $A$.  Combined with the work of Graham and Li for spiral
Schubert varieties, this implies the Lookup Conjecture for $\tilde{A}_2$.  
\end{abstract}

\maketitle

\tableofcontents

\section{Introduction}
The local topology of Schubert varieties at torus-fixed points has been of interest since the foundational paper \cite{KaLu:79} of Kazhdan and Lusztig connecting this
topology with representation theory.  The connection is simplest when the fixed point is rationally smooth, which holds automatically if the point is smooth, and there has been considerable interest
in understanding the loci of smooth and rationally smooth points in Schubert varieties (see e.g.\ \cite{BiLa:00}).  For example, the loci of non-smooth (i.e.\ singular)
and non-rationally smooth (nrs) points are closed,
and identifying the torus-fixed points in these loci which are maximal (in terms of the Bruhat order) would make it relatively easy to test whether any fixed point
is smooth or rationally smooth.  

In this paper, we study the loci of rationally smooth points and of smooth points in Schubert varieties in the flag variety for the Kac-Moody
group of type $\tilde{A}_2$.  We give explicit descriptions of these loci in terms of the action of the (affine) Weyl group $W$ on the plane $\R^2$.  In particular, the results of this paper, combined with the results of \cite{GrLi:21} and \cite{GrLi:15} for spiral elements, allow us
identify the maximal singular and nrs points in any Schubert variety of type $\tilde{A}_2$.  This provides a computationally efficient way to test whether any torus-fixed point is smooth or rationally smooth.  Comparable results are not known for any other infinite-dimensional Kac-Moody flag variety
(except for type $\tilde{A}_1$, where every Schubert variety is rationally smooth \cite{BiCr:12}).  By studying type $\tilde{A}_2$ in detail, we obtain insight which should be helpful in studying other Kac-Moody flag varieties.

Our methods allow us to complete the proof of the Lookup Conjecture (see \cite{BoGr:03}) for type $\tilde{A}_2$.  In particular, 
we prove that if $x$ is a point in a non-spiral\footnote{Here the spiral elements are as defined in \cite{GrLi:21}.  As discussed in that paper, the term spiral was adopted from \cite{BiMi:10}, but the definitions of \cite{GrLi:21} and \cite{BiMi:10} differ slightly.}
Schubert variety $X_w$, then $x$ is nrs in $X_w$ if and only if there are more than
$\dim X_w$ curves in $X_w$ through $x$ which are stable under the action of a maximal torus,
as is true for Schubert varieties in (finite) type $A$.  In other words, for non-spiral Schubert varieties of type $\tilde{A}_2$,
only the trivial case of
the Lookup Conjecture (see \cite{BoGr:03}) occurs.  For spiral Schubert varieties, the Lookup Conjecture was verified in \cite{GrLi:21}.  We conclude
that the Lookup Conjecture is true in type $\tilde{A}_2$.  

Since we can identify the loci of singular and nrs points in any Schubert variety in type $\tilde{A}_2$, we can identify the Schubert varieties for which these loci are empty -- that is, the Schubert varieties which are rationally smooth or smooth.  These Schubert varieties had previously been identified by Billey and Crites \cite{BiCr:12} in the case of $\tilde{A}_2$, in somewhat different terms.  See Corollaries \ref{c:ratsmoothSchubertvar} and \ref{c:smoothSchubertvar} and Remark \ref{r:BilleyCrites}.

It is well-known that the singular locus of a Schubert variety has codimension at least 2 (in other words, Schubert varieties are non-singular in codimension 1).  We prove a converse to this result in type $\tilde{A}_2$: if $X_w$ is not one of the 31 smooth Schubert varieties or 33 singular Schubert varieties of small dimension, then the singular locus of $X_w$ has codimension exactly 2.
Moreover, any Schubert subvariety of codimension at least 6 is contained in the singular locus.  If $X_w$ is not a rationally smooth Schubert variety, then the nrs locus of $X_w$ has codimension at most $4$.  More precise statements are given in Sections \ref{s:nrs} and \ref{s:smooth}.  
The smooth points of a non-spiral Schubert variety of type $\tilde{A}_2$ are identified in Theorem \ref{t:smoothLocus}.  We see 
 that a ``generic" Schubert variety has exactly $36$ smooth torus-fixed smooth points; for non-generic Schubert varieties, this number can be smaller.  

We now describe our results in more detail.  Let $\cg$ be a Kac-Moody group and $\cb$ a standard Borel subgroup, with
$T \subset \cb$ a maximal torus.  Let $W$ denote
the corresponding Weyl group.  The flag variety $X = \cg/\cb$ is not necessarily finite dimensional, but it can be written as a union of finite dimensional
algebraic varieties.  Corresponding to each $w \in W$ there is a subset $X^0_w = \cb w \cb$ of $X$, whose closure is the Schubert
variety $X_w$.  This is a finite dimensional algebraic variety whose dimension (over $\C$) is the length $\ell(w)$ of $w$.  The point $x \cb$ is in $X_w$
if and only if $x \leq w$ with respect to the Bruhat order on $W$.  The variety $X_w$ is the (finite) union of all $X_x^0$ for $x \leq w$ .  The local topology
of $X_w$ is the same at any point of a given $X_x^0$, so to determine the loci of smooth or rationally
smooth points in $X_w$, it suffices to determine which $T$-fixed points $x \cb$ are smooth or rationally smooth.

Non-rational smoothness can be detected by the following criterion, due to Carrell and Peterson, which is equivalent to a condition that had happeared
in Jantzen's work on highest weight modules (\cite{Car:94}, \cite{Jan:79}; see \cite[Theorem 12.2.14]{Kum:02}).
Given $x \leq w$, let $n^w_x$ denote the
number of reflections $r$ in $W$ such that $rx \leq w$, and let $q^w_x = n^w_x - \ell(w)$.  The numbers $q^w_x$ are all non-negative.
The point $x \cb$ is nrs in $X_w$ if and only if $q^w_y > 0$ for some $y \in W$ with $x \leq y \leq w$.

Given $w$ in $W$, we say that the Lookup Conjecture holds for $w$ (or for $X_w$), if, for all $x\le w$, $x \cb$ is nrs in $X_w$ if and only if $q^w_y > 0$ for $y = x$ or for some $x < y = rx \leq w$, where $r$ is a reflection in $W$.   The Lookup Conjecture states that the Lookup Conjecture holds for all $w \in W$.
If true, the Lookup Conjecture would greatly reduce the number of $q^w_y$ that
need to be computed to detect non-rational smoothness.  

The nontrivial case of the Lookup Conjecture occurs when $x \cb$ is nrs in $X_w$ and $q^w_x = 0$. The trivial case
occurs when $x \cb$ is nrs in $X_w$ and $q^w_x > 0$.  The statement that only the trivial case of the Lookup Conjecture holds for $w$ means that $x \cb$ is nrs in $X_w$ if and only if $q^w_x > 0$.  If this holds, the Lookup Conjecture holds for $w$.
 In (finite) type $A_n$, only the trivial case of the Lookup Conjecture occurs for any $w \in W$ (see \cite{Deo:85}), and so the Lookup Conjecture
is true in type $A_n$. 

 In type $\tilde{A}_2$, Graham and Li \cite{GrLi:21} proved that the Lookup Conjecture holds (nontrivially) for so-called spiral elements $w \in W$: those which have a reduced expression $s_{i}s_{j}s_{k}s_{i}s_{j}s_{k}s_{i} \dots$, where $i, j, k$ are distinct.  To prove the Lookup Conjecture
 in type $\tilde{A}_2$, it therefore suffices to prove that if $w$ is not spiral, then $x \cb$ is nrs in $X_w$ if and only if $q^w_x > 0$, 
 as had been observed by Graham and Li
in examples. One of our main results is that this is true, and therefore the Lookup Conjecture holds in type $\tilde{A}_2$.  

As in \cite{GrLi:21}, we make use of the action of $W$ on a Euclidean space $V \cong \R^2$.  Let
$\Phi$ denote the finite root system of type $A_2$, so $\Phi$ consists of the roots $\ga_1, \ga_2, \tilde{\ga} = \ga_1 + \ga_2$,
and their negatives.  Corresponding to each $\ga \in \Phi$  and each $k \in \Z$, there is an affine root hyperplane (in fact, a line)
$H_{\ga,k}$ in $V$, consisting of the $v \in V$ whose inner product with $\ga$ is $k$.  The Weyl group $W$ is the group of (affine) transformations
of $V$ generated by the (affine) reflections $s_{\ga,k}$ in the $H_{\ga,k}$.
The connected components of the complement of
the union of the root hyperplanes are called alcoves.  

We fix an alcove $A_0$, which we call the fundamental alcove, and let $q$ denote the center
of $A_0$.
The group $W$ acts simply transitively on the set of alcoves, or equivalently,  on the
set of alcove centers, so we have bijections
$$
W \leftrightarrow \{ \mbox{alcoves} \} \leftrightarrow \{ \mbox{alcove centers} \} 
$$
$$
w \leftrightarrow w A_0 \leftrightarrow wq.
$$
The spiral elements correspond to the alcoves lying in regions we call fundamental root strips.  The complement of the union of these strips has 6 components, which we call chambers.  Unlike Weyl chambers, these are not all equivalent; rather, they are of two types, which we call even and odd.  The details are in Section \ref{s.preliminaries}.

Our results rely on a description of the Bruhat order in terms of alcove geometry.  We identify Weyl group elements with alcove
centers, as above.  Given a non-spiral element $w \in W$, we show that the set of $x \in W$ such that $x \leq w$ coincides with the
set $\ch_w$ of elements lying in a convex hexagon, which we define
explicitly by specifying its vertices (see Theorem \ref{t.hexagon}).  Our methods can be used to recover the analogous result for spiral elements from \cite{GrLi:21}, with a simpler proof.
The hexagon is not regular, but it has greater symmetry if $w$ lies in an even chamber.  As a consequence, many of our results have simpler
statements for such $w$ than for $w$ in an odd chamber.  
After proving Theorem \ref{t.hexagon}, we learned that this result, as well as the analogous result for spiral
elements, were conjectured in \cite{Jit:23}.  That paper also contains an outline of a proof of the conjecture, which makes use of ideas similar to some ideas in our proof of Theorem \ref{t.hexagon}.

Using the Bruhat hexagon $\ch_w$, we describe the set of $x$ such that $q^w_x >0$, and deduce that if $x \leq y \leq w$ and $q^w_y > 0$, then 
$q^w_x>0$.  This implies that $x \cb$ is nrs in $X_w$ if and only if $q^w_x>0$.  The proofs are inductive.  They make use of the relationship
between $\ch_w$ and $\ch_{w'}$, where $w'$ is obtained from $w$ by a simple reflection or by ``translation into the chamber."  In fact, we obtain a fairly explicit description of $q^w_x$ in terms of the location
of $x$ in $\ch_w$.  We identify the maximal nrs points of $X_w$, and deduce that if $X_w$ is nrs, then the codimension of the nrs locus is either 3 or 4.

Not all the rationally smooth points in $X_w$ are smooth.  To determine the smooth locus, as in \cite{GrLi:15}, we use Kumar's criterion
in terms of equivariant multiplicities $e^w_x$.  These geometric invariants contain more information than the integers $q^w_x$, but they
are considerably more difficult to calculate.  However, knowing the nrs locus in $X_w$ allows us to determine the smooth locus by
calculating $e^w_x$ for a small number of $x$ for which $\ell(w) - \ell(x) \leq 3$.  This computation turns out to be manageable.  We describe
the smooth locus of a non-spiral $X_w$ by identifying 
six elements of $\ch_{w}$ near each vertex of the hexagon, such that $x \cb$ is a smooth point of $X_w$ if and only if $x$ is one of these elements.
As a consequence, we see that any Schubert variety of type $\tilde{A}_2$ has at most 36 smooth torus-fixed points.
We identify the maximal singular points, and conclude that except for finitely many Schubert varieties of dimension at most 9, 
 the singular locus of $X_w$ has codimension 2.  Moreover, any Schubert subvariety of $X_w$ of codimension 7 or more is contained in the singular locus.


\section{Affine Weyl groups and Schubert varieties} \label{s.preliminaries}
In this section we summarize some facts about affine Weyl groups and Schubert varieties for Kac-Moody groups, and prove an extension of the
Setup Move of \cite{BoGr:03} to equivariant multiplicities.  Our main references are
\cite{Hum:90}, \cite{BjBr:05}, \cite{Kum:02}, and \cite{BoGr:03}; some of the exposition is from \cite{GrLi:21}.  

\subsection{Affine Weyl groups} \label{ss.affineWeyl}
Let $V$ be a Euclidean space; that is, a real vector space equipped with a positive definite
inner product $(\cdot, \cdot)$.  We will use the inner product to identify $V$ with $V^*$, and
view both roots and coroots as elements of $V$.  Let $\Phi \subset V$ be an irreducible root system, and 
$\Phi^{\vee} \subset V$ the dual root system;
the coroot $\ga^{\vee}$ is related to the root $\ga$ by
$\ga^{\vee} = 2 \ga /(\ga, \ga)$.  
If there is only one root length, we scale the inner product so that $(\ga, \ga) = 2$ for each root,
and then $\Phi= \Phi^{\vee}$.  Let $L(\Phi^{\vee}) \subset V$ denote the abelian group
generated by $\Phi^{\vee}$.  We will identify $L(\Phi^{\vee})$ with the corresponding group of translations of $V$; the translation corresponding to $\gamma \in L(\Phi^{\vee})$ is denoted by $t(\gamma)$.

Given $\alpha \in \Phi$ and $k\in\Z$, let $H_{\ga, k} = \{ v \in V \mid (\ga, v) = k \}$, called a root hyperplane.  Let $s_{\ga, k}: V \to V$ denote the map given by reflection across this hyperplane;
we write $s_{\ga} = s_{\ga, 0}$, and then $s_{\ga, k} = t(k \ga^{\vee})s_{\ga}$.
The affine Weyl group $W$ associated to $\Phi$ is the
group of affine transformations of $V$ generated by the elements $s_{\ga, k }$, and the finite Weyl group
$W_{f}$ is the subgroup of $W$ generated by the $s_{\ga}$.  The group $W_f$ acts on $L(\Phi^{\vee})$ by $w \cdot t(\gamma) = t(w(\gamma))$, and $W$ 
 can be identified
with the semidirect product $L(\Phi^{\vee}) \rtimes W_f $.   We have the following useful formulas:
\begin{equation} \label{e.usefulaffine}
w s_{\gb} w^{-1} = s_{w \gb} \quad \text{and} \quad t(\lambda) s_{\gb, k} t(- \lambda) = s_{\gb, k+(\lambda,\beta)}
\end{equation}
for $\beta \in \Phi$, $w \in W_f$, $k\in\Z$, and $\lambda \in V$
satisfying $(\gl, \ga) \in \Z$ for all $\alpha \in \Phi$ (see \cite[Prop.~4.1]{Hum:90}).  
We will later use the corollary that if $(\beta, \lambda )=k$, then
\begin{equation} \label{e.betaktga}
s_{\beta,k}t(\lambda)=t(\lambda )s_{\beta}.
\end{equation}

Choose a set of simple roots $\{ \ga_1, \ldots, \ga_{n} \}$ for $\Phi$, and write $s_i = s_{\ga_i}$.
 Let $\tilde{\ga}$ denote the highest root and let $s_0 = s_{\tilde{\ga},1}$.  The elements
 $s_0, s_{1}, \ldots, s_n$ are called simple reflections, and they generate $W$ (even as a Coxeter group).  
 The length $\ell(w)$ of an element $w \in W$ is equal to $\ell$ if $w$ can be written as a product
 of $\ell$ simple reflections, but no fewer.  Given $x,y \in W$, let $H(x,y)$ denote the set of hyperplanes $H_{\ga, k}$ separating $x A_0$ and $y A_0$, for $\ga \in \Phi$ and $k \in \Z$. 
Then $\ell(w) = | H(w,e) |$, where $e$ denotes the identity element of $W$ \cite[Theorem 4.5]{Hum:90}.
 
  The reflections in $W$ are defined to be the conjugates in $W$ of the $s_i$; if there is only one root length, then
 the reflections are exactly the elements
 $s_{\ga,k}$ for $\ga \in \Phi$ and $k \in \Z$ (see for example \cite[Section 2]{GrLi:21}).
 Let $R$ denote the set of reflections in $W$.  The Bruhat order on $W$ is the partial order generated
 by the relations $w < rw$ where $r \in R$ and $\ell(w) < \ell(rw)$.  

Given $w \in W$, let $R(w)$ (resp.~$L(w)$)  denote the subgroup of $W$ generated by the simple reflections $s$ such that
$w s < w$ (resp.~$sw < w$).  The group $R(w)$ (resp.~$L(w)$) is isomorphic to the finite Weyl group corresponding to the Coxeter graph
obtained from the Coxeter graph of $W$ by deleting the simple reflections not in $R(w)$ (resp.~$L(w)$).  If $w$ is not the identity
element $e$, then $R(w)$ contains at least one simple reflection.  If $\Phi$ is of type $A_2$, then $R(w)$ contains exactly one simple reflection $s$, in which case it is the group $\{e,s \}$, or $R(w)$ contains two simple
reflections $s,t$, in which case it is isomorphic to the Weyl group of type $A_2$ generated by $s$ and $t$.

The alcoves are the connected components of $V\, \setminus\, \bigcup_{\ga \in \Phi, k\in \Z}\ H_{\ga ,k}$. 
The alcove bounded by the hyperplanes $H_{\ga_1,0},  \dots, H_{\ga_n,0}, 
H_{\tilde \alpha,1}$
is called the fundamental alcove, and denoted
$A_{\circ}$. Thus 
$A_{\circ}=\{v\in V \mid (v,\alpha_i)> 0, i = 1, \ldots, n,  \mbox{ and } (v,\tilde \alpha)< 1 \}$. 
Let $q$ denote the center of the alcove $A_{\circ}$, so $wq$ is the center of the alcove $w A_{\circ}$, $w\in W$. The Weyl group $W$ acts simply transitively on the set of alcoves, so there are bijections
$W\leftrightarrow \{ \mbox{alcoves} \} \leftrightarrow Wq $ given by 
$w\leftrightarrow wA_{\circ}  \leftrightarrow wq$.
Using these bijections, we will frequently identify elements of $W$ with alcoves and alcove centers.

\subsection{Schubert varieties} \label{ss:Schubert}
Associated to the root system $\Phi$ there is a Kac-Moody group $\cg$ of affine type
(over the ground field $\C$)
with Borel subgroup $\cb$, whose Weyl group is the affine Weyl group $W$ associated
to $\Phi$.
The flag variety is $X = \cg/\cb$, which has the structure
of a projective ind-variety.  Given $w \in W$, the Schubert variety $X_w$ 
is defined as the union
$\cup_{x \leq w} \cb x \cb / \cb \subset X$; it is the closure of the $\cb$-orbit $X_{w}^{0} = \cb \cdot w \cb$,
and has the structure of a finite dimensional projective variety
of dimension $\ell(w)$.  We will be concerned with the loci of smooth and rationally smooth points in $X_w$.  Rational smoothness is a topological condition related to smoothness, but weaker: a smooth point is rationally smooth, but the converse can fail.  For the definition of rational smoothness, see for example
\cite[Definition 12.2.7]{Kum:02}.  We will not need this definition; our main tools for determining rational smoothness and smoothness
will be (respectively) the Carrell-Peterson criterion described in the introduction, and Kumar's smoothness criterion \cite[Theorems~12.2.14, 12.1.11]{Kum:02}.  For brevity, we write nrs to mean not rationally smooth.

We have $x \cb \in X_w$ $\Leftrightarrow$ $X_x \subseteq X_w$ $\Leftrightarrow$
$x \leq w$.  The singular locus and nrs locus in $X_w$ are closed and $\cb$-invariant,
so they are unions of Schubert varieties; 
moreover, $X_x$ is in the singular or nrs locus of $X_w$ if and only if the point
$x \cb$ is.  Hence, if $x \leq y \leq w$ and the point $y \cb$ is singular (resp.~nrs) in $X_w$, then so $x \cb$.  Hence, to determine the loci of singular or nrs points in $X_w$, it suffices to determine which points $x \cb$ are singular or nrs in $X_w$.  
 The singular locus of any Schubert variety $X_w$ has codimension at least $2$ \cite[Prop.~12.1.1]{Kum:02}, so $w \cb$ is smooth in $X_w$, and if $x < w$ and $\ell(x) = \ell(w) -1$, then $x \cb$ is smooth in $X_w$.  

If $x \leq w$ and $u \in R(w)$, there is a close relationship between invariants of $X_w$ at the points $x \cb$ and $xu \cb$,
 for example, the values of the integers $q^w_\bullet$,
rational smoothness, and smoothness.  There are analogous relationships for $u \in L(w)$.  We use the term ``Simple Move" (adapted from \cite{BoGr:03}) to refer
to such a relationship.

The next proposition states some of the simple move relationships for $u \in R(w)$.  The proposition is an extension of results from \cite[Section 4]{BoGr:03}; see also \cite{GrLi:21} and \cite{GrLi:15}.  There is an analogous version for $u \in L(w)$, which we will also use.

\begin{Prop} \label{p.summary} 
Let $w \in W$ and  $u \in R(w)$.
\begin{enumerate}
\item[(a)] If $x \in W$, then $x\le w\Leftrightarrow xu \le w$.

\item[(b)] If $x \leq w$, then $q^w_x=q^w_{xu}$.

\item[(c)] Let $x \leq w$.  The point $x \cb$ is rationally smooth (resp.~smooth) in $X_w$ $\Leftrightarrow$ the point $xu \cb$ is rationally smooth (resp.~smooth) in $X_w$.

\item[(d)] Let $x \leq w$.  The Lookup Conjecture holds for the pair $x \leq w$ $\Leftrightarrow$ the Lookup Conjecture holds for the pair $xu \leq w$.
\end{enumerate}
\end{Prop} 

\begin{proof} The proposition holds if $u$ is a simple reflection, by \cite{BoGr:03} or \cite{GrLi:21}, except for the assertion about smoothness, which
is \cite[Proposition 2.5(d)]{GrLi:15}.  In general, we can write
$u$ as a product of simple generators of $R(w)$; the result then follows by induction on $\ell(u)$.
\end{proof}

\subsection{The Setup Move} \label{ss.setup}
The Setup Move was introduced in \cite{BoGr:03} to study rational smoothness.
In this paper, we will use the Setup Move 
to help determine the smooth locus of $X_w$.  
Suppose that 
\begin{equation} \label{e:setup}
x \leq w \text{ and } s \text{ is a simple reflection such that } w < ws \text{ but } xs \not\leq w.
\end{equation}
In this situation, we say $x,w,$ and $s$ satisfy the hypotheses of the Setup Move (on the right); in this case, $xs < ws$ by \cite[Lemma 2.3]{BoGr:03} (the Maximum Principle).
We will say that the pair $(xs, ws)$ is obtained from the pair $(x,w)$ by the Setup Move.
There is also a Setup Move on the left, where the simple reflection $s$ multiplies on the left, and the analogous results hold.   

\begin{Prop} \label{p:setup}
Suppose $x,w$ and $s$ satisfy \eqref{e:setup}.  The point $x \cb$ is rationally smooth (resp.~smooth) in $X_w$ $\Leftrightarrow$ 
the point $xs \cb$ is rationally smooth (resp.~smooth) in $X_{ws}$.
\end{Prop} 

The assertion about rational smoothness is proved exactly as in 
 \cite[Prop.~4.12]{BoGr:03}.  To prove the assertion about smoothness, we need to use Kumar's criterion
for smoothness in terms of equivariant multiplicities.  
This is the only point in the paper where we use equivariant multiplicities, and we will need to briefly recall some background about these.
The reader who is willing to accept this proposition can skip the rest of this section.  

The real roots of the Kac-Moody group 
are certain linear combinations of the simple roots $\gb_0, \gb_1, \ldots, \gb_n$.
Here $\gb_i = \ga_i$ for $i = 1, \ldots, n$, 
where $\ga_1, \ldots, \ga_n$ are simple roots for the finite
Lie algebra, together with an additional root $\gb_0$. 
(We introduce the notation $\gb_i$, rather than using $\ga_i$, because we wish to reserve the notation
$\ga_0$ for the longest finite root.)
  The Weyl group acts on the set of real roots;
this action is characterized by the formula $s_i(\gb_j) = \gb_j - a_{ij} \beta_i$ for $i\ne j$ in $\{0,1, \ldots, n\}$. 
Here the $a_{ij}$ are the entries of the generalized Cartan matrix constructed from the root system $\Phi$; see
\cite[Section 13.1]{Kum:02}.  If $\Phi$ is of type $A_2$, the case of interest in this paper, then $n=2$ and $a_{ij} = -1$ for $i \neq j$,
so $s_i(\beta_j) = \beta_j + \beta_i$ for $i \neq j$.
Corresponding to each real root $\gb$ there is a reflection $s_{\gb}$.  Under our identification of $W$ with transformations of $\R^n$,
each $s_{\gb}$ corresponds to $s_{\ga,n}$ for some $\ga \in \Phi$ and some $n \in \Z$.  For type $\tilde{A}_2$, this correspondence
is described in \cite[Prop.~2.6]{GrLi:15}.
Let $\Psi^w_x$ be the set of positive real roots $\gb$ such that $s_{\gb} x \leq w$.  We denote by $\prod \Psi^w_x$ the product of the elements in $\Psi^w_x$.

The following result is analogous to \cite[Prop.~4.11]{BoGr:03}, and has a similar proof, which we omit.

\begin{Prop} \label{p:setup-Psi}
Suppose $x,w$ and $s = s_{\beta}$ satisfy \eqref{e:setup}.  Then
\begin{equation} \label{e.Setup-Psi}
\Psi^{ws}_{xs} = \Psi^w_x \sqcup \{ x \beta \}.
\end{equation}
\end{Prop}

Given $x \leq w$ in $W$, the equivariant multiplicity
$e^w_x$ is a rational function in the $\gb_i$, defined as follows.  Write $s_i = s_{\gb_i}$.  Let $\cs = (s_{i_1}, \ldots, s_{i_\ell})$ be a fixed reduced expression
for $w$; that is, $s_{i_{1}}\cdots s_{i_{\ell}}=w$ and $\ell=\ell(w)$ .
A subexpression of $\cs$ is a sequence $\sigma = (\sigma_1, \ldots, \sigma_\ell)$ such that each $\sigma_j$ equals either $e$ or $s_{i_j}$.
Let $\cs(x)$ denote the set of subexpressions $\sigma$ of $\cs$ such that $\sigma_1 \sigma_2 \cdots \sigma_\ell = x$ (``$\sigma$ multiplies to $x$'').  Then
$$
e^w_x = (-1)^{\ell(w)} \sum_{\sigma \in \cs(x)}\ \prod_{j=1}^\ell \frac{1}{\sigma_1 \cdots \sigma_j (\beta_{i_j})} 
$$
(see \cite[Theorem 11.1.2]{Kum:02}).  The point $x \cb$ is smooth in $X_w$ if and only if $e^w_x = {(-1)^{\ell(w) - \ell(x)}} / {\prod \Psi^w_x}$ \cite[Theorem 12.1.11]{Kum:02}.

{\em Proof of Proposition \ref{p:setup}.}
As noted immediately after the statement of the proposition, we need only prove the assertion about smoothness.  We assume $x,w$ and $s = s_{\gb}$ satisfy \eqref{e:setup}, with $\cs$ as above.  Observe that
\begin{equation} \label{e.Setup-equiv}
e^{ws}_{xs} = \frac{1}{x \gb} e^w_x.
\end{equation}
Indeed, $\cs' = (s_{i_1}, \ldots, s_{i_\ell}, s)$ is a reduced expression for $ws$.  Moreover, since $xs \not\leq w$, and in particular $x<xs$,
$\cs'(xs)$ consists of the subexpressions $(\sigma_1, \ldots, \sigma_\ell, s)$, where $(\sigma_1, \ldots, \sigma_\ell) \in \cs(x)$.  Equation \eqref{e.Setup-equiv}
follows from these observations and the definition of equivariant multiplicities. 

It follows from Proposition \ref{p:setup-Psi} that
\begin{equation} \label{e.prod-Setup-Psi}
\prod \Psi^{ws}_{xs} = x \beta \prod \Psi^w_x.
\end{equation}
By Kumar's criterion, $x \cb$ is smooth in $X_w$ if and only if $e^w_x = \frac{(-1)^{\ell(w) - \ell(x)}}{\prod \Psi^w_x}$.  
By \eqref{e.Setup-equiv}, \eqref{e.prod-Setup-Psi}, and Kumar's criterion again, this holds if and only if
$$
e^{ws}_{xs} =  \frac{1}{x \gb} \cdot \frac{(-1)^{\ell(w) - \ell(x)}}{\prod \Psi^w_x} = \frac{ (-1)^{\ell(ws) - \ell(xs)} }{ \prod \Psi^{ws}_{xs} } \iff xs \cb \text{ is smooth in } X_{ws}.
$$
This proves the proposition. \hfill $\Box$

\begin{Rem} \label{r.Setup-left}
As noted above, the analogue of Proposition \ref{p:setup} holds for the Setup Move on the left.   Under the hypotheses of the Setup Move on the left,
the analogue of Proposition \ref{p:setup-Psi} is
\begin{equation} \label{e.Setup-Psi-left}
\Psi^{sw}_{sx} = s(\Psi^w_x) \sqcup \{ \gb \},
\end{equation}
and the analogue of \eqref{e.Setup-equiv} is
\begin{equation} \label{e.Setup-equiv-left}
e^{sw}_{sx} = \frac{1}{\gb} s(e^w_x).
\end{equation}
\end{Rem}

\begin{Rem} \label{r.Simple-equiv}
There are analogues of Proposition \ref{p:setup-Psi} and equation \eqref{e.Setup-equiv} for the Simple Move (see \cite[Prop.~2.5]{GrLi:15}).
\end{Rem}

\begin{Rem} \label{r.subexpression}
Given $x$ and $w$ in $W$, it is well-known that $x \leq w$ if and only if any reduced expression for $w$ contains a subexpression multiplying to $x$
(see \cite[Theorem 5.10]{Hum:90}).  This connection between the Bruhat order and subexpressions plays a key role in the proof of Theorem \ref{t.hexagon}, so to emphasize it,
we will sometimes use the phrase ``$x$ is a subexpression of $w$" to mean $x \leq w$.
\end{Rem}

\section{Root strips, root strings and chambers} \label{s:chambers}
In this section, we define root strips, root strings and chambers, which play a major role in this paper.
For the remainder of the paper, $W$ will denote the affine Weyl group of type $\tilde{A}_2$.  

\subsection{Definitions and length results} 
This section contains definitions, as well as some results about how the length of an element changes under various operations, which will be needed later.

Fix $\ga \in \{\ga_1, \ga_2, \tilde\ga\ =: \ga_{0}\}$,
and $k\in\Z$. An $\ga$ root string is a line parallel to $\ga$ passing through some alcove center $wq$.  

The region $\{ v \in V \mid k \le (\ga,v) \le k+1 \}$ is called an $\ga$ root strip. 
It is a union of alcove closures. We will also sometimes refer to the set of centers of those alcoves, or the corresponding elements of $W$, as a root strip. The fundamental root strips are those with $k=0$. The elements $w$ belonging to the fundamental root strips are the spiral elements of \cite{GrLi:21}. For the purposes of this paper, we will refer to chambers as the connected components of the complement of the union of the fundamental root strips (or sometimes, their closures).  The boundary of a chamber is the union of two rays along root hyperplanes (with $k=0$ or 1); we will refer to these root hyperplanes (or sometimes, the rays) as the walls of the chamber.  
\begin{figure}[htbp!]
%

\pgfmathsetmacro{\cols}{5}
\pgfmathsetmacro{\rows}{6}
\pgfmathsetmacro{\slant}{cot(60)}
\pgfmathsetmacro{\height}{0.5 * \rows * tan(60)}
\pgfmathsetmacro{\triht}{sin(60)}
\pgfmathsetmacro{\upmid}{0.25 * sec(30)}
\pgfmathsetmacro{\downmid}{\triht - \upmid}

\begin{tikzpicture}[every label/.append style = {fill=white}]
    

\input{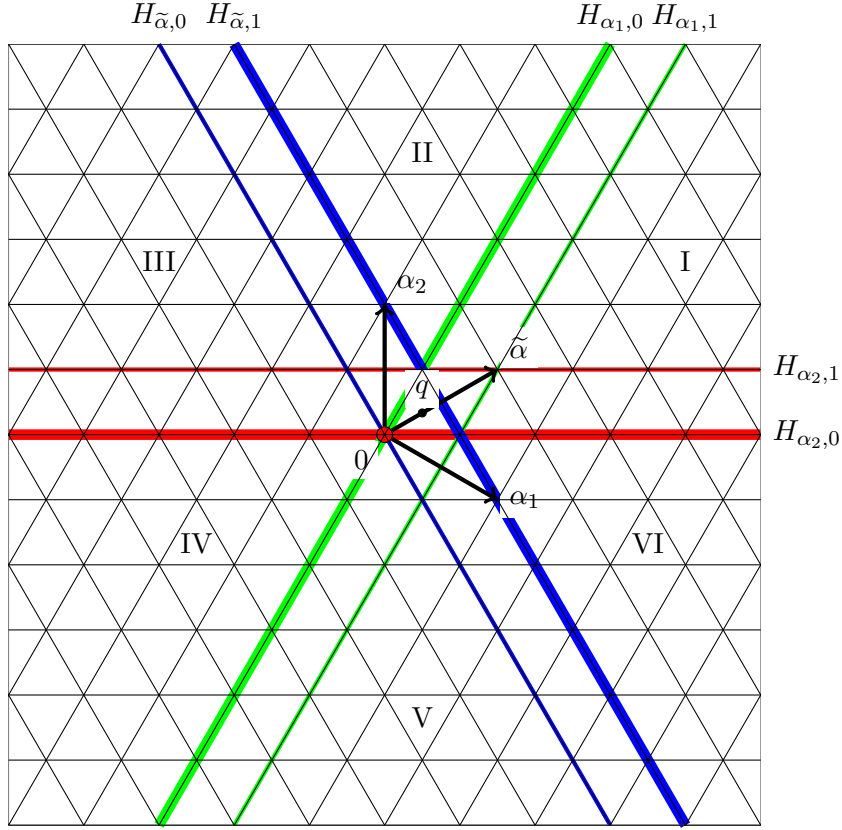}

	\begin{pgfonlayer}{nodelayer}
		\node [style=none] (0) at (-5, 0) {};
		\node [style=none, label={right:$H_{\alpha_{2},0}$}] (1) at (5, 0) {};
		\node [style=none] (2) at (-5, 1*\triht) {};
		\node [style=none, label={right:$H_{\alpha_{2},1}$}] (3) at (5, 1*\triht) {};
		\node [style=red dot, label={below left:$0$}] (4) at (0, 0) {};
		\node [style=none] (5) at (-3, -6*\triht) {};
		\node [style=none, label={above:$H_{\alpha_{1},0}$}] (6) at (3, 6*\triht) {};
		\node [style=none] (7) at (-2, -6*\triht) {};
		\node [style=none, label={above:$H_{\alpha_{1},1}$}] (8) at (4, 6*\triht) {};
		\node [style=none] (9) at (3, -6*\triht) {};
		\node [style=none, label={above:$H_{\tilde\alpha,0}$}] (10) at (-3, 6*\triht) {};
		\node [style=none] (11) at (4, -6*\triht) {};
		\node [style=none, label={above:$H_{\tilde\alpha,1}$}] (12) at (-2, 6*\triht) {};
		\node [style=small black dot, label={above:$q$}] (13) at (0.5, \upmid) {};
		\node [style=none, label={right:$\alpha_1$}] (14) at (1.5, -1*\triht) {};
		\node [style=none, label={above right:$\alpha_2$}] (15) at (0, 2*\triht) {};
		\node [style=none, label={above right:$\widetilde\alpha$}] (16) at (1.5, 1*\triht) {};
		\node [style=none] (17) at (4, 2*\triht + \downmid) {I};
		\node [style=none] (18) at (0.5, 4*\triht + \upmid) {II};
		\node [style=none] (19) at (-3, 2*\triht + \downmid) {III};
		\node [style=none] (20) at (-2.5, -2*\triht + \upmid) {IV};
		\node [style=none] (21) at (0.5, -5*\triht + \downmid) {V};
		\node [style=none] (22) at (3.5, -2*\triht + \upmid) {VI};

	\end{pgfonlayer}
	\begin{pgfonlayer}{edgelayer}
		\draw [style=thick red line] (0.center) to (1.center);
		\draw [style=red line] (2.center) to (3.center);
		\draw [style=thick green line] (5.center) to (6.center);
		\draw [style=green line] (7.center) to (8.center);
		\draw [style=thick blue line] (12.center) to (11.center);
		\draw [style=blue line] (10.center) to (9.center);
		\draw [style=black arrow] (4) to (14.center);
		\draw [style=black arrow] (4) to (15.center);
		\draw [style=black arrow] (4) to (16.center);

	\end{pgfonlayer}


    \clip       (-\cols, -\height) rectangle (\cols, \height);
    \draw[gray] (-\cols, -\height) rectangle (\cols, \height);

    \pgfmathsetmacro{\from}{-2 *\cols}
    \pgfmathsetmacro{\to}{2 * \cols}
    \foreach\i in {\from, ..., \to} {
        \draw[xslant=\slant]  (\i, -\height) -- (\i, \height);
        \draw[xslant=-\slant] (\i, -\height) -- (\i, \height);
    }

    \foreach\j in {-\rows, ..., \rows} {
        \pgfmathsetmacro{\y}{0.5 * \j * tan(60)}
        \draw (-\cols, \y) -- (\cols, \y);
    }

\end{tikzpicture}

\caption{Alcove geometry for $\tilde A_{2}$}  \label{fig:overview}
\end{figure}

Figure \ref{fig:overview} shows the alcove and chamber geometry. (In this and subsequent figures, the origin is denoted by a large red dot.) For the sake of normalizing figures and geometric terminology, we assume that the simple roots $\ga_{1}$ and $\ga_{2}$ make angles of $-\pi/6$ and $\pi/2$ with the positive $x$-axis, respectively. Then the root hyperplanes (which are lines) and root strips are at angles which are integral multiples of $\pi/3$. The alcoves are equilateral triangles with one edge horizontal and the opposite vertex either above or below; we call these ``Up'' and ``Down,'' respectively. (In \cite{GrLi:21} they were referred to as ``Even'' and ``Odd.'') We number the chambers I, II, \dots, VI
beginning with the chamber in the first quadrant and proceeding counterclockwise, as shown in Figure \ref{fig:overview}.
The alcove centers on a fundamental strip form a line through $q$ (the center of $A_{\circ}$); we will call the portion of the strip whose alcove centers lie on one of the two rays on this line (starting at $q$) a fundamental half-strip.  The chambers numbered I, III, and V (resp.~ II, IV, VI) are called odd (resp.~even) chambers.

Assume $w\ne e$. Then there is at least one simple reflection $s$ (resp.\ $t$) such that $w<ws$ (resp.\ $w>wt$). We call $w$ of Type $\tau$ ($\tau=1, 2$) if the number of simple reflections $s$ such that $w<ws$ is exactly $\tau$. Evidently when $w$ is of Type $\tau$ there are exactly $3-\tau$ simple reflections $t$ such that $w>wt$. Note also that when $w<ws$ (for $s$ simple), if $w$ is of Type 1 then $ws$ is of Type 2; the converse is true if and only if $w$ is not a spiral element.  This is the reason some of the inductive procedures of this paper do not work for spiral elements.

Recall that for $x,y$ in $W$, $H(x,y)$ denotes the set of root hyperplanes separating the alcoves $xA_0$ and $yA_0$.
The sets $H(x,e)$ and $H(y,e)$ can differ only by the hyperplanes in $H(x,y)$.  Precisely, we have:

\begin{Lem} \label{l:Hyperplanesep}
Let $x,y \in W$.  If $E$ is a root hyperplane not in $H(x,y)$, then either $E$ is in both of the sets $H(x,e)$ and $H(y,e)$, or in neither of them.
\end{Lem}

\begin{proof}
If $E$ is not in $H(x,y)$, then $x A_0$ and $y A_0$ lie on the same side of $E$.  If $A_0$ lies on this side of $E$, then $E$ is in neither of $H(x,e)$ and $H(y,e)$;
if $A_0$ is on the other side of $E$, then $E$ is in both of them.
\end{proof}

The following useful lemma describes how the length of an element changes when we reflect across the edge of a certain root strip.  

\begin{Lem} \label{l:lengthChange}
Let $E$ be an edge of a fundamental root strip $\cs$. There are three chambers for which this edge is the closer of the two edges of $\cs$ to the chamber.
Let $\cc$ be one of these chambers, let $z\in W$ be in $\cc$, and let $z'\in W$ be obtained by reflecting $z$ across $E$. Then $\ell(z') < \ell(z)$, so
$z'<z$.  More precisely, if $E \cap A_0$ is the vertex of $\cc$, then $\ell(z') = \ell(z) - 3$.  In all other cases, $\ell(z') = \ell(z) - 1$; in particular, this is always the case when a ray along $E$ is a boundary of $\cc$. 
\end{Lem}

\begin{proof}
The fact that there are three chambers for which $E$ is the closer of the two edges of $\cs$ to the chamber can be seen by inspection: they are the three chambers on the opposite side of $E$ from $\cs$.

Let $r$ be reflection across $E$, so $z' = rz$.  We have $H(rz,r) = r H(z,e)$.  Moreover, $\ell(z) = |H(z,e)|$ and $\ell(z') = |H(rz,e)|$. 
Our hypothesis implies that
$E$ is in $H(z,e)$, since $e$ belongs to $\cs$.
Since $z A_0$ and $A_0$ are on opposite sides of $E$, $E$ is not in $H(rz,e)$.

First suppose $E \cap A_0$ is an edge of $A_0$.  Then $E$ is the only root hyperplane separating $A_0$ and $rA_0$ (that is, in $H(r,e)$).
By the preceding paragraph and Lemma \ref{l:Hyperplanesep},
$H(rz,r) = H(rz,e) \sqcup \{E \}$.
Hence
$$
\ell(z') = |H(rz,e)| = |H(rz,r)| -1 = |H(z,e)| - 1 = \ell(z) - 1,
$$
proving the lemma in this case.

Next suppose that $E \cap A_0$ is a vertex $p$ of $A_0$.  There are $3$ hyperplanes $E, H, H'$ between $A_0$ and $rA_0$.  We have $rE = E$ and
$rH = H'$.  If $p$ is the vertex of the chamber $\cc$ containing $z$, then each of $E,H$, and $H'$ is in $H(z,e)$, and $H(rz,r) = H(rz,e) \sqcup \{E,H,H' \}$.
Reasoning as above, we see that $\ell(z') = \ell(z) - 3$.  If $p$ is not the vertex of $\cc$, then one of the hyperplanes $H,H'$ is in $H(z,e)$, and the other is
not.  By relabelling, we may assume  $H$ is in $H(z,e)$.  Then
$$
H(rz,r) = (H(rz,e) \sqcup \{H',E\}) \smallsetminus \{ H \}.
$$
Hence $\ell(z') = \ell(z) - 1$.
\end{proof}

\begin{Rem}
In the setting of Lemma \ref{l:lengthChange},
there is at most one chamber $\cc$ such that $\ell(z') = \ell(z) - 3$.  Such a chamber exists exactly when $E$ intersects the fundamental alcove $A_0$ in a vertex $p$
of $A_0$; then $p$ is also the vertex of $\cc$.
\end{Rem}

One application of Lemma \ref{l:lengthChange} is the following proposition, which relates the chambers to the groups $L(w)$.

\begin{Prop} \label{p:evenoddchamber}
Let $\cc$ be a chamber, and let $w \in \cc$.

\noindent (a)  If $\cc$ is an even chamber, then the reflection across each wall of the chamber is a simple reflection (on the left).  
If $s$ is either of those simple reflections, then $s w < w$.  Thus, $L(w)$ is a Coxeter group of type $A_2$ 
generated by these two simple reflections.

\noindent (b)  If $\cc$ is an odd chamber, then neither of the reflections across a wall of this chamber is simple (on the left), so
$L(w)$ is the group with two elements.
 \end{Prop}
 
 \begin{proof}
 Lemma \ref{l:lengthChange} implies that if $r$ is the reflection across a wall of $\cc$, then $r w<w$.
One can verify by inspection that if $\cc$ is even (resp.~odd) then both (resp.~neither) of the reflections across a wall of the chamber 
is simple. (The simple reflections on the left correspond to the thick green, red, and blue lines in Figure \ref{fig:overview}.) Thus, if $\cc$ is even, then $L(w)$ contains both of these simple reflections, and the remaining assertions about $L(w)$ follow
(see Section \ref{ss.affineWeyl}).  If $\cc$ is odd, and if $s$ is a reflection across a wall of $\cc$, then this wall is the nearer edge of 
a root strip to $\cc$.  The reflection $t$ across the other edge of this
root strip is simple.  Applying Lemma
\ref{l:lengthChange} to the element $tw$, we see that $t(tw)< tw$, so $w < tw$.  Hence $t$ is not in $L(w)$.  This excludes
two of the simple reflections from $L(w)$, so $L(w)$ contains only one simple reflection (across a hyperplane not parallel to a wall of $\cc$).
Hence $L(w)$ is a Weyl group of type $A_1$ with two elements.
 \end{proof} 
 
\begin{Lem} \label{l:rootchamber}
Let $\cc$ be a chamber.  There is a unique $\ga \in \Phi$ with the following property:
If $z \in \cc$, then $t(\ga) z = z + \ga \in \cc$.  In this case, we say $\ga$ points into the chamber $\cc$.
\end{Lem}

\begin{proof}
The root $\ga$ is the one parallel to the angle bisector of the chamber, in the direction from the chamber vertex into the chamber. See Figure \ref{fig:overview}. Specifically, for chambers I, II, III, $\ga = \tilde\ga,\ \ga_{2}, -\ga_{1}$, respectively, and for chambers IV, V, VI, the negatives of those. 
\end{proof}

There is a bijection between the set of chambers and $\Phi$, which takes a chamber
$\cc$ to the unique root $\ga \in \Phi$ which points into $\cc$.  Given a non-spiral element $w$,
there is an operation of ``translation into the chamber" which takes $w$ to $w' = t(\alpha) w$, where $\alpha$ is the root pointing into the
chamber containing $w$.  In Section \ref{s:translation}, we will make use of the following result.

\begin{Prop} \label{p:lengthChange}
Suppose $w$ is a non-spiral element of $W$.  Let $w'$ be the element obtained by translating $w$ into the chamber.
Then $w' > w$, and $\ell(w') = \ell(w) + 4$.  
\end{Prop}
 
To prove this proposition, we use the following lemma. 

\begin{Lem} \label{l:lengthChange2}
Suppose $w \in W$ is a non-spiral element.  Let $\ga$ be the root pointing into the chamber $\cc$ containing $w$, and suppose the alcove $w A_0$ is between the hyperplanes $H_{\ga,k-1}$ and $H_{\ga,k}$.  Then $k \geq 1$.  
If $w A_0 \cap H_{\ga,k}$ is an edge of $wA_0$, then $\ell(s_{\ga,k} w) = \ell(w) + 1$.  If $w A_0 \cap H_{\ga,k}$ is a vertex of $wA_0$, then
$\ell(s_{\ga,k} w) = \ell(w) + 3$.
\end{Lem}

\begin{proof}
The inequality $k \geq 1$ holds because for any $p \in \cc$, we have $(\ga, p) \geq 0$, as can be seen by inspection.
Let $w' = s_{\ga,k} w$.  If $w A_0 \cap H_{\ga,k}$ is an edge of $wA_0$, then $H(w,w') = \{ H_{\ga,k} \}$, and $H(w',e) = H(w,e) \sqcup \{ H_{\ga,k} \}$, so
$\ell(s_{\ga,k} w) = \ell(w) + 1$.  If $w A_0 \cap H_{\ga,k}$ is a vertex of $wA_0$, then the edges of $w A_0$ are the intersections of $w A_0$ with
$H_{\ga, k+1}$ and two other hyperplanes $H, H'$.  The hyperplanes $H$ and $H'$ are parallel to the edges of $\cc$, and
$H(w,w') = \{ H_{\ga,k}, H, H' \}$.  Moreover, none of the hyperplanes in $H(w,w')$ are
in $H(w,e)$.
Since $H(w',e) = H(w,e) \sqcup  \{ H_{\ga,k}, H, H' \}$, we have
$\ell(s_{\ga,k} w) = \ell(w) + 3$.
\end{proof}

 {\em Proof of Proposition \ref{p:lengthChange}.}
 Let $\ga$ be the root pointing into the chamber of $w$.  Then $w A_0$ is between $H_{\ga,k-1}$ and $H_{\ga,k}$ for some $k \geq 1$.
Let $u = s_{\ga,k} w$; then $u A_0$ is between $H_{\ga,k}$ and $H_{\ga,k+1}$.  Moreover, $w' = t(\ga^{\vee}) w = s_{\ga,k+1} s_{\ga,k} w =
s_{\ga,k+1} u$.  By Lemma \ref{l:lengthChange2}, either $\ell(u) = \ell(w)+1$ and $\ell(w') = \ell(u) + 3$, or $\ell(u) = \ell(w)+3$ and $\ell(w') = \ell(u) + 1$.  
In either case, $\ell(w') = \ell(w) + 4$, and by the definition of the Bruhat order, $w' > w$. \hfill $\Box$

\subsection{Alcove edge labels} 

\begin{figure}[htbp!]
%

\pgfmathsetmacro{\cols}{5}
\pgfmathsetmacro{\rows}{6}
\pgfmathsetmacro{\slant}{cot(60)}
\pgfmathsetmacro{\height}{0.5 * \rows * tan(60)}
\pgfmathsetmacro{\triht}{sin(60)}
\pgfmathsetmacro{\upmid}{0.25 * sec(30)}
\pgfmathsetmacro{\downmid}{\triht - \upmid}

\begin{tikzpicture}[every label/.append style = {fill=white}]
    

\input{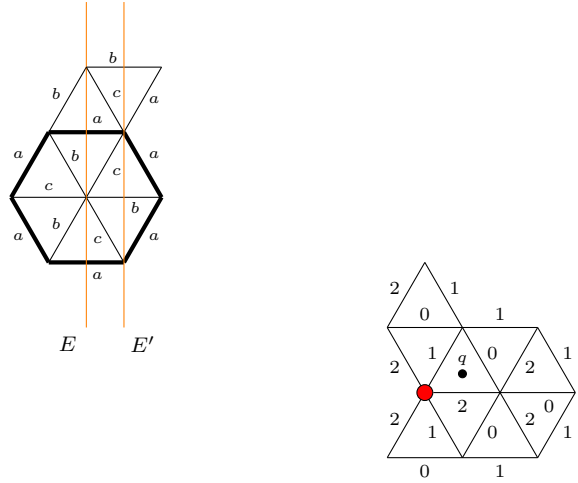}

	\begin{pgfonlayer}{nodelayer}
		\node [style=red dot] (4) at (0, 0) {};
		\node [style=small black dot] (13) at (0.5, 0.25) {};
		\node [style=none] (14) at (1.5, -1 * \triht) {};
		\node [style=none] (15) at (0, 2 * \triht) {};
		\node [style=none] (16) at (1.5, 1 * \triht) {};
		\node [style=none] (17) at (2, 0) {};
		\node [style=none] (18) at (0.5, -1 * \triht) {};
		\node [style=none] (19) at (-0.5, -1 * \triht) {};
		\node [style=none] (20) at (-0.5, 1 * \triht) {};
		\node [style=none] (21) at (0.5, 1 * \triht) {};
		\node [style=none] (22) at (-4, 2 * \triht) {};
		\node [style=none] (23) at (-5, 2 * \triht) {};
		\node [style=none] (24) at (-4.5, 3 * \triht) {};
		\node [style=none] (25) at (-5, 4 * \triht) {};
		\node [style=none] (26) at (-4, 4 * \triht) {};
		\node [style=none] (27) at (-4.5, 5 * \triht) {};
		\node [style=none] (28) at (-3.5, 3 * \triht) {};
		\node [style=none] (29) at (0.9, 0.6 * \triht) {\tiny 0};
		\node [style=none] (31) at (0.1, 0.6 * \triht) {\tiny 1};
		\node [style=none] (33) at (0.5, -0.2 * \triht) {\tiny 2};
		\node [style=none] (35) at (1, 1.2 * \triht) {\tiny 1};
		\node [style=none] (37) at (1.4, 0.4 * \triht) {\tiny 2};
		\node [style=none] (39) at (-0.4, 0.4 * \triht) {\tiny 2};
		\node [style=none] (41) at (0, 1.2 * \triht) {\tiny 0};
		\node [style=none] (43) at (0.4, 1.6 * \triht) {\tiny 1};
		\node [style=none] (45) at (-0.4, 1.6 * \triht) {\tiny 2};
		\node [style=none] (47) at (0.1, -0.6 * \triht) {\tiny 1};
		\node [style=none] (51) at (0.9, -0.6 * \triht) {\tiny 0};
		\node [style=none] (53) at (-0.4, -0.4 * \triht) {\tiny 2};
		\node [style=none] (55) at (0, -1.2 * \triht) {\tiny 0};
		\node [style=none] (56) at (1, -1.2 * \triht) {\tiny 1};
		\node [style=none] (57) at (1.4, -0.4 * \triht) {\tiny 2};
		\node [style=none] (58) at (1.65, -0.2 * \triht) {\tiny 0};
		\node [style=none] (59) at (1.9, -0.6 * \triht) {\tiny 1};
		\node [style=none] (60) at (1.9, 0.6 * \triht) {\tiny 1};
		\node [style=none] (61) at (-4.5, 6 * \triht) {};
		\node [style=none] (62) at (-4.5, 1 * \triht) {};
		\node [style=none] (63) at (-4.65, 3.65 * \triht) {$\scriptscriptstyle b$};
		\node [style=none] (64) at (-4.35, 4.2 * \triht) {$\scriptscriptstyle a$};
		\node [style=none] (65) at (-4.1, 3.4 * \triht) {$\scriptscriptstyle c$};
		\node [style=none] (66) at (-4.9, 4.6 * \triht) {$\scriptscriptstyle b$};
		\node [style=none] (67) at (-4.1, 4.6 * \triht) {$\scriptscriptstyle c$};
		\node [style=none] (68) at (-3.6, 3.65 * \triht) {$\scriptscriptstyle a$};
		\node [style=none] (69) at (-3.85, 2.85 * \triht) {$\scriptscriptstyle b$};
		\node [style=none] (70) at (-3.6, 2.4 * \triht) {$\scriptscriptstyle a$};
		\node [style=none] (71) at (-4.35, 2.35 * \triht) {$\scriptscriptstyle c$};
		\node [style=none] (72) at (-4.35, 1.8 * \triht) {$\scriptscriptstyle a$};
		\node [style=none] (73) at (-4.9, 2.6 * \triht) {$\scriptscriptstyle b$};
		\node [style=none] (74) at (0.5, 0.5 * \triht) {$\scriptscriptstyle q$};
		\node [style=none] (75) at (-5.5, 3 * \triht) {};
		\node [style=none] (76) at (-5.4, 3.65 * \triht) {$\scriptscriptstyle a$};
		\node [style=none] (77) at (-5.4, 2.4 * \triht) {$\scriptscriptstyle a$};
		\node [style=none] (78) at (-5, 3.17 * \triht) {$\scriptscriptstyle c$};
		\node [style=none] (79) at (-3.5, 5 * \triht) {};
		\node [style=none] (80) at (-4, 6 * \triht) {};
		\node [style=none] (81) at (-4, 1 * \triht) {};
		\node [style=none] (82) at (-4.15, 5.15 * \triht) {$\scriptscriptstyle b$};
		\node [style=none] (83) at (-3.6, 4.5 * \triht) {$\scriptscriptstyle a$};
		\node [style=none] (84) at (-4.75, 0.75 * \triht) {$\scriptstyle E$};
		\node [style=none] (85) at (-3.75, 0.75 * \triht) {$\scriptstyle E'$};
	\end{pgfonlayer}
	\begin{pgfonlayer}{edgelayer}
		\draw (26.center) to (23.center);
		\draw (25.center) to (22.center);
		\draw [style=black line] (28.center) to (22.center);
		\draw [style=black line] (25.center) to (26.center);
		\draw (24.center) to (28.center);
		\draw [style=black line] (23.center) to (22.center);
		\draw (27.center) to (25.center);
		\draw (15.center) to (14.center);
		\draw (20.center) to (18.center);
		\draw (21.center) to (19.center);
		\draw (16.center) to (18.center);
		\draw (17.center) to (14.center);
		\draw (15.center) to (20.center);
		\draw (16.center) to (17.center);
		\draw (20.center) to (16.center);
		\draw (4) to (17.center);
		\draw (19.center) to (14.center);
		\draw [style=orange line] (61.center) to (62.center);
		\draw [style=black line] (25.center) to (75.center);
		\draw [style=black line] (75.center) to (23.center);
		\draw (24.center) to (75.center);
		\draw (27.center) to (26.center);
		\draw [style=black line] (26.center) to (28.center);
		\draw (27.center) to (79.center);
		\draw (79.center) to (26.center);
		\draw [style=orange line] (80.center) to (81.center);
	\end{pgfonlayer}


%
%

\end{tikzpicture}

\caption{Labels on alcove walls}  \label{fig:wallLabels}
\end{figure}

In this subsection we briefly describe the \emph{right} action of simple reflections on alcoves in terms of a certain labeling of the alcove walls (see Figure \ref{fig:wallLabels}). Begin by labeling the three walls of $A_{\circ}$ by 1, 2, 0 along $H_{\ga_{1},0}$, $H_{\ga_{2},0}$, and $H_{\tilde\ga,1}$, respectively. Recursively, having labeled the walls of an alcove $A$ by 1, 2, 0, and given an alcove $A'$ which is the reflection of $A$ across one of its walls, label the walls of $A'$ by reflecting the labels on the walls of $A$. See the lower right portion of Figure \ref{fig:wallLabels}, which shows the labels on the walls of several alcoves near $A_{\circ}$ (recall its center is denoted $q$). The fact that the labels are consistent when moving around a small hexagon implies that these alcove wall labels are well defined. 

The upper left portion of Figure \ref{fig:wallLabels} shows the behavior of labels on the walls of some alcoves lying on or near two parallel root strings (depicted in orange). Note that in the small hexagon outlined in bold, the exterior edges all have the same label, $a$, and the interior edges alternate between the other two labels, $b$ and $c$. This is a general phenomenon. The proof of the following lemma is straightforward and is left to the reader; see Figure \ref{fig:wallLabels}.

\begin{Lem} \label{l:wallLabels}  (a) Let $A, A'$ be two alcoves sharing a wall with label $a\in\{0,1,2\}$ and corresponding to elements $w, w'\in W$, respectively. Then $w'=w s_{a}$. \hspace{2em}

\noindent (b) Let $A, A'$ be two alcoves such that $A'$ is the reflection of $A$ across some root hyperplane. Then the labels on $A'$ are the reflections of the labels on $A$. \hspace{2em}

\noindent (c) Given six alcoves sharing a common vertex and forming a small hexagon, the labels on the six exterior edges of the hexagon are all the same, say $a$. The labels on the edges incident to the center of the hexagon alternate between $b$ and $c$, where $\{a,b,c\} = \{0,1,2\}$. \hspace{2em}

\noindent (d) Given two alcoves lying on a fixed root string, their edges on the same side of the string have the same label. \label{l:wallLabels(d)} \hfill $\Box$
\end{Lem}

\subsection{Paths and spiral factorizations} 
It will be convenient to have a geometric criterion for when a product of simple reflections is reduced. Given an expression $w = t_{1} t_{2 }\dots t_{m}$, where each $t_{i}$ is a simple reflection, we can define a sequence of alcoves $A_{0} =A_{\circ},\ A_{1}= t_{1} A_{\circ},\ A_{2}= t_{1} t_{2} A_{\circ},\ \dots,\ A_{m}=w A_{\circ}$, in which each alcove $A_{i}$ is obtained from the previous one $A_{i-1}$ by reflecting it across the wall of $A_{i-1}$ labeled by $t_{i}$.
If we shade in all the alcoves $A_{0}, \dots, A_{m}$, we obtain a connected, piecewise-linear, directed path $\cp$ (whose width is the distance between two adjacent root hyperplanes), from $A_{\circ}$ to $w A_{\circ}$. We say $\cp$ is a reduced path if the expression for $w$ is reduced.  Such paths are related to the galleries defined in 
\cite{GaLi:05}.

\begin{figure}[htbp!]
%

\pgfmathsetmacro{\cols}{11}
\pgfmathsetmacro{\rows}{9}
\pgfmathsetmacro{\slant}{cot(60)}
\pgfmathsetmacro{\ymultiplier}{0.5 * tan(60)}
\pgfmathsetmacro{\height}{\rows * \ymultiplier}
\pgfmathsetmacro{\triht}{sin(60)}
\pgfmathsetmacro{\upmid}{0.25 * sec(30)}
\pgfmathsetmacro{\downmid}{\triht - \upmid}

\begin{tikzpicture}[every label/.append style = {fill=white}]
    

\input{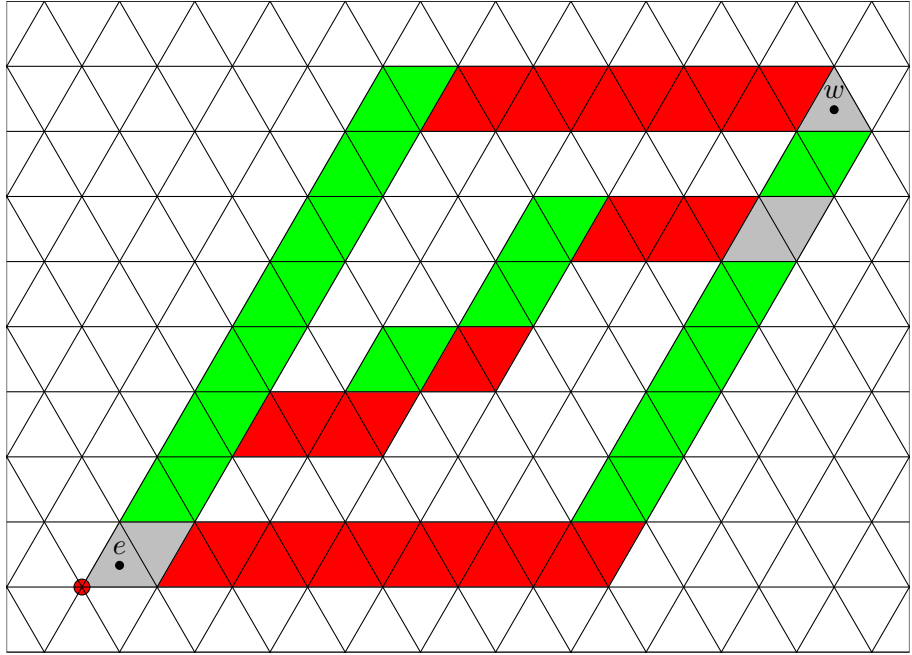}

	\begin{pgfonlayer}{nodelayer}
		\node [style=red dot] (0) at (0, 0) {};
		\node [style=none] (1) at (0.5, 1 * \triht) {};
		\node [style=none] (2) at (1.5, 1 * \triht) {};
		\node [style=none] (3) at (1, 0) {};
		\node [style=none] (4) at (7, 0) {};
		\node [style=none] (5) at (7.5, 1 * \triht) {};
		\node [style=none] (6) at (6.5, 1 * \triht) {};
		\node [style=none] (7) at (8.5, 5 * \triht) {};
		\node [style=none] (8) at (9.5, 5 * \triht) {};
		\node [style=none] (9) at (9, 6 * \triht) {};
		\node [style=none] (10) at (10, 6 * \triht) {};
		\node [style=none] (11) at (9.5, 7 * \triht) {};
		\node [style=none] (12) at (10.5, 7 * \triht) {};
		\node [style=none] (13) at (10, 8 * \triht) {};
		\node [style=none] (14) at (4, 8 * \triht) {};
		\node [style=none] (15) at (5, 8 * \triht) {};
		\node [style=none] (16) at (4.5, 7 * \triht) {};
		\node [style=none] (17) at (2, 2 * \triht) {};
		\node [style=none] (19) at (4, 2 * \triht) {};
		\node [style=none] (20) at (4.5, 3 * \triht) {};
		\node [style=none] (23) at (5, 4 * \triht) {};
		\node [style=none] (24) at (5.5, 3 * \triht) {};
		\node [style=none] (25) at (6, 4 * \triht) {};
		\node [style=none] (26) at (6, 6 * \triht) {};
		\node [style=none] (27) at (7, 6 * \triht) {};
		\node [style=none] (28) at (6.5, 5 * \triht) {};
		\node [style=none] (29) at (2.5, 3 * \triht) {};
		\node [style=none] (30) at (4, 4 * \triht) {};
		\node [style=none] (31) at (3.5, 3 * \triht) {};
		\node [style=small black dot] (32) at (10, 7 * \triht + \upmid) {};
		\node [style=none] (33) at (10, 7 * \triht + \downmid - 0.05) {$w$};
		\node [style=small black dot] (34) at (0.5, \upmid) {};
		\node [style=none] (34) at (0.5, \downmid - 0.05) {$e$};		
	\end{pgfonlayer}
	\begin{pgfonlayer}{edgelayer}
		\draw [style=filled grey] (1.center)
			 to (2.center)
			 to (3.center)
			 to (0.center)
			 to cycle;
		\draw [style=filled red] (3.center)
			 to (2.center)
			 to (5.center)
			 to (4.center)
			 to cycle;
		\draw [style=filled green] (14.center)
			 to (15.center)
			 to (2.center)
			 to (1.center)
			 to cycle;
		\draw [style=filled red] (16.center)
			 to (15.center)
			 to (13.center)
			 to (11.center)
			 to cycle;
		\draw [style=filled grey] (12.center)
			 to (11.center)
			 to (13.center)
			 to cycle;
		\draw [style=filled green] (9.center)
			 to (10.center)
			 to (12.center)
			 to (11.center)
			 to cycle;
		\draw [style=filled grey] (9.center)
			 to (7.center)
			 to (8.center)
			 to (10.center)
			 to cycle;
		\draw [style=filled green] (5.center)
			 to (6.center)
			 to (7.center)
			 to (8.center)
			 to cycle;
		\draw [style=filled green] (20.center) to (23.center);
		\draw [style=filled green] (27.center)
			 to (26.center)
			 to (23.center)
			 to (25.center)
			 to cycle;
		\draw [style=filled red] (25.center)
			 to (24.center)
			 to (20.center)
			 to (23.center)
			 to cycle;
		\draw [style=filled red] (27.center)
			 to (9.center)
			 to (7.center)
			 to (28.center)
			 to cycle;
		\draw [style=filled red] (29.center)
			 to (17.center)
			 to (19.center)
			 to (20.center)
			 to cycle;
		\draw [style=filled green] (31.center)
			 to (30.center)
			 to (23.center)
			 to (20.center)
			 to cycle;
	\end{pgfonlayer}


    \clip       (-1, -\ymultiplier) rectangle (\cols, \height);
    \draw[gray] (-1, -\ymultiplier) rectangle (\cols, \height);

    \pgfmathsetmacro{\from}{-2 *\cols}
    \pgfmathsetmacro{\to}{2 * \cols}
    \foreach\i in {\from, ..., \to} {
        \draw[xslant=\slant]  (\i, -\height) -- (\i, \height);
        \draw[xslant=-\slant] (\i, -\height) -- (\i, \height);
    }

    \foreach\j in {-1, ..., \rows} {
        \pgfmathsetmacro{\y}{0.5 * \j * tan(60)}
        \draw (-1, \y) -- (\cols, \y);
    }

\end{tikzpicture}

\caption{Three reduced paths from $e$ to $w$}  \label{fig:paths}
\end{figure}

Figure \ref{fig:paths} shows three (reduced) paths from $A_{\circ}$ to $w A_{\circ}$ for a certain $w$ in Chamber~I. Path segments parallel to $\ga_{1}$ (resp.\ $\tilde\ga$) are colored red (resp.\ green), while segments which would be red in one path and green in another are grey.  In coloring this path, we follow the convention implied in
Remark \ref{r:spiralFactorizations} below: in a path of the form $uv$, where $u$ and $v$ are spiral, we make $u$ as long as possible, and $v$ moves out of the root strip containing $u$. In a zigzag path, we apply this convention iteratively from the left, so the color changes only when the path changes direction.

A chamber is bounded by two fundamental root strip ``rays,'' each of which is a union of alcoves beginning at $A_{\circ}$ and extending to infinity.  In contrast to our usual convention, for the purposes of the following lemma we will consider the boundary strips to be part of the chamber.

\begin{Prop} \label{p:reducedPath}
Suppose $w=t_{1} \dots t_{m}$ as above with associated path $\cp$, and with $w A_{\circ}$ belonging to a chamber $\cc$. If the expression is reduced, then every directed linear segment of $\cp$ is in the direction of one of the two boundary strip rays of $\cc$.
\end{Prop}

\begin{proof}
If $w$ is spiral, then as noted in \cite{GrLi:15}, $w$ has a unique reduced expression; the associated path is the portion of the fundamental (spiral) strip from $A_{\circ}$ to $w A_{\circ}$. (If this were not true, then the path would leave the strip at some point, and at some later point it would have to return, say from $A_{i-1}$ to $A_{i}$. But the fact that length equals number of hyperplanes separating the associated alcove from $A_{\circ}$ implies that $\ell(t_{1}\dots t_{i}) = \ell(t_{1}\dots t_{i-1})-1$, contradicting the fact that the expression is reduced.)

For $w$ non-spiral, the proof is by induction on $\ell(w)$. Set $w'=w t_{m}$, where $\ell(w')=\ell(w)-1$. It is easy to see that $w'$ also lies in $\cc$ (because we expanded $\cc$ to include both its boundary strips). If $w'$ is of Type 2 then the path from $w'$ to $w$ is in one of the two $\cc$ edge strip directions from $w'A_{\circ}$ to $wA_{\circ}$, whereas if $w'$ is of Type 1 then the path from $w'$ to $w$ can be viewed as lying in both of the $\cc$ edge strip directions from $w'A_{\circ}$ to $wA_{\circ}$. Since by induction the linear segments of the path from $A_{\circ}$ to $w'A_{\circ}$ all lie in one of the $\cc$ edge directions, the same is true for the entire path from $A_{\circ}$ to $w A_{\circ}$.
\end{proof}

\begin{Cor} \label{c:reducedParallelogram}
For $w$ non-spiral, the path $\cp$ associated to any reduced expression for $w$ lies entirely within the parallelogram with vertices $A_{\circ}$ and $w A_{\circ}$ and edges parallel to the fundamental root strips bounding the chamber of $w$.
\end{Cor}

\begin{Rem} \label{r:spiralFactorizations}
Each non-spiral $w$ thus has two canonical reduced expressions, corresponding to the two paths from $A_{\circ}$ to $w A_{\circ}$ along the boundary of the parallelogram described in Corollary \ref{c:reducedParallelogram}. These canonical paths are illustrated in Figure \ref{fig:paths}. Since the sequence of alcove edge labels crossed along any segment of a root strip has the pattern $a b c a b c \dots$, we have two canonical factorizations $w=uv$ where $u$ and $v$ are spiral and $\ell(u)+\ell(v)=\ell(w)$. We will always assume that $u$ is as long as possible in its fundamental root strip, so that $u v_{1}$ lies outside that strip, where $v_{1}$ is the first simple reflection in the (unique!) reduced expression for $v$.
\end{Rem}

\begin{Lem} \label{l:edgeTranslation}
Let $E$ and $E'$ be two parallel root strings. Then there is a fixed (spiral) element $v\in W$ such that, for any $x$ on $E$, the ``corresponding'' alcove $x'$ on $E'$ (meaning that $x$ and $x'$ lie in the same root strip perpendicular to $E$) is given by $x'=xv$.
\end{Lem}

\begin{proof}
Since the sequence of edges crossed when moving along any root strip always defines a spiral element of $W$, it suffices to show that, for any two alcoves $x$ and $y$ on $E$, the first edges crossed when moving from $x$ or $y$ along root strips perpendicular to $E$ and towards $E'$ have the same label.  But this follows from Lemma \ref{l:wallLabels}(d). (See Figure \ref{fig:wallLabels}, where the edges joining $E$ to the adjacent root string $E'$ all have the same label, $c$.)
\end{proof}


\section{Bruhat hexagons} \label{s.hexagons}
In this section, we will show that for any (non-spiral) $w\in W$, the set of alcove centers $\{ xq \mid x \le w \}$ is the intersection of the lattice $Wq$ with a certain closed convex hexagon in the plane, which we will describe explicitly in two different ways.  Although these hexagons are not generally regular, they have more symmetry for even-chamber $w$ than for odd-chamber $w$.

\begin{figure}[htbp!]
%

\pgfmathsetmacro{\cols}{5}
\pgfmathsetmacro{\rows}{5}
\pgfmathsetmacro{\slant}{cot(60)}
\pgfmathsetmacro{\height}{0.5 * \rows * tan(60)}
\pgfmathsetmacro{\triht}{sin(60)}
\pgfmathsetmacro{\upmid}{0.25 * sec(30)}
\pgfmathsetmacro{\downmid}{\triht - \upmid}

\begin{tikzpicture}
    

\input{sample.tikzstyles}

	\begin{pgfonlayer}{nodelayer}
		\node [style=red dot] (4) at (0, 0) {};
		\node [style=small black dot, label={right:$w_3$}] (23) at (3, 2 * \triht + \downmid) {};
		\node [style=small black dot, label={above:$w_4$}] (24) at (0.5, 4 * \triht + \upmid) {};
		\node [style=small black dot, label={above left:$w_5$}] (25) at (-3.5, 1 * \triht + \downmid) {};
		\node [style=small black dot, label={below left:$w=w_0$}] (26) at (-3.5, -2 * \triht + \upmid) {};
		\node [style=small black dot, label={below:$w_1$}] (27) at (0.5, -5 * \triht + \downmid) {};
		\node [style=small black dot, label={below right:$w_2$}] (28) at (3, -3 * \triht + \upmid) {};
		\node [style=none, label={left:$H_{\beta,j}$}] (29) at (-5, 0) {};
		\node [style=none] (30) at (5, 0) {};
		\node [style=none] (31) at (2.5, 5 * \triht) {};
		\node [style=none, label={below:$H_{\alpha,i}$}] (32) at (-2.5, -5 * \triht) {};
		\node [style=none] (33) at (-2.5, 5 * \triht) {};
		\node [style=none, label={below:$H_{\gamma,k}$}] (34) at (2.5, -5 * \triht) {};
	\end{pgfonlayer}
	\begin{pgfonlayer}{edgelayer}
		\draw (23) to (24);
		\draw (24) to (25);
		\draw (25) to (26);
		\draw (26) to (27);
		\draw (27) to (28);
		\draw (28) to (23);
		\draw [style=red line] (29.center) to (30.center);
		\draw [style=green line] (32.center) to (31.center);
		\draw [style=blue line] (33.center) to (34.center);
		\draw [style=dashed line] (26) to (23);
		\draw [style=dashed line] (24) to (27);
		\draw [style=dashed line] (25) to (28);
	\end{pgfonlayer}

    \clip       (-\cols, -\height) rectangle (\cols, \height);
    \draw[gray] (-\cols, -\height) rectangle (\cols, \height);

    \pgfmathsetmacro{\from}{-2 *\cols}
    \pgfmathsetmacro{\to}{2 * \cols}
    \foreach\i in {\from, ..., \to} {
        \draw[xslant=\slant]  (\i, -\height) -- (\i, \height);
        \draw[xslant=-\slant] (\i, -\height) -- (\i, \height);
    }

    \foreach\j in {-\rows, ..., \rows} {
        \pgfmathsetmacro{\y}{0.5 * \j * tan(60)}
        \draw (-\cols, \y) -- (\cols, \y);
    }

\end{tikzpicture}

\caption{A hexagon $\ch_{w}$ for chamber IV}  \label{fig:hexagoneven}
\end{figure}

\begin{figure}[htbp!]
%

\pgfmathsetmacro{\cols}{6}
\pgfmathsetmacro{\rows}{6}
\pgfmathsetmacro{\slant}{cot(60)}
\pgfmathsetmacro{\height}{0.5 * \rows * tan(60)}
\pgfmathsetmacro{\triht}{sin(60)}
\pgfmathsetmacro{\upmid}{0.25 * sec(30)}
\pgfmathsetmacro{\downmid}{\triht - \upmid}

\begin{tikzpicture}
    

\input{sample.tikzstyles}

\begin{pgfonlayer}{nodelayer}
		\node [style=none] (2) at (-6, 1 * \triht) {};
		\node [style=none, label={right:$H_{\beta,j}$}] (3) at (6, 1 * \triht) {};
		\node [style=red dot] (4) at (0, 0) {};
		\node [style=none] (7) at (-2, -6 * \triht) {};
		\node [style=none, label={above:$H_{\alpha,i}$}] (8) at (4, 6 * \triht) {};
		\node [style=none] (11) at (4, -6 * \triht) {};
		\node [style=none, label={above:$H_{\gamma,k}$}] (12) at (-2, 6 * \triht) {};
		\node [style=small black dot, label={right:$w_0=w$}] (13) at (4, 3 * \triht + \upmid) {};
		\node [style=small black dot, label={above:$w_1$}] (14) at (2, 4 * \triht + \downmid) {};
		\node [style=small black dot, label={right:$w_5$}] (19) at (4, -2 * \triht+ \downmid) {};
		\node [style=small black dot, label={below:$w_4$}] (20) at (0.5, -4 * \triht + \upmid) {};
		\node [style=small black dot, label={left:$w_3$}] (21) at (-3, -2 * \triht+ \downmid) {};
		\node [style=small black dot, label={left:$w_2$}] (22) at (-3, 1 * \triht+ \upmid) {};
		\node [style=purple dot] (23) at (0.5, 3 * \triht + \downmid) {};
		\node [style=purple dot] (24) at (-1.5, 2 * \triht + \upmid) {};
		\node [style=purple dot] (25) at (2, -3 * \triht + \upmid) {};
		\node [style=purple dot] (26) at (2.5, -3 * \triht + \downmid) {};
	\end{pgfonlayer}
	\begin{pgfonlayer}{edgelayer}
		\draw [style=red line] (2.center) to (3.center);
		\draw [style=green line] (7.center) to (8.center);
		\draw [style=blue line] (12.center) to (11.center);
		\draw (13) to (14);
		\draw (19) to (13);
		\draw (20) to (19);
		\draw (21) to (20);
		\draw (14) to (22);
		\draw (22) to (21);
		\draw [style=dashed line] (13) to (21);
		\draw [style=dashed line] (20) to (23);
		\draw [style=dashed line] (19) to (24);
		\draw [style=dashed line] (14) to (25);
		\draw [style=dashed line] (22) to (26);
		\draw [style=thick purple line] (23) to (24);
		\draw [style=thick purple line] (25) to (26);
	\end{pgfonlayer}
	

    \clip       (-\cols, -\height) rectangle (\cols, \height);
    \draw[gray] (-\cols, -\height) rectangle (\cols, \height);

    \pgfmathsetmacro{\from}{-2 *\cols}
    \pgfmathsetmacro{\to}{2 * \cols}
    \foreach\i in {\from, ..., \to} {
        \draw[xslant=\slant]  (\i, -\height) -- (\i, \height);
        \draw[xslant=-\slant] (\i, -\height) -- (\i, \height);
    }

    \foreach\j in {-\rows, ..., \rows} {
        \pgfmathsetmacro{\y}{0.5 * \j * tan(60)}
        \draw (-\cols, \y) -- (\cols, \y);
    }

\end{tikzpicture}

\caption{A hexagon $\ch_{w}$ for chamber I}  \label{fig:hexagon}
\end{figure}

Assume $w$ is not a spiral element, so $w$ belongs to some chamber $\cc$. Let $H_{\ga,i}$ and $H_{\gb,j}$ be the boundaries of $\cc$, and $H_{\gg,k}$ be the nearer of the two hyperplanes $H_{\gg,0},\ H_{\gg,1}$ to $w$, where $\{\ga,\gb,\gg\} = \{\ga_{1},\ga_{2}, \tilde\ga\}$ and $i,j,k \in \{0,1\}$. For definiteness, assume that $H_{\ga,i}$ is counterclockwise from $w$, so $H_{\gb,j}$ is clockwise from $w$. We define a hexagon $\ch_{w}$ with vertices $w_{0}=w, w_{1}, \dots, w_{5}$ labeled counterclockwise from $w$, as follows (see Figures \ref{fig:hexagoneven}, \ref{fig:hexagon}):
\begin{equation} \label{e:hexagonVertices}
w_{1}=s_{\ga,i}\,w, \quad w_{5}=s_{\gb,j}\,w, \quad w_{3}=s_{\gg,k}\,w, \quad w_{2}=s_{\gg,k}\,w_{1}, \quad w_{4}=s_{\gg,k}\,w_{5}.
\end{equation}

In this paper we will encounter many hexagons like these, with edges parallel to the three positive roots. We will refer to their edges, starting at the right vertical edge and moving counterclockwise, by the compass directions E, NE, NW, W, SW, and SE.  For simplicity, we will often say that an alcove lies on an edge if its center does.

\begin{Thm} \label{t.hexagon}
Assume $w$ is not a spiral element. The set $\{\, x\in W \mid x\le w \,\}$ equals the set of $x$ whose alcove centers $xq$ lie on or inside $\ch_{w}$ (``points in the hexagon'').
\end{Thm}

\begin{Rem}
For spiral elements, the set of $x \leq w$ is again the set of $x$ whose alcove centers lie in a convex region which we again denote by $\ch_{w}$, although it is not a hexagon.
More precisely, for spiral elements of length at least $2$,
the hexagon degenerates to a quadrilateral: two of its edges become points. For completeness, we describe one of these; the descriptions for spiral elements in the other five fundamental half-strips are analogous.
Assume $w$ lies in the fundamental $\ga_{2}$ half-strip to the left of $A_{\circ}$; these are the spiral elements considered in \cite{GrLi:21}. We will also assume $\ell(w)>1$; otherwise, $\{x\le w\} = \{1,w\}$. Then the set $\{ x\in W \mid x\le w \}$ has four vertices
(again numbered clockwise):
$$
w_{0}=w, \quad w_{1}=s_{\ga_{1}}w, \quad w_{3}=s_{\tilde\ga,0}w, \quad w_{2}=s_{\ga_{1}}w_{3}.
$$ 
All the quadrilateral vertices except (if $\ell(w) \geq 3$) $w_{3}$ lie in fundamental root strips.
\end{Rem}

The following observations about the geometry of the hexagons will be necessary in describing the rationally smooth locus of $X_{w}$.  If $v$ is a vertex of $\ch_{w}$, then two of the root strings through $v$ contain edges of the hexagon, and the third passes through the interior of the hexagon. We call the portion of this root string inside the hexagon a diagonal.  If a diagonal intersects an edge, we refer to this as the opposite edge to the vertex.  The diagonal starting at $w=w_{0}$ always ends at the opposite vertex, $w_{3}$; this is part of the next proposition.

\begin{Prop} \label{p.diagonals}
(a) If $w$ belongs to an even chamber, the vertices
of $\ch_w$ are the orbit $L(w) w$.  The hexagon is symmetric about each of the lines $H_{\ga,i}, H_{\gb,j}, H_{\gg,k}$ used to construct its vertices 
 (although it need not be a regular hexagon; see Figure \ref{fig:hexagoneven}).    Every diagonal passes through a pair of opposite vertices. \hspace{2em}
 
\noindent (b) If $w$ belongs to an odd chamber, then the hexagon is symmetric about the line $H_{\gg,k}$.
The diagonals through $w_{i}$ for $i\ne 0, 3$ do not pass through the opposite vertex, but instead intersect the opposite edge at the center of an alcove lying two alcoves away (along a root string) from a vertex alcove.  That is, the diagonals through the opposite vertices $w_{i}$ and $w_{i+3}$, $i=1,2$, are parallel and three root strings apart.
\end{Prop}

\begin{proof}
Since $w_3 = s_{\gg,k} w_0$, the vertices $w_0$ and $w_3$ lie on a diagonal parallel to $\gg$, independent of whether
the chamber is even or odd.  Suppose $w$ belongs to an even chamber $\cc$.  In the notation of \eqref{e:hexagonVertices}, 
Proposition \ref{p:evenoddchamber} implies that $s = s_{\ga, i}$ and $t = s_{\gb,j}$ are simple reflections, and they generate $L(w)$, which
is of type $A_2$.  The hyperplanes $H_{\ga,i}$, $H_{\gb,j}$ and $H_{\gg,k}$ meet at an alcove vertex (see Figure \ref{fig:overview}), so $s_{\gg,k} = sts = tst$.
Thus every vertex of $\ch_w$ is obtained from $w$ by applying an element of the group generated by $s$ and $t$.
Hence the vertices of $\ch_w$ are the orbit $L(w) w$.  Since each of the reflections $s_{\ga, i}$, $s_{\gb,j}$, and $s_{\gg,k}$ is in $L(w)$,
these reflections preserve the set $L(w) w$ of vertices, so they preserve the hexagon $\ch_w$.  Hence the hexagon is symmetric
about the lines $H_{\ga,i}, H_{\gb,j}, H_{\gg,k}$.  Finally, using the formulas from 
 \eqref{e:hexagonVertices}, we see that $w_4 = s_{\gb,j} w_1$ and $w_5 = s_{\ga,i} w_2$.
Therefore, $w_1$ and $w_4$ lie on a diagonal parallel to $\gb$, and $w_2$ and $w_5$ lie on a diagonal parallel to $\ga$.   This proves (a).

Next suppose $w$ belongs to an odd chamber.  By construction, $s_{\gg,k}$ interchanges the vertices $w_0$ and $w_3$, $w_1$ and $w_2$, and $w_4$ and $w_5$.
Hence $s_{\gg,k}$ preserves $\ch_w$, so the hexagon is symmetric about $H_{\gg,k}$.   For the remainder of the proof, assume $w$ is in chamber I; the
argument will be similar for the other odd chambers.  Our assumption implies that
 $s_{\ga,i} = s_{\ga_1,1}$, $s_{\gb,j} = s_{\ga_2,1}$, and $s_{\gg,k} = s_{\tilde{\ga},1}$.  We first consider the diagonal through $w_4$, which is a portion of an $\ga_2$ root string.
Direct calculation shows 
that $s_2 w_4 = t(-\tilde{\ga}) w_1$, which lies two alcoves away from $w_2$ on an edge parallel to $\tilde{\ga}$. 
This verifies the result for the diagonal through $w_4$; similar calculations show the result for the other diagonals.
\end{proof}

It can happen that the diagonals emanating from two adjacent vertices ($w_{1}$ and $w_{2}$, or $w_{4}$ and $w_{5}$) cross in the interior of the hexagon, and intersect the opposite edge $E$, say at the centers of alcoves $A$ and $A'$ (two alcoves in from the ends of $E$). This happens precisely when $E$ is at least 6 alcoves long.  We call the alcoves along $E$ from $A$ to $A'$ (inclusive) a \emph{special segment}. By Proposition \ref{p.diagonals}, a special segment can only exist when $w$ belongs to an odd chamber. See Figure \ref{fig:hexagon}, where the two special segments are highlighted in purple. When $w$ is in an odd chamber, we call $w_1 w_2$ (resp.\ $w_4 w_5$) a \emph{special edge}, even if its special segment is empty (which happens when the opposite diagonals do not cross in the hexagon interior). 

To prove Theorem \ref{t.hexagon}, we will need some preliminary lemmas.

\begin{Lem} \label{lem:sameSide}
In the setup of \eqref{e:hexagonVertices}, $w_{1}$ and $w_{5}$ are on the same side of $H_{\gg,k}$ as $w$.
\end{Lem}

\begin{proof}
We will give the proof for $w_{1}$; the proof for $w_{5}$ is almost identical. Notice that for any of the six chambers, the line $H_{\gg,k}$ is either disjoint from the chamber closure (for odd chambers), or intersects it only at its vertex (for even chambers). In the second case, the reflection $s_{\ga,i}$ takes $H_{\gb,j}$ to $H_{\gg,k}$, and hence takes points in the chamber (such as $w$) to points on the same side of $H_{\gg,k}$. In the first case, the reflection takes $H_{\gb,j}$ to the translate of $H_{\gg,k}$ passing through the vertex of the chamber closure (i.e., closer to the chamber than $H_{\gg,k}$). So again, reflections of points in the chamber certainly stay on the same side of $H_{\gg,k}$.
\end{proof}

\begin{Lem} \label{lem:reflDown}
Let $w\in W$ be a non-spiral element, $\gd\in\{\tilde\ga,\ga_{1},\ga_{2}\}$, and $H_{\gd,\gre}$ for $\gre\in\{0,1\}$ be the nearer of $H_{\gd,0}, H_{\gd,1}$ to $wq$. Then $s_{\gd,\gre}w < w$.
\end{Lem}

\begin{proof}
This follows immediately from Lemma \ref{l:lengthChange}.
\end{proof}

\begin{proof}[Proof of Theorem \ref{t.hexagon}]
We first show that all points $xq$ in the hexagon satisfy $x\le w$. A key ingredient is the Endpoint Theorem of Graham-Li, \cite[Theorem 5.5]{GrLi:21}, which says that if $xq, yq, zq$ lie on a root string with $yq$ between $xq$ and $zq$, then either $y\le x$ or $y\le z$. By Lemma \ref{lem:reflDown}, we have $w_{1}, w_{3}, w_{5} < w$. But Lemma \ref{lem:sameSide} implies that we may again use Lemma \ref{lem:reflDown} to conclude that $w_{2}< w_{1}$ and $w_{4}< w_{5}$. So the hexagon vertices satisfy $w_{i}\le w$ for $0\le i\le 5$. Any hexagon edge point lies on a root string interval determined by the endpoints of the edge, so the edge points are all $\le w$ by the Endpoint Theorem. Finally, any interior point of the hexagon lies on a root string interval (in fact, three of them) with endpoints on the hexagon edges. So these are all $\le w$ by the Endpoint Theorem again.

For the converse, we must show that any $xq$ outside the hexagon has $x \nleq w$. Evidently it is enough to show that no element lying on one of the six root strings just outside the hexagon, and parallel to the adjacent hexagon edge, is dominated by $w$. But it is clear geometrically that the spiral element on each such root string is less than every other element on the string: one can move from the spiral alcove to any other alcove on the string by a sequence of reflections across a wall of each alcove reached, and the number of hyperplanes separating the alcove from $A_{\circ}$ increases at each step, by an argument similar to the one used in the proof of Lemma \ref{lem:reflDown}. (See \cite[Theorem 5.1]{GrLi:21} for a more algebraic argument.)  So it suffices to prove that the spiral elements just outside the hexagon are $\nleq w$. For this we will use an inductive argument which will be a key technique for the remainder of the paper. 

We fix an arbitrary fundamental half-strip and in what follows, we will assume that all spiral elements considered lie in this fixed half-strip. Our induction hypothesis is that, for some given (non-spiral) $w$, if $y$ is the spiral element (in the fixed fundamental half-strip) on the boundary of $\ch_{w}$, and $z$ is the next spiral element just outside $\ch_{w}$, then $y\le w$ but $z\nleq w$. 
The induction step will be to prove the same for the element(s) $w'$ obtained by reflecting $w$ across a wall of $wA_{\circ}$ and having $\ell(w')=\ell(w)+1$. (The base cases, where $w$ belongs to the lowest alcove in one of the six chambers, can be checked by a direct computation.)

Before beginning the induction step, we make a simple observation. If $w=t_{1}\dots t_{r}$ is a reduced expression (here the $t_{j}$ are simple reflections in $W$), and if $y=s_{1}\dots s_{k}$ is a subexpression of $w$, say $s_{i} = t_{j_{i}}$ for some indices $1\le j_{1}<\dots<j_{k}\le r$, then, for $i = 1, 2, \dots, k$, we may choose each $j_{i}$ in turn so that $t_{j_{i}}$ is the first occurrence of $s_{i}$ to the right of $t_{j_{{i-1}}}$ in the given expression for $w$. (When $i=1$, omit the phrase ``to the right of $t_{j_{{i-1}}}$''.) We call this the leftmost subexpression $s_{1}\dots s_{k}$ in $t_{1}\dots t_{r}$. We will assume in what follows that all subexpressions chosen are the leftmost ones. Recall also that each spiral element $y$ has a unique reduced expression. This means that $y\le w$ if and only if the unique reduced expression for $y$ is a subexpression of any (and every) reduced expression for $w$.

\begin{figure}[htbp!]
%

\pgfmathsetmacro{\cols}{6}
\pgfmathsetmacro{\rows}{6}
\pgfmathsetmacro{\slant}{cot(60)}
\pgfmathsetmacro{\height}{0.5 * \rows * tan(60)}
\pgfmathsetmacro{\triht}{sin(60)}
\pgfmathsetmacro{\upmid}{0.25 * sec(30)}
\pgfmathsetmacro{\downmid}{\triht - \upmid}

\begin{tikzpicture}
    

\input{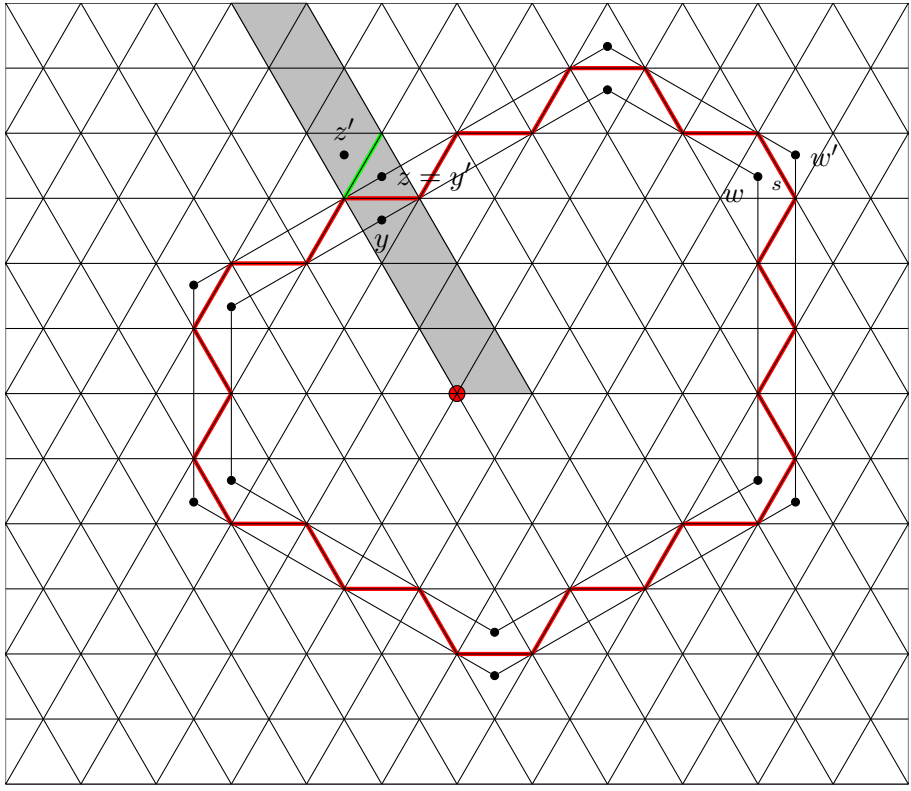}

	\begin{pgfonlayer}{nodelayer}
		\node [style=red dot] (4) at (0, 0) {};
		\node [style=small black dot, label={below left:$w$}] (13) at (4, 3 * \triht + \upmid) {};
		\node [style=small black dot] (14) at (2, 4 * \triht + \downmid) {};
		\node [style=small black dot] (19) at (4, -2 * \triht + \downmid) {};
		\node [style=small black dot] (20) at (0.5, -4 * \triht + \upmid) {};
		\node [style=small black dot] (21) at (-3, -2 * \triht + \downmid) {};
		\node [style=small black dot] (22) at (-3, 1 * \triht + \upmid) {};
		\node [style=small black dot, label={right:$w'$}] (23) at (4.5, 3 * \triht + \downmid) {};
		\node [style=small black dot] (24) at (2, 5 * \triht + \upmid) {};
		\node [style=small black dot] (25) at (-3.5, 1 * \triht + \downmid) {};
		\node [style=small black dot] (26) at (-3.5, -2 * \triht + \upmid) {};
		\node [style=small black dot] (27) at (0.5, -5 * \triht + \downmid) {};
		\node [style=small black dot] (28) at (4.5, -2 * \triht + \upmid) {};
		\node [style=none] (29) at (4, 4 * \triht) {};
		\node [style=none] (30) at (3, 4 * \triht) {};
		\node [style=none] (31) at (2.5, 5 * \triht) {};
		\node [style=none] (32) at (1.5, 5 * \triht) {};
		\node [style=none] (33) at (1, 4 * \triht) {};
		\node [style=none] (34) at (0, 4 * \triht) {};
		\node [style=none] (35) at (-0.5, 3 * \triht) {};
		\node [style=none] (36) at (-1.5, 3 * \triht) {};
		\node [style=none] (37) at (-2, 2 * \triht) {};
		\node [style=none] (38) at (-3, 2 * \triht) {};
		\node [style=none] (39) at (-3.5, 1 * \triht) {};
		\node [style=none] (40) at (-3, 0) {};
		\node [style=none] (41) at (-3.5, -1 * \triht) {};
		\node [style=none] (42) at (-3, -2 * \triht) {};
		\node [style=none] (43) at (-2, -2 * \triht) {};
		\node [style=none] (44) at (-1.5, -3 * \triht) {};
		\node [style=none] (45) at (-0.5, -3 * \triht) {};
		\node [style=none] (46) at (0, -4 * \triht) {};
		\node [style=none] (47) at (1, -4 * \triht) {};
		\node [style=none] (48) at (1.5, -3 * \triht) {};
		\node [style=none] (49) at (2.5, -3 * \triht) {};
		\node [style=none] (50) at (3, -2 * \triht) {};
		\node [style=none] (51) at (4, -2 * \triht) {};
		\node [style=none] (52) at (4.5, -1 * \triht) {};
		\node [style=none] (53) at (4, 0) {};
		\node [style=none] (54) at (4.5, 1 * \triht) {};
		\node [style=none] (55) at (4, 2 * \triht) {};
		\node [style=none] (56) at (4.5, 3 * \triht) {};
		\node [style=none] (57) at (4.25, 3.2 * \triht) {$\scriptstyle s$};
		\node [style=small black dot, label={below:$y$}] (58) at (-1, 2 * \triht + \downmid) {};
		\node [style=small black dot, label={right:$z=y'$}] (59) at (-1, 3 * \triht + \upmid) {};
		\node [style=small black dot, label={above:$z'$}] (60) at (-1.5, 3 * \triht + \downmid) {};
		\node [style=none] (61) at (1, 0) {};
		\node [style=none] (62) at (-2, 6 * \triht) {};
		\node [style=none] (63) at (-3, 6 * \triht) {};
		\node [style=none] (64) at (-1, 4 * \triht) {};
		\node [style=none] (65) at (-1.5, 3 * \triht) {};
	\end{pgfonlayer}
	\begin{pgfonlayer}{edgelayer}
		\draw [style=filled grey] (4.center)
			 to (61.center)
			 to (62.center)
			 to (63.center)
			 to cycle;
		\draw (13) to (14);
		\draw (19) to (13);
		\draw (20) to (19);
		\draw (21) to (20);
		\draw (14) to (22);
		\draw (22) to (21);
		\draw (23) to (24);
		\draw (24) to (25);
		\draw (25) to (26);
		\draw (26) to (27);
		\draw (27) to (28);
		\draw (28) to (23);
		\draw [style=red line] (56.center) to (29.center);
		\draw [style=red line] (29.center) to (30.center);
		\draw [style=red line] (30.center) to (31.center);
		\draw [style=red line] (31.center) to (32.center);
		\draw [style=red line] (32.center) to (33.center);
		\draw [style=red line] (33.center) to (34.center);
		\draw [style=red line] (34.center) to (35.center);
		\draw [style=red line] (35.center) to (36.center);
		\draw [style=red line] (36.center) to (37.center);
		\draw [style=red line] (37.center) to (38.center);
		\draw [style=red line] (38.center) to (39.center);
		\draw [style=red line] (39.center) to (40.center);
		\draw [style=red line] (40.center) to (41.center);
		\draw [style=red line] (41.center) to (42.center);
		\draw [style=red line] (42.center) to (43.center);
		\draw [style=red line] (43.center) to (44.center);
		\draw [style=red line] (44.center) to (45.center);
		\draw [style=red line] (45.center) to (46.center);
		\draw [style=red line] (46.center) to (47.center);
		\draw [style=red line] (47.center) to (48.center);
		\draw [style=red line] (48.center) to (49.center);
		\draw [style=red line] (49.center) to (50.center);
		\draw [style=red line] (50.center) to (51.center);
		\draw [style=red line] (51.center) to (52.center);
		\draw [style=red line] (52.center) to (53.center);
		\draw [style=red line] (53.center) to (54.center);
		\draw [style=red line] (54.center) to (55.center);
		\draw [style=red line] (55.center) to (56.center);
		\draw [style=green line] (64.center) to (65.center);
	\end{pgfonlayer}
	

    \clip       (-\cols, -\height) rectangle (\cols, \height);
    \draw[gray] (-\cols, -\height) rectangle (\cols, \height);

    \pgfmathsetmacro{\from}{-2 *\cols}
    \pgfmathsetmacro{\to}{2 * \cols}
    \foreach\i in {\from, ..., \to} {
        \draw[xslant=\slant]  (\i, -\height) -- (\i, \height);
        \draw[xslant=-\slant] (\i, -\height) -- (\i, \height);
    }

    \foreach\j in {-\rows, ..., \rows} {
        \pgfmathsetmacro{\y}{0.5 * \j * tan(60)}
        \draw (-\cols, \y) -- (\cols, \y);
    }

\end{tikzpicture}

\caption{Relation between hexagons for $w$ (Type 1) and $w'=ws>w$}  \label{fig:move1s}
\end{figure}

\noindent
{\bf Case 1}: Suppose there is a unique simple reflection $s$ such that $w':=ws > w$. Then each vertex of $\ch_{w}$ has the same property (for the same $s$!), and the vertices of $\ch_{w'}$ are given by $w'_{i}=w_{i}s,\ 0\le i \le 5$. Moreover, the alcoves $A'$ whose centers lie along an edge of $\ch_{w'}$ are precisely the alcoves $A$ whose centers lie along the edges of $\ch_{w}$, reflected across their $s$ edges. See Figure \ref{fig:move1s}, where the red edges are all labeled by $s$, and one of the six spiral half-strips is shaded grey.

Recall the spiral elements $y$ and $z$, on and just outside the boundary of $\ch_{w}$. Then $y'=z=ys$ is the spiral element on the boundary of $\ch_{w'}$, and $z'=zt$ (for one of the simple reflections $t\ne s$) is the spiral element just outside $\ch_{w'}$. Clearly $y'$ is a subexpression of $w'$, but the leftmost (and, indeed, every) such subexpression must use the final $s$ of $w'=ws$, else $y'$ would be a subexpression of $w$. And clearly $z'$ is not a subexpression of $w'$ since there is no factor $t$ beyond the final $s$ of $w'$ and its subexpression $y'$.

\begin{figure}[htbp!]
%

\pgfmathsetmacro{\cols}{6}
\pgfmathsetmacro{\rows}{6}
\pgfmathsetmacro{\slant}{cot(60)}
\pgfmathsetmacro{\height}{0.5 * \rows * tan(60)}
\pgfmathsetmacro{\triht}{sin(60)}
\pgfmathsetmacro{\upmid}{0.25 * sec(30)}
\pgfmathsetmacro{\downmid}{\triht - \upmid}

\begin{tikzpicture}
    

\input{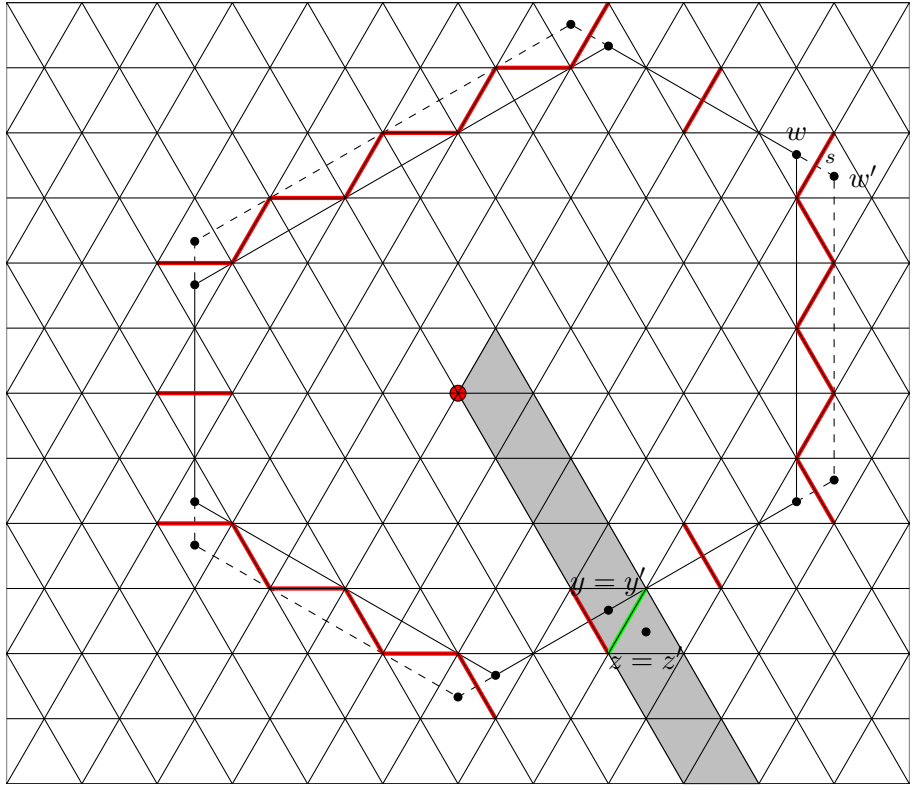}

	\begin{pgfonlayer}{nodelayer}
		\node [style=red dot] (4) at (0, 0) {};
		\node [style=small black dot, label={above:$w$}] (23) at (4.5, 3 * \triht + \downmid) {};
		\node [style=small black dot] (24) at (2, 5 * \triht + \upmid) {};
		\node [style=small black dot] (25) at (-3.5, 1 * \triht + \downmid) {};
		\node [style=small black dot] (26) at (-3.5, -2 * \triht + \upmid) {};
		\node [style=small black dot] (27) at (0.5, -5 * \triht + \downmid) {};
		\node [style=small black dot] (28) at (4.5, -2 * \triht + \upmid) {};
		\node [style=small black dot, label={right:$w'$}] (29) at (5, 3 * \triht + \upmid) {};
		\node [style=small black dot] (30) at (1.5, 5 * \triht + \downmid) {};
		\node [style=small black dot] (31) at (-3.5, 2 * \triht + \upmid) {};
		\node [style=small black dot] (32) at (-3.5, -3 * \triht + \downmid) {};
		\node [style=small black dot] (33) at (0, -5 * \triht + \upmid) {};
		\node [style=small black dot] (34) at (5, -2 * \triht + \downmid) {};
		\node [style=none] (35) at (5, 4 * \triht) {};
		\node [style=none] (36) at (4.5, 3 * \triht) {};
		\node [style=none] (37) at (5, 2 * \triht) {};
		\node [style=none] (38) at (4.5, 1 * \triht) {};
		\node [style=none] (39) at (5, 0) {};
		\node [style=none] (40) at (4.5, -1 * \triht) {};
		\node [style=none] (41) at (5, -2 * \triht) {};
		\node [style=none] (42) at (3.5, 5 * \triht) {};
		\node [style=none] (43) at (3, 4 * \triht) {};
		\node [style=none] (44) at (2, 6 * \triht) {};
		\node [style=none] (45) at (1.5, 5 * \triht) {};
		\node [style=none] (46) at (0.5, 5 * \triht) {};
		\node [style=none] (47) at (0, 4 * \triht) {};
		\node [style=none] (48) at (-1, 4 * \triht) {};
		\node [style=none] (49) at (-1.5, 3 * \triht) {};
		\node [style=none] (50) at (-2.5, 3 * \triht) {};
		\node [style=none] (51) at (-3, 2 * \triht) {};
		\node [style=none] (52) at (-4, 2 * \triht) {};
		\node [style=none] (53) at (-3, 0) {};
		\node [style=none] (54) at (-4, 0) {};
		\node [style=none] (55) at (-4, -2 * \triht) {};
		\node [style=none] (56) at (-3, -2 * \triht) {};
		\node [style=none] (57) at (-2.5, -3 * \triht) {};
		\node [style=none] (58) at (-1.5, -3 * \triht) {};
		\node [style=none] (59) at (-1, -4 * \triht) {};
		\node [style=none] (60) at (0, -4 * \triht) {};
		\node [style=none] (61) at (0.5, -5 * \triht) {};
		\node [style=none] (62) at (2, -4 * \triht) {};
		\node [style=none] (63) at (1.5, -3 * \triht) {};
		\node [style=none] (64) at (3.5, -3 * \triht) {};
		\node [style=none] (65) at (3, -2 * \triht) {};
		\node [style=none] (66) at (3, -6 * \triht) {};
		\node [style=none] (67) at (4, -6 * \triht) {};
		\node [style=none] (68) at (0.5, 1 * \triht) {};
		\node [style=small black dot, label={above:$y=y'$}] (69) at (2, -4 * \triht + \downmid) {};
		\node [style=small black dot, label={below:$z=z'$}] (70) at (2.5, -4 * \triht + \upmid) {};
		\node [style=none] (71) at (2.5, -3 * \triht) {};
		\node [style=none] (72) at (4.95, 3.6 * \triht) {$\scriptstyle s$};
	\end{pgfonlayer}
	\begin{pgfonlayer}{edgelayer}
		\draw [style=filled grey] (66.center)
			 to (67.center)
			 to (68.center)
			 to (4.center)
			 to cycle;
		\draw (23) to (24);
		\draw (24) to (25);
		\draw (25) to (26);
		\draw (26) to (27);
		\draw (27) to (28);
		\draw (28) to (23);
		\draw [style=dashed line] (23) to (29);
		\draw [style=dashed line] (30) to (24);
		\draw [style=dashed line] (25) to (31);
		\draw [style=dashed line] (32) to (26);
		\draw [style=dashed line] (27) to (33);
		\draw [style=dashed line] (34) to (28);
		\draw [style=dashed line] (29) to (34);
		\draw [style=dashed line] (30) to (31);
		\draw [style=dashed line] (32) to (33);
		\draw [style=red line] (35.center) to (36.center);
		\draw [style=red line] (36.center) to (37.center);
		\draw [style=red line] (37.center) to (38.center);
		\draw [style=red line] (38.center) to (39.center);
		\draw [style=red line] (39.center) to (40.center);
		\draw [style=red line] (40.center) to (41.center);
		\draw [style=red line] (42.center) to (43.center);
		\draw [style=red line] (44.center) to (45.center);
		\draw [style=red line] (45.center) to (46.center);
		\draw [style=red line] (46.center) to (47.center);
		\draw [style=red line] (47.center) to (48.center);
		\draw [style=red line] (48.center) to (49.center);
		\draw [style=red line] (49.center) to (50.center);
		\draw [style=red line] (50.center) to (51.center);
		\draw [style=red line] (51.center) to (52.center);
		\draw [style=red line] (53.center) to (54.center);
		\draw [style=red line] (55.center) to (56.center);
		\draw [style=red line] (56.center) to (57.center);
		\draw [style=red line] (57.center) to (58.center);
		\draw [style=red line] (58.center) to (59.center);
		\draw [style=red line] (59.center) to (60.center);
		\draw [style=red line] (60.center) to (61.center);
		\draw [style=red line] (63.center) to (62.center);
		\draw [style=red line] (65.center) to (64.center);
		\draw [style=green line] (62.center) to (71.center);
	\end{pgfonlayer}
	

    \clip       (-\cols, -\height) rectangle (\cols, \height);
    \draw[gray] (-\cols, -\height) rectangle (\cols, \height);

    \pgfmathsetmacro{\from}{-2 *\cols}
    \pgfmathsetmacro{\to}{2 * \cols}
    \foreach\i in {\from, ..., \to} {
        \draw[xslant=\slant]  (\i, -\height) -- (\i, \height);
        \draw[xslant=-\slant] (\i, -\height) -- (\i, \height);
    }

    \foreach\j in {-\rows, ..., \rows} {
        \pgfmathsetmacro{\y}{0.5 * \j * tan(60)}
        \draw (-\cols, \y) -- (\cols, \y);
    }

\end{tikzpicture}

\caption{Relation between hexagons for $w$ (Type 2) and $w'=ws>w$}  \label{fig:move2tu}
\end{figure}

\noindent
{\bf Case 2}: Suppose there are two simple reflections $s$
 for which $ws>w$. Fix one of them, and set $w':=ws$. Each of three non-intersecting edges of $\ch_{w'}$ contains an edge of $\ch_{w}$; the remaining three edges of $\ch_{w'}$ are one root string outside the parallel edges of $\ch_{w}$. (See Figure \ref{fig:move2tu}, where the red edges are all labeled by $s$, and one of the six spiral half-strips is shaded grey.) The alcoves along an edge of $\ch_{w'}$ which ``moved out'' are obtained by reflecting across their $s$ edges the alcoves along the adjacent parallel edge of $\ch_{w}$, as in Case 1. For the spiral elements along these three edges, the analysis is exactly the same as before.  

Reflection in the $s$ edge stabilizes the set of alcoves on the remaining three edges of $\ch_{w}$ (and on the collinear edges of $\ch_{w'}$).
For the spiral element $y$ on one of these edges of $\ch_{w}$, $y'=y$ and $z'=z=yt$ for some simple reflection $t\ne s$. Of course $y'=y\le w<w'$ by the induction hypothesis. And since $y$ is a subexpression of $w$ but $z=yt$ is not, $z'=yt$ could only be a subexpression of $w'=ws$ if $t = s$, which is not the case.
\end{proof}


\section{The translation move} \label{s:translation}
We will study the integers $q^w_x$ for $x \leq w$ using the following strategy.  The operation of translation into the chamber introduced
in Section \ref{s:chambers} gives us a ``Translation Move".  The main result of this section is that if $w$ is obtained from $w'$ by the Translation Move,
and $x$ is not in an outer shell of $\ch_w$, then $q^w_x = q^{w'}_x + 2$ (see Theorem \ref{t:translationmove}).  In Section \ref{s:outershells}, we
describe $q^w_x$ if $x$ is on an outer shell of $\ch_w$.  In Section \ref{s:basecases}, we describe $q^w_x$ if $w$ is a ``base case", that is, a non-spiral element
which is not obtained by a translation move.  The results of Sections \ref{s:translation},  \ref{s:outershells}, and \ref{s:basecases} together yield
a description of $q^w_x$ for any non-spiral element $w$.

\subsection{Hexagon shells} \label{ss.shells}
We call the union of the (generically, six) edges of $\ch_{w}$ its 0-shell. These edges lie on certain root strings, and the hexagon together with its interior consists of the intersection of the closed half-spaces containing $q$ determined by these six root strings. Now consider the root strings passing through the interior of the hexagon and one string in from an edge (the ``reference edge''). The intersection of the half-spaces on the side of these root strings not containing the reference edge is a convex region inside the hexagon, whose boundary (consisting of intervals lying along some or all of those root strings) we call the 1-shell of $\ch_{w}$. Repeating this process we obtain the 2-shell, and so on. Figure \ref{fig:trans-odd} 
shows the 0-, 1-, 2-, and 3-shells of a typical hexagon $\ch_{w'}$. See also Figure \ref{fig:edge3shell} for a situation where $\ch_{w}$ has some very short edges. 

\subsection{Translations} \label{ss:translation}
Recall from Section \ref{s:chambers} that for each chamber $\cc$, there is a unique root $\ga \in \Phi$ pointing into $\cc$.  Moreover,
given a non-spiral $w$,
there is an operation ``translation into the chamber" which takes $w$ to $w' = t(\alpha) w$, where $\alpha$ is the root pointing into the
chamber containing $w$.

\begin{figure}[htbp!]
%

\pgfmathsetmacro{\cols}{7}
\pgfmathsetmacro{\rows}{8}
\pgfmathsetmacro{\slant}{cot(60)}
\pgfmathsetmacro{\height}{0.5 * \rows * tan(60)}
\pgfmathsetmacro{\triht}{sin(60)}
\pgfmathsetmacro{\upmid}{0.25 * sec(30)}
\pgfmathsetmacro{\downmid}{\triht - \upmid}

\begin{tikzpicture}[scale=0.9]
    

\input{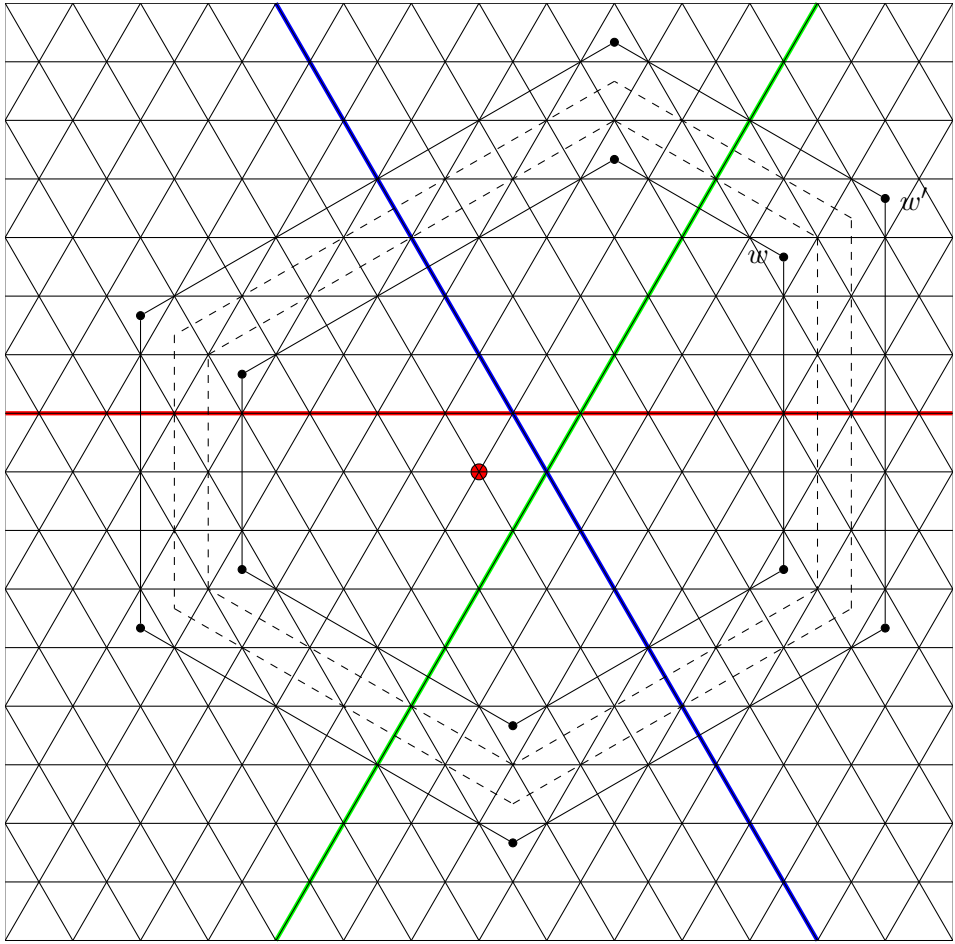}

	\begin{pgfonlayer}{nodelayer}
		\node [style=none] (2) at (-7, 1 * \triht) {};
		\node [style=none] (3) at (7, 1 * \triht) {};
		\node [style=red dot] (4) at (0, 0) {};
		\node [style=none] (7) at (-3, -8 * \triht) {};
		\node [style=none] (8) at (5, 8 * \triht) {};
		\node [style=none] (11) at (5, -8 * \triht) {};
		\node [style=none] (12) at (-3, 8 * \triht) {};
		\node [style=small black dot, label={right:$w'$}] (23) at (6, 4 * \triht + \downmid) {};
		\node [style=small black dot] (24) at (2, 7 * \triht + \upmid) {};
		\node [style=small black dot] (25) at (-5, 2 * \triht + \downmid) {};
		\node [style=small black dot] (26) at (-5, -3 * \triht + \upmid) {};
		\node [style=small black dot] (27) at (0.5, -7 * \triht + \downmid) {};
		\node [style=small black dot] (28) at (6, -3 * \triht + \upmid) {};
		\node [style=none] (29) at (5.5, 4 * \triht + \upmid) {};
		\node [style=none] (30) at (2, 6 * \triht + \downmid) {};
		\node [style=none] (31) at (-4.5, 2 * \triht + \upmid) {};
		\node [style=none] (32) at (-4.5, -3 * \triht + \downmid) {};
		\node [style=none] (33) at (0.5, -6 * \triht + \upmid) {};
		\node [style=none] (99) at (5.5, -3 * \triht + \downmid) {};
		\node [style=none] (100) at (5, 4 * \triht ) {};
		\node [style=small black dot, label={left:$w$}] (101) at (4.5, 3 * \triht + \downmid) {};
		\node [style=small black dot] (102) at (2, 5 * \triht + \upmid) {};
		\node [style=small black dot] (105) at (-3.5, 1 * \triht + \downmid) {};
		\node [style=small black dot] (106) at (-3.5, -2 * \triht + \upmid) {};
		\node [style=small black dot] (107) at (0.5, -5 * \triht + \downmid) {};
		\node [style=small black dot] (108) at (4.5, -2 * \triht + \upmid) {};
		\node [style=none] (109) at (2, 6 * \triht) {};
		\node [style=none] (110) at (-4, 2 * \triht ) {};
		\node [style=none] (112) at (0.5, -5 * \triht) {};
		\node [style=none] (113) at (5, -2 * \triht) {};
		\node [style=none] (114) at (-4, -2 * \triht) {};
	\end{pgfonlayer}
	\begin{pgfonlayer}{edgelayer}
		\draw [style=red line] (2.center) to (3.center);
		\draw [style=green line] (7.center) to (8.center);
		\draw [style=blue line] (12.center) to (11.center);
		\draw (23) to (24);
		\draw (24) to (25);
		\draw (25) to (26);
		\draw (26) to (27);
		\draw (27) to (28);
		\draw (28) to (23);
		\draw [style=dashed line] (29.center) to (30.center);
		\draw (102) to (101);
		\draw (101) to (108);
		\draw (102) to (105);
		\draw (105) to (106);
		\draw (106) to (107);
		\draw (107) to (108);
		\draw [style=dashed line] (109.center) to (100.center);
		\draw [style=dashed line] (100.center) to (113.center);
		\draw [style=dashed line] (113.center) to (112.center);
		\draw [style=dashed line] (112.center) to (114.center);
		\draw [style=dashed line] (114.center) to (110.center);
		\draw [style=dashed line] (110.center) to (109.center);
		\draw [style=dashed line] (30.center) to (31.center);
		\draw [style=dashed line] (31.center) to (32.center);
		\draw [style=dashed line] (32.center) to (33.center);
		\draw [style=dashed line] (33.center) to (99.center);
		\draw [style=dashed line] (99.center) to (29.center);
	\end{pgfonlayer}	


    \clip       (-\cols, -\height) rectangle (\cols, \height);
    \draw[gray] (-\cols, -\height) rectangle (\cols, \height);

    \pgfmathsetmacro{\from}{-2 *\cols}
    \pgfmathsetmacro{\to}{2 * \cols}
    \foreach\i in {\from, ..., \to} {
        \draw[xslant=\slant]  (\i, -\height) -- (\i, \height);
        \draw[xslant=-\slant] (\i, -\height) -- (\i, \height);
    }

    \foreach\j in {-\rows, ..., \rows} {
        \pgfmathsetmacro{\y}{0.5 * \j * tan(60)}
        \draw (-\cols, \y) -- (\cols, \y);
    }

\end{tikzpicture}
\caption{Translation and outer shells}  \label{fig:trans-odd}
\end{figure}

\begin{Thm} \label{t:translationmove}
Suppose  $w$ is a non-spiral element of $W$.  Let $w'$ be the element obtained by translating $w$ into the chamber.
If $x \leq w$, then $q_x^{w'}= q_x^w + 2$.
\end{Thm}

\begin{proof}
By Proposition \ref{p:lengthChange}, $w' > w$, and $\ell(w') = \ell(w) + 4$.
Let $w = w_0, \ldots, w_5$ (resp.\ $w' = w'_0, \ldots, w'_5$ be the vertices of $\ch_w$ (resp.\ $\ch_{w'}$), numbered as Section \ref{s.hexagons}.  Let $\cc_0, \ldots, \cc_5$ be the chambers,
numbered counterclockwise, with $w_0 \in \cc_0$.  (If $w_i$ is not spiral, then $w_i$ lies in $\cc_i$.)  Let $\gg_i$ be the root pointing into $\cc_i$.
From the formulas for $w_i$ and 
$w'_i$, we see that $w'_i = t(\gg_i) w_i$.

Let $\ga \in \Phi$, and let $L$ be an $\ga$-root string which intersects $\ch_w$.  Let $p$ and $p'$ be the extreme points of $L \cap \ch_w$, labelled so that
$p-p'$ is a positive multiple of $\ga$.  Then 
we claim that the extreme points of $L \cap \ch_{w'}$ are $p+ \ga$ and $p' - \ga$.  Indeed, suppose $p$ lies on the edge $E=w_{i-1}w_{i}$ of $\ch_{w}$, where  subscripts are interpreted mod 6. (If $p$ happens to be a vertex of $\ch_{w}$, choose $E$ so that it is not parallel to $\ga$.) Suppose $\ga=\gg_{i}$. Note that $w_{i}+\ga=w'_{i}$, an endpoint of the edge $E'=w'_{i-1}w'_{i}$ of $\ch_{w'}$, and  $E$ and $E'$ are parallel, with $E'$ longer. It follows by elementary geometry that for any point, such as $p$, on $E$, the point $p+\ga$ lies on $E'$. The proof is similar if $\ga=\gg_{i-1}$ using the other endpoints of $E$ and $E'$. The third option, $\ga=\pm\gg_{i+1}$, is not actually a possibility because $\gg_{i+1}$ is parallel to the edge $E$, and we chose $E$ to exclude this situation. The proof for $p'$ is analogous.  This proves the claim.

We claim that each of the intervals (i.e., half-open line segments) $(p, p + \ga]$
and $(p',p'-\ga]$ contains two alcove centers, one of each orientation.  Indeed, consider the interval $(p, p + \ga]$, the argument
for $(p',p'-\ga]$ being similar.  Then either $p$ is the center of an Up alcove, the center of a Down alcove, or 
an alcove vertex.  By symmetry, it suffices to check the claim for a single point $p$ for
each of these $3$ possibilities, and this can be done by inspection.  

Given $x \leq w$, there are $3$ root strings passing through $x$, each of which intersects $\ch_w$.  If $L$ is such a root string, the preceding paragraph
implies that the alcove centers on $L \cap \ch_{w'}$ are the alcove centers on $L \cap \ch_w$, together with $4$ additional alcoves, $2$ of each
orientation.  Therefore, there are $6$ additional reflections ($2$ for each $L$) which, when applied to $x$, stay in the hexagon $\ch_{w'}$ (as compared to the hexagon
$\ch_w$).   But $\ell(w) - \ell(w') =  -4$ by Proposition \ref{p:lengthChange}, so $q^{w'}_x = q^w_x + 2$, as desired.
\end{proof}

\begin{Cor} \label{c:qpositive}
Let $w'$ be a non-spiral element of $W$, and assume $w'$ is not a base case. 

\noindent (a) For all $x$ on or inside the 3-shell of $\ch_{w'}$, $q^{w'}_{x}\ge 2$. 

\noindent (b) If $w'$ is of Type 2, then for all $x$ on the 2-shell of $\ch_{w'}$, $q^{w'}_{x}\ge 1$.

\noindent
In particular, the only $x\le w'$ for which $q^{w'}_{x}=0$ are on the 0-, 1-, or 2-shell, and are explicitly described in Theorem \ref{t:outershells}.
\end{Cor}

\begin{proof}
Any $w'$ as in the statement of the corollary is the translation into the chamber of some $w$ in the chamber, as in Theorem \ref{t:translationmove}. Moreover the boundary of $\ch_{w}$ is the 3-shell of $\ch_{w'}$. So part (a) of the corollary follows immediately from the last statement of the theorem, along with the fact that $q^{w}_{x}\ge 0$ for all $x\le w$. Part (b) follows from Theorem \ref{t:outershells} (c) and (d), below.
\end{proof}

\begin{Rem} \label{r:base-rs}
This corollary implies that if a non-spiral $X_w$ is rationally smooth, then $w$ is necessarily a base case.  The analysis of the base cases in Section \ref{s:basecases} 
then allows us to determine the rationally smooth Schubert varieties in type $\tilde{A}_2$ (Corollary \ref{c:ratsmoothSchubertvar}).
\end{Rem}

Figure \ref{fig:trans-odd} 
illustrates the translation move in Chamber I.
Rays along the red and green lines give the walls of the chamber containing $w$ and $w'$, and the blue line is the third hyperplane used to obtain the hexagons $\ch_w$ and $\ch_{w'}$.  

\begin{figure}[htbp!]
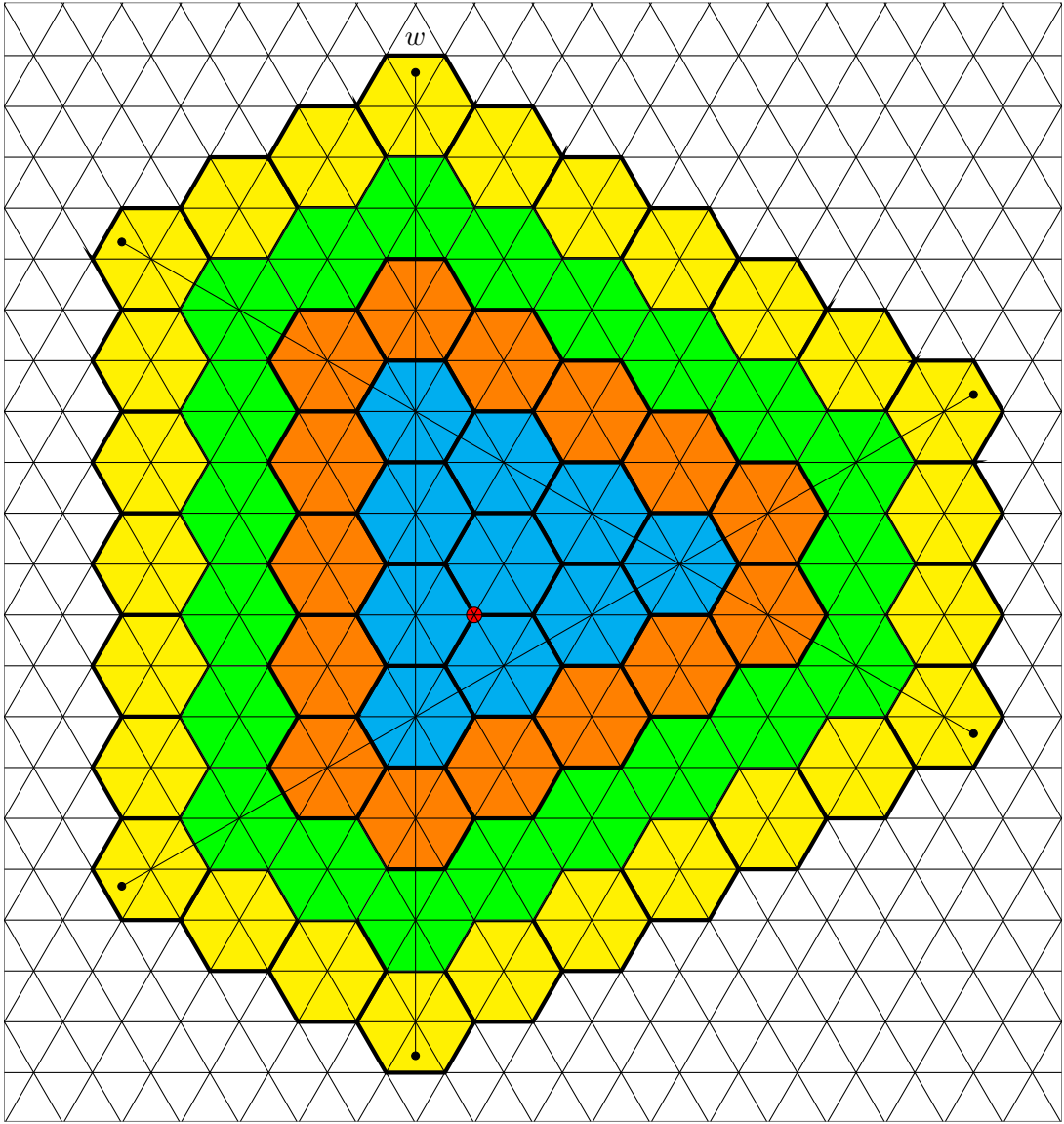

\include{Fig-qShellsEvenType1}
\caption{Shell behavior of $q$}  \label{fig:qshells}
\end{figure}

Following the strategy outlined at the beginning of this section, we can describe $q^w_x$ in a fairly explicit way for any non-spiral $w$.  In any chamber $\cc$, there are two types of base cases:
the alcoves which share an edge with a fundamental strip, or the alcoves which share a vertex with a fundamental strip.   So there are essentially four types of 
base cases in total: two for odd chambers, and two for even chambers.  Any $w$ in $\cc$ is either a base case, or is obtained from a base case by a sequence of translations by the root  $\alpha$ pointing into $\cc$.  Therefore we can proceed as follows.  First, if $w$ is a base case, we will describe (in Section \ref{s:basecases}) $q^w_x$ for all $x \leq w$.  Second, if $w' \in \cc$ and $w' = t(\alpha) w$, then for $x \leq w$, we can assume $q^w_x$ is known by induction. Therefore, for $x \leq w'$, if $x$ is not on the $0$-, $1$-, or $2$-shell of $\ch_{w'}$, then $x \leq w$ and $q^{w'}_x = q^w_x + 2$. To complete the picture,
we will describe $q^{w'}_x$ for $x$ on the $0$-, $1$-, or $2$-shell of 
$\ch_{w'}$ (in Section \ref{s:outershells}). 

Figure \ref{fig:qshells} shows the behavior of $q$ for the simplest situation: $w$ is of Type 1 in an even chamber. This means $w$ is obtained by translation from Base Case I.  The values of $q$ are constant on the orbits $x R(w)$, which are small $W(A_{2})$-hexagons: yellow 0, green 2, orange 4, and blue 6. More generally, $q$ is constant equal to $2k$ on the shells numbered $3k, 3k+1, 3k+2$ for $k=0, 1, \dots$, except that $q$ remains constant once the $W(A_{2})$-hexagon of $x$ meets the triangle formed by the hexagon diagonals.  If $w$ is of Type 2 in an even chamber, $w$ is obtained by translation from Base Case 2.  Here, the orbits $x R(w)$ are small $W(A_1)$-diamonds, and
the value of $q$ is constant on these orbits.  
More precisely, $q$ has value equal to $2k$ on shells numbered $3k$ and $3k+1$, and has value $2k+1$ on shells numbered $3k+2$, for $k=0, 1, \dots$, until $q$ remains constant once $x$ meets the triangle formed by the hexagon diagonals.
For $w$ in an odd chamber, the behavior of $q$ is more complicated.  Figures \ref{fig:basecase3} and \ref{fig:basecase4} illustrate the two base cases for $w$ in an odd chamber.


\section{Values of \texorpdfstring{$q$}{q} on outer hexagon shells} \label{s:outershells}

In this section, we fix $w\in W$ and compute explicitly the values $q_{x}^{w}$ for $x$ on or near the boundary of the hexagon $\ch_{w}$. This will handle the new elements $x$ arising in the Translation Move. More precisely, our goal here is to compute $q^{w}_{x}$ for $x$ on the 0-, 1-, or 2-shell of $\ch_{w}$.

We begin with some observations which will be used in the proof of the main result of this section, Theorem \ref{t:outershells}.
Assume $w$ is non-spiral of Type 1, with $w<ws=:w'$ for a (unique) simple reflection $s$. Then the vector from $w$ to $w'$ is a positive multiple of $\ga$, the root which points into the chamber containing $w$ (see Lemma \ref{l:rootchamber}). The root hyperplanes used to construct the vertices of $\ch_{w}$ and $\ch_{w'}$ are identical. Therefore these two hexagons are concentric and one root string apart along each edge. Moreover, by Lemma \ref{l:wallLabels} there is a bijection between alcoves $x$ along the edge of $\ch_{w}$ and alcoves $x':=xs$ along the edge of $\ch_{w'}$, as shown in Figure \ref{fig:move1s}, where the red edges all correspond to $s$.

\begin{Thm} \label{t:outershells}
Fix $w\in W$ a non-spiral element and assume $x\le w$ lies on the 0-, 1- or 2-shell of $\ch_{w}$.

\noindent (a) If $w$ is of Type 1 in an even chamber, then $q^{w}_{x}=0$.

\noindent (b) If $w$ is of Type 1 in an odd chamber, then $q^{w}_{x}=1$ if $x R(w)$ intersects a special segment of $\ch_{w}$; otherwise $q^{w}_{x}=0$.

\noindent (c) If $w$ is of Type 2 in an even chamber, then $q^{w}_{x}=0$ if $x$ belongs to the 0- or 1-shell of $\ch_{w}$; otherwise $q^{w}_{x}=1$.

\noindent (d) Assume $w$ is of Type 2 in an odd chamber. First suppose $x$ belongs to the 0- or 1-shell. Then $q^{w}_{x}=1$ if $x R(w)$ intersects a special segment; otherwise $q^{w}_{x}=0$. Next suppose $x$ belongs to the 2-shell. If $x$ belongs to a 2-shell edge which is parallel to and two root strings in from a special edge of $\ch_{w}$, then $q^{w}_{x}=2$; otherwise $q^{w}_{x}=1$.
\end{Thm}

As hinted in the proof of Theorem \ref{t.hexagon}, we will use an inductive argument involving reflection across alcove walls. More precisely, we will assume the result for a Type 1 element $w$ (such as the lowest alcove in each chamber), and prove it for $w':=ws > w$ (``Move 1''). Now $w'$ is of Type 2, and we will do another inductive step to prove the result for $w'':= w't > w'$ (``Move 2''). Here $s$ and $t$ are distinct simple reflections. This completes the cycle since $w''$ is again of Type 1. An important point is that we can reach any non-spiral element (except for the minimal element in each chamber) by a sequence of Moves 1 and 2 within a single chamber.

In order to keep track of the change in $q$ during each Move, we will introduce root-string-based integer coordinates for the alcoves in $\ch_{w}$ (where we are now assuming $w$ is of Type 1). See Figures \ref{fig:moves1and2even} and  \ref{fig:moves1and2odd}.
The $\ga_{0}$ strings
 in the hexagon will be labeled according to where they intersect the W or SW edges of the hexagon. First, label the alcoves whose centers lie on these edges, with label 0 at the SW vertex and increasing counterclockwise, resulting in positive labels along the SW edge and negative along the W edge. (For compactness of notation, we denote odd labels by a bar above instead of a negative sign in front.) So far we have labels for roughly 2/3 of the $\ga_{0}$ strings, but there are also strings which intersect the edge not at an alcove center but at an alcove vertex. We label these strings by using the label of the adjacent odd-numbered string, decorated with a * superscript. So the sequence of labels moving counterclockwise through the vertex will be $\dots, \bar 2, \bar 1^{*}, \bar 1, 0, 1, 1^{*}, 2, 3, 3^{*}, \dots$.

The $\ga_{1}$ strings are labeled according to their intersection with the SE or E edges, with 0 at the SE vertex, increasing counterclockwise for the alcove centers along the edge, and with additional odd* labels for the strings which interesect the edge in alcove vertices. And the $\ga_{2}$ strings are labeled similarly by their intersection with the NE and NW edges, with 0 at the N vertex and increasing counterclockwise. Finally, each alcove $x$ in the hexagon $\ch_{w}$ is labeled $(i_{0},i_{1},i_{2})$ where $i_{j}$ is the label on the $\ga_{j}$ string passing through $x$, for $j=0, 1, 2$. We will refer to each coordinate $i_{j}$ as either even, odd, or odd*; in particular ``odd'' implicitly means ``without a *''.

\begin{proof}
We give the proof (mainly) for the Up alcoves. The proof for the Down alcoves is similar, and the details will mostly be left to the reader.

\noindent {\bf Move 1}: Fix $w$ of Type 1, with $w<ws=:w'$, and assume the result is true for $w$. Note that $R(w') = \langle s \rangle$. 

Assume $x\le w$ is an Up alcove having coordinates $(i_{0},i_{1},i_{2})$. Recall that $q_{x}^{w}$ counts the total number of reflections $r$ such that $rx\le w$, minus $\ell(w)$. The reflections $r$ have the form $s_{\ga,k}$ for $\ga \in\{\ga_1, \ga_2, \tilde{\ga} \}$, and $k\in\Z$. 
For a fixed such $\ga$, the reflections $s_{\ga,k}$ contributing to $q_{x}^{w}$ are in bijection with the Down alcoves inside $\ch_{w}$ on the $\ga$ string through $x$. As described in the paragraph before Theorem \ref{t:outershells}, the hexagon $\ch_{w'}$ is obtained from $\ch_{w}$ by moving every edge out to the next adjacent root string. Thus the reflections of type $\ga$ contributing to $q_{x}^{w'}$ correspond to the same Down alcoves as for $w$, plus any additional Down alcoves at the ends of the $\ga$ string through $x$ in $\ch_{w'}$. It is easy to see that there is one such additional Down alcove when the coordinate $i_{j}$ is either even or odd*, and none when the coordinate is odd. Therefore the change in $q_{x}^{\bullet}$ for Move 1 from $w$ to $w'$ is

\begin{equation} \label{e:Move1deltaq}
\gD q_{1} := q_{x}^{w'} - q_{x}^{w} = \sum_{j=0}^{2} (\gd_{i_{j},e} + \gd_{i_{j},*}) - 1,
\end{equation}
where $\gd_{i_{j},e}$ (resp.\ $\gd_{i_{j},o}$, $\gd_{i_{j},*})$ is 1 if $i_{j}$ is even (resp.\ odd, odd*), and 0 otherwise. The ``$-1$'' is present because $\ell(w')=\ell(w)+1$. In words, $\Delta q_{1}$ for an Up alcove is one less than the number of even entries plus the number of odd* entries among the coordinates of $x$.

A similar analysis for Down alcoves yields almost the same formula. However, when $w$ is in an odd chamber, there are two pairs of root strings lying between parallel diagonals of the hexagon. These are the strings in Figure \ref{fig:moves1and2odd} labeled $i_{1}=\bar 1^{*}, \bar 1$ or $i_{2}=1, 1^{*}$. On these 4 strings, the roles of odd and odd* are reversed. The upshot is that $\Delta q_{1}$ for a Down alcove is one less than the number of even entries, plus the number of odd* entries not between two parallel diagonals, plus the number of odd entries between two parallel diagonals. We sometimes refer to as ``unusual'' these odd entries which contribute to $\Delta q_{1}$.

Continue to assume $x$ is an Up alcove. Suppose now that $x$ lies on a boundary edge $E$ of $\ch_{w}$; i.e., on the 1-shell of $\ch_{w'}$. Since the end alcoves of any hexagon edge are related by a reflection through a line which passes through an edge or a vertex (but not the center) of an alcove along the edge, the number of alcoves on the edge is even. Since one end alcove is labeled 0, the other end alcove has an odd label. It follows that the coordinate which is constant along $E$ is odd. The coordinate for an Up alcove along $E$ is always odd, because the alcoves labeled 0---the N, SE, and SW corners of $\ch_{w}$---are always Down for $w$ of Type~1. Finally, the remaining coordinate of an Up alcove on $E$ is always even. Therefore the coordinates of $x$ are some permutation of $(o,o,e)$, and it follows from \eqref{e:Move1deltaq} that $\gD q_{1} = 0$; i.e., $q_{x}^{w'} = q_{x}^{w}$. We leave it to the reader to check that the same formula is true for the Down alcoves on $E$.

By the induction hypothesis, $q_{x}^{w}=0$, unless $x$ lies on a special segment of $E$, in which case $q_{x}^{w}=1$. Recall from the paragraph before Theorem \ref{t:outershells} that each alcove $x'$ of the corresponding edge $E'$ of $\ch_{w'}$ is of the form $x'=xs>x$ for an $x$ on $E$. By the Simple Move (Proposition \ref{p.summary}), since $w'>w's$, $q_{x'}^{w'}=q_{x}^{w'}$. Moreover, since the diagonals which determine the special segments for $\ch_{w'}$ and $\ch_{w}$ are collinear, $x'$ is on a special segment of $E'$ if and only if $x$ is on a special segment of $E$. Combining the formulas $q_{x'}^{w'}=q_{x}^{w'}$ and $q_{x}^{w'} = q_{x}^{w}$ thus gives the desired values of $q$ on the 0- and 1-shells for $w'$.

It remains to compute $q_{x}^{w'}$ for $x$ on the 2-shell of $\ch_{w'}$; i.e., $x$ is on some edge $E$ of the 1-shell of $\ch_{w}$. 

First, assume $x$ does not lie between two parallel diagonals of $\ch_{w}$. By a similar analysis as before, the coordinates of $x$ are some permutation of $(e,o,o^{*})$, so $\Delta q = 1+1-1 = 1$. If $x$ is not (resp.\ is) a reflection across one of its alcove walls from a special segment alcove of $\ch_{w}$, we have $q^{w}_{x} = 0$ (resp.\ $1$), and thus $q^{w'}_{x}=1$ (resp.\ $2$). 

Finally, assume $x$ \emph{does} lie between two parallel diagonals of $\ch_{w}$. Such alcoves occur in pairs forming a ``diamond'' at one end or the other of $E$ (and $E$ must be one root string in from and parallel to a special edge of $\ch_{w}$). Assuming $E$ has length at least 4, there will be two such diamonds on $E$. The outer (end) alcoves $x$ have $q^{w}_{x}=0$. The outer Up alcove has coordinates (some permutation of) $(e,o^{*},o^{*})$, so $\Delta q=1+1+1-1=2$, whence $q^{w'}_{x}=0+2=2$. The outer Down alcove has coordinates (some permutation of) $(e,o^{*},o)$, but the ``$o$'' is unusual (counts towards $\Delta q$) because $x$ is a Down alcove between parallel diagonals. So again $\Delta q=2$ and $q^{w'}_{x}=2$. The inner (next-to-the-end) diamond alcoves $x$ have $q^{w}_{x}=1$ (because $xR(w)$ intersects a special segment). The Up alcove coordinates are (some permutation of) $(e,o,o^{*})$, so $\Delta q=1+1-1=1$, whence $q^{w'}_{x}=1+1=2$. The Down alcove coordinates are (some permutation of) $(e,o,o)$, but one ``$o$'' is unusual, and again $q^{w'}_{x}=1+1=2$. Lastly, if $E$ has length 2, then there is just one diamond of alcoves between parallel diagonals, but they lie between \emph{both} pairs of parallel diagonals. The adjacent special edge of $\ch_{w}$ has empty special segment, so $q^{w}_{x}=0$. The Up alcove coordinates are (some permutation of) $(e,o^{*},o^{*})$, so $\Delta q=1+1+1-1=2$, and $q^{w'}_{x}=0+2=2$. The Down alcove coordinates are (some permutation of) $(e,o,o)$, but \emph{both} $o$'s are unusual, so again $\Delta q=2$, and $q^{w'}_{x}=0+2=2$.

This completes the verification of the values of $q$ for the Type 2 element $w'$ on its 0-, 1- and 2-shells.

\noindent {\bf Move 2}: We assume the result for $w$ and $w'=ws$ as above, and prove it for $w'':= w't > w'$ for one of the two possible simple reflections $t$. Because of the differing geometry of the hexagon in even and odd chambers, we must treat each case separately. We give the argument for a representative chamber of each parity; the result for the other chambers of the same parity follows similarly. Likewise we give the argument only for one of the two possible $t$'s, the other following after the obvious adjustments.

\begin{figure}[htbp!]
%

\pgfmathsetmacro{\cols}{5}
\pgfmathsetmacro{\rows}{5}
\pgfmathsetmacro{\slant}{cot(60)}
\pgfmathsetmacro{\height}{0.5 * \rows * tan(60)}
\pgfmathsetmacro{\triht}{sin(60)}
\pgfmathsetmacro{\upmid}{0.25 * sec(30)}
\pgfmathsetmacro{\downmid}{\triht - \upmid}

\begin{tikzpicture}
    

\input{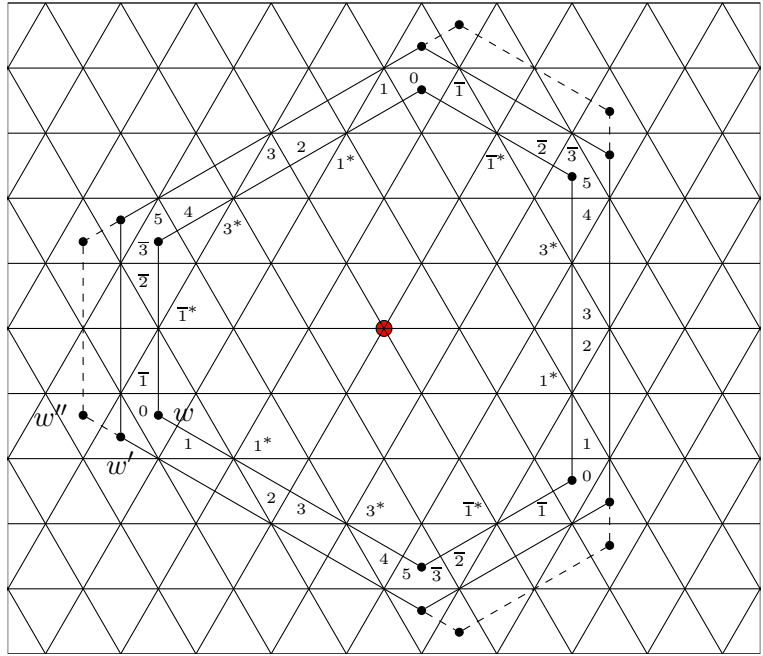}

	\begin{pgfonlayer}{nodelayer}
		\node [style=red dot] (4) at (0, 0) {};
		\node [style=small black dot] (13) at (2.5, 2 * \triht + \upmid) {};
		\node [style=small black dot] (14) at (0.5, 3 * \triht + \downmid) {};
		\node [style=small black dot] (19) at (2.5, -3 * \triht + \downmid) {};
		\node [style=small black dot] (20) at (0.5, -4 * \triht + \upmid) {};
		\node [style=small black dot, label={right:$w$}] (21) at (-3, -2 * \triht + \downmid) {};
		\node [style=small black dot] (22) at (-3, 1 * \triht + \upmid) {};
		\node [style=small black dot] (23) at (3, 2 * \triht + \downmid) {};
		\node [style=small black dot] (24) at (0.5, 4 * \triht + \upmid) {};
		\node [style=small black dot] (25) at (-3.5, 1 * \triht + \downmid) {};
		\node [style=small black dot, label={below:$w'$}] (26) at (-3.5, -2 * \triht + \upmid) {};
		\node [style=small black dot] (27) at (0.5, -5 * \triht + \downmid) {};
		\node [style=small black dot] (28) at (3, -3 * \triht + \upmid) {};
		\node [style=small black dot] (29) at (3, 3 * \triht + \upmid) {};
		\node [style=small black dot] (30) at (1, 4 * \triht + \downmid) {};
		\node [style=small black dot] (31) at (-4, 1 * \triht + \upmid) {};
		\node [style=small black dot, label={left:$w''$}] (32) at (-4, -2 * \triht + \downmid) {};
		\node [style=small black dot] (33) at (1, -5 * \triht + \upmid) {};
		\node [style=small black dot] (34) at (3, -4 * \triht + \downmid) {};
		\node [style=none] (35) at (-3, \triht + \upmid + .3) {$\scriptscriptstyle 5$};
		\node [style=none] (36) at (-2.6, \triht + \upmid + .4) {$\scriptscriptstyle 4$};
		\node [style=none] (37) at (-2, \triht + \upmid + .2) {$\scriptscriptstyle 3^{*}$};
		\node [style=none] (38) at (-1.5, 2* \triht + \upmid + .3) {$\scriptscriptstyle 3$};
		\node [style=none] (39) at (-1.1, 2* \triht + \upmid + .4) {$\scriptscriptstyle 2$};
		\node [style=none] (40) at (-0.5, 2* \triht + \upmid + .2) {$\scriptscriptstyle 1^{*}$};
		\node [style=none] (41) at (0, 3* \triht + \upmid + .3) {$\scriptscriptstyle 1$};
		\node [style=none] (42) at (0.4, 3* \triht + \upmid + .45) {$\scriptscriptstyle 0$};
		\node [style=none] (46) at (1, 3* \triht + \upmid + .3) {$\scriptscriptstyle \bar 1$};
		\node [style=none] (47) at (1.5, 2* \triht + \upmid + .2) {$\scriptscriptstyle \bar 1^{*}$};
		\node [style=none] (48) at (2.1, 2* \triht + \upmid + .4) {$\scriptscriptstyle \bar 2$};
		\node [style=none] (49) at (2.5, 2* \triht + \upmid + .3) {$\scriptscriptstyle \bar 3$};
		\node [style=none] (50) at (-3.2, \triht + .2) {$\scriptscriptstyle \bar 3$};
		\node [style=none] (51) at (-3.2, \upmid + .35) {$\scriptscriptstyle \bar 2$};
		\node [style=none] (52) at (-2.6, .2) {$\scriptscriptstyle \bar 1^{*}$};
		\node [style=none] (53) at (-3.2, -1* \triht + .2) {$\scriptscriptstyle \bar 1$};
		\node [style=none] (54) at (-3.2, -2 * \triht + \upmid + .35) {$\scriptscriptstyle 0$};
		\node [style=none] (55) at (-2.6, -2 * \triht + .2) {$\scriptscriptstyle 1$};
		\node [style=none] (56) at (-1.6, -2 * \triht + .2) {$\scriptscriptstyle 1^{*}$};
		\node [style=none] (57) at (-1.5, -3 * \triht + .35) {$\scriptscriptstyle 2$};
		\node [style=none] (58) at (-1.1, -3 * \triht + .2) {$\scriptscriptstyle 3$};
		\node [style=none] (59) at (-0.1, -3 * \triht + .2) {$\scriptscriptstyle 3^{*}$};
		\node [style=none] (60) at (0, -4 * \triht + .4) {$\scriptscriptstyle 4$};
		\node [style=none] (61) at (0.3, -4 * \triht + .2) {$\scriptscriptstyle 5$};
		\node [style=none] (62) at (0.7, -4 * \triht + .2) {$\scriptscriptstyle \bar 3$};
		\node [style=none] (63) at (1, -4 * \triht + .4) {$\scriptscriptstyle \bar 2$};
		\node [style=none] (64) at (1.2, -3 * \triht + .2) {$\scriptscriptstyle \bar 1^{*}$};
		\node [style=none] (65) at (2.1, -3 * \triht + .2) {$\scriptscriptstyle \bar 1$};
		\node [style=none] (66) at (2.7, -3 * \triht + \upmid + .35) {$\scriptscriptstyle 0$};
		\node [style=none] (70) at (2.7, -2* \triht + .2) {$\scriptscriptstyle 1$};
		\node [style=none] (71) at (2.2, -1 * \triht + .2) {$\scriptscriptstyle 1^{*}$};
		\node [style=none] (72) at (2.7, -1 * \triht +\upmid + .35) {$\scriptscriptstyle 2$};
		\node [style=none] (73) at (2.7,  .2) {$\scriptscriptstyle 3$};
		\node [style=none] (74) at (2.2, 1 * \triht + .2) {$\scriptscriptstyle 3^{*}$};
		\node [style=none] (75) at (2.7, 1 * \triht + \upmid + .35) {$\scriptscriptstyle 4$};
		\node [style=none] (76) at (2.7, 2* \triht + .2) {$\scriptscriptstyle 5$};
	\end{pgfonlayer}
	\begin{pgfonlayer}{edgelayer}
		\draw (13) to (14);
		\draw (19) to (13);
		\draw (20) to (19);
		\draw (21) to (20);
		\draw (14) to (22);
		\draw (22) to (21);
		\draw (23) to (24);
		\draw (24) to (25);
		\draw (25) to (26);
		\draw (26) to (27);
		\draw (27) to (28);
		\draw (28) to (23);
		\draw [style=dashed line] (23) to (29);
		\draw [style=dashed line] (29) to (30);
		\draw [style=dashed line] (30) to (24);
		\draw [style=dashed line] (25) to (31);
		\draw [style=dashed line] (31) to (32);
		\draw [style=dashed line] (32) to (26);
		\draw [style=dashed line] (27) to (33);
		\draw [style=dashed line] (33) to (34);
		\draw [style=dashed line] (34) to (28);
	\end{pgfonlayer}
	

    \clip       (-\cols, -\height) rectangle (\cols, \height);
    \draw[gray] (-\cols, -\height) rectangle (\cols, \height);

    \pgfmathsetmacro{\from}{-2 *\cols}
    \pgfmathsetmacro{\to}{2 * \cols}
    \foreach\i in {\from, ..., \to} {
        \draw[xslant=\slant]  (\i, -\height) -- (\i, \height);
        \draw[xslant=-\slant] (\i, -\height) -- (\i, \height);
    }

    \foreach\j in {-\rows, ..., \rows} {
        \pgfmathsetmacro{\y}{0.5 * \j * tan(60)}
        \draw (-\cols, \y) -- (\cols, \y);
    }

\end{tikzpicture}

\caption{Coordinates for Moves 1 and 2, even chamber.}  \label{fig:moves1and2even}
\end{figure}

{\bf Even chambers}: Assume first that $w, w'$ and $w''$ belong to chamber IV. Because the hexagon diagonal endpoints are a pair of opposite vertices, there are no special edges. We will take $t$ to be such that $w''$ is outward and to the right from $w'$, when viewed from the origin. (This is the situation, for example, when $w, w'$ and $w''$ all lie in the strip $-1 < (\ga_{2},v) < 0$, just below the upper boundary of chamber IV, as in Figure \ref{fig:moves1and2even}; note that although $w''$ is to the left of $w'$ in the picture, it is to the right when viewed from the origin.) Then the E, NW, and SW edges of $\ch_{w'}$ and of $\ch_{w''}$ are collinear; the NE, W, and SE edges of $\ch_{w''}$ are one root string farther out than the corresponding edges of $\ch_{w'}$. We have $w''>w''s$, and for $y$ a boundary point of $\ch_{w''}$, the six alcoves $yu$ for $u\in R(w'') = \langle\, s, t \, \rangle \cong W(A_{2})$ are arranged in a small hexagon all sharing a common vertex, and lying inside $\ch_{w''}$ (along its 0-, 1-, and 2-shells). By Proposition \ref{p.summary}, 
\begin{equation} \label{e:smallHexagon}
q_{y}^{w''}=q_{yu}^{w''} \quad \text{for all such }u.
\end{equation}
Moreover, the union of these small hexagons equals the 0-, 1- and 2-shells of $\ch_{w''}$, and every such hexagon contains at least one alcove along the boundary of $\ch_{w}$. In other words, to verify the theorem for $w''$, it suffices to verify it for $x$ on the 0-shell of $\ch_{w}$.

Repeating the analysis of the second paragraph of Move 1, we find that an Up alcove $x$ in $\ch_{w}$ with coordinates $(i_{0}, i_{1}, i_{2})$ acquires an additional Down alcove (in $\ch_{w''}$ but not $\ch_{w'}$) at the end of the $\ga_{j}$ string through $x$ ($j=$0, 1, or 2) if and only if $i_{j}$ is either odd and negative, or odd* and negative. Therefore the change in $q_{x}^{\bullet}$ for Move 2 from $w'$ to $w''$ is
\begin{equation} \label{e:Move2EvenDeltaq}
\gD q_{2} := q_{x}^{w''} - q_{x}^{w'} = \sum_{j=0}^{2} (\gd_{i_{j},o<0} + \gd_{i_{j},*<0}) - 1,
\end{equation}
where the $\delta$ functions have the obvious meanings.

Consider $x$ an Up alcove on the 0-shell of $\ch_{w}$. As observed in Move 1, the coordinates of $x$ are some permutation of $(o,o,e)$ and $q_{x}^{w'}=0$ (since there are no special segments). In fact, one odd coordinate is positive and the other is negative, so $\gD q_{2} = 0$ and thus $q_{x}^{w''}=0$. The reader can check that the same is true for the Down alcoves. By the remark in the previous paragraph, this completes the proof for Type 1 elements in even chambers.

\begin{figure}[htbp!]
%

\pgfmathsetmacro{\cols}{6}
\pgfmathsetmacro{\rows}{6}
\pgfmathsetmacro{\slant}{cot(60)}
\pgfmathsetmacro{\height}{0.5 * \rows * tan(60)}
\pgfmathsetmacro{\triht}{sin(60)}
\pgfmathsetmacro{\upmid}{0.25 * sec(30)}
\pgfmathsetmacro{\downmid}{\triht - \upmid}

\begin{tikzpicture}
    

\input{sample.tikzstyles}

	\begin{pgfonlayer}{nodelayer}
		\node [style=red dot] (4) at (0, 0) {};
		\node [style=small black dot, label={below left:$w$}] (13) at (4, 3 * \triht + \upmid) {};
		\node [style=small black dot] (14) at (2, 4 * \triht + \downmid) {};
		\node [style=small black dot] (19) at (4, -1 * \triht - \upmid) {};
		\node [style=small black dot] (20) at (0.5, -3 * \triht - \downmid) {};
		\node [style=small black dot] (21) at (-3, -1 * \triht - \upmid) {};
		\node [style=small black dot] (22) at (-3, 1 * \triht + \upmid) {};
		\node [style=small black dot, label={above:$w'$}] (23) at (4.5, 3 * \triht + \downmid) {};
		\node [style=small black dot] (24) at (2, 5 * \triht + \upmid) {};
		\node [style=small black dot] (25) at (-3.5, 1 * \triht + \downmid) {};
		\node [style=small black dot] (26) at (-3.5, -1 * \triht - \downmid) {};
		\node [style=small black dot] (27) at (0.5, -4 * \triht - \upmid) {};
		\node [style=small black dot] (28) at (4.5, -1 * \triht - \downmid) {};
		\node [style=small black dot, label={above right:$w''$}] (29) at (5, 3 * \triht + \upmid) {};
		\node [style=small black dot] (30) at (1.5, 5 * \triht + \downmid) {};
		\node [style=small black dot] (31) at (-3.5, 2 * \triht + \upmid) {};
		\node [style=small black dot] (32) at (-3.5, -2 * \triht - \upmid) {};
		\node [style=small black dot] (33) at (0, -4 * \triht - \downmid) {};
		\node [style=small black dot] (34) at (5, -1 * \triht - \upmid) {};
		\node [style=none] (35) at (-3, \triht + \upmid + .3) {$\scriptscriptstyle 7$};
		\node [style=none] (36) at (-2.6, \triht + \upmid + .4) {$\scriptscriptstyle 6$};
		\node [style=none] (37) at (-2, \triht + \upmid + .2) {$\scriptscriptstyle 5^{*}$};
		\node [style=none] (38) at (-1.5, 2* \triht + \upmid + .3) {$\scriptscriptstyle 5$};
		\node [style=none] (39) at (-1.1, 2* \triht + \upmid + .4) {$\scriptscriptstyle 4$};
		\node [style=none] (40) at (-0.5, 2* \triht + \upmid + .2) {$\scriptscriptstyle 3^{*}$};
		\node [style=none] (41) at (0, 3* \triht + \upmid + .3) {$\scriptscriptstyle 3$};
		\node [style=none] (42) at (0.4, 3* \triht + \upmid + .4) {$\scriptscriptstyle 2$};
		\node [style=none] (43) at (1, 3* \triht + \upmid + .2) {$\scriptscriptstyle 1^{*}$};
		\node [style=none] (44) at (1.5, 4* \triht + \upmid + .3) {$\scriptscriptstyle 1$};
		\node [style=none] (45) at (1.9, 4* \triht + \upmid + .45) {$\scriptscriptstyle 0$};
		\node [style=none] (46) at (2.5, 4* \triht + \upmid + .3) {$\scriptscriptstyle \bar 1$};
		\node [style=none] (47) at (3, 3* \triht + \upmid + .2) {$\scriptscriptstyle \bar 1^{*}$};
		\node [style=none] (48) at (3.6, 3* \triht + \upmid + .4) {$\scriptscriptstyle \bar 2$};
		\node [style=none] (49) at (4, 3* \triht + \upmid + .3) {$\scriptscriptstyle \bar 3$};
		\node [style=none] (50) at (-3.2, \triht + .2) {$\scriptscriptstyle \bar 3$};
		\node [style=none] (51) at (-3.2, \upmid + .35) {$\scriptscriptstyle \bar 2$};
		\node [style=none] (52) at (-2.6, .2) {$\scriptscriptstyle \bar 1^{*}$};
		\node [style=none] (53) at (-3.2, -1* \triht + .2) {$\scriptscriptstyle \bar 1$};
		\node [style=none] (54) at (-3.2, -2 * \triht + \upmid + .35) {$\scriptscriptstyle 0$};
		\node [style=none] (55) at (-2.6, -2 * \triht + .2) {$\scriptscriptstyle 1$};
		\node [style=none] (56) at (-1.6, -2 * \triht + .2) {$\scriptscriptstyle 1^{*}$};
		\node [style=none] (57) at (-1.5, -3 * \triht + .35) {$\scriptscriptstyle 2$};
		\node [style=none] (58) at (-1.1, -3 * \triht + .2) {$\scriptscriptstyle 3$};
		\node [style=none] (59) at (-0.1, -3 * \triht + .2) {$\scriptscriptstyle 3^{*}$};
		\node [style=none] (60) at (0, -4 * \triht + .4) {$\scriptscriptstyle 4$};
		\node [style=none] (61) at (0.3, -4 * \triht + .2) {$\scriptscriptstyle 5$};
		\node [style=none] (62) at (0.7, -4 * \triht + .2) {$\scriptscriptstyle \bar 5$};
		\node [style=none] (63) at (1, -4 * \triht + .4) {$\scriptscriptstyle \bar 4$};
		\node [style=none] (64) at (1.2, -3 * \triht + .2) {$\scriptscriptstyle \bar 3^{*}$};
		\node [style=none] (65) at (2.1, -3 * \triht + .2) {$\scriptscriptstyle \bar 3$};
		\node [style=none] (66) at (2.5, -3 * \triht + .35) {$\scriptscriptstyle \bar 2$};
		\node [style=none] (67) at (2.7, -2 * \triht + .2) {$\scriptscriptstyle \bar 1^{*}$};
		\node [style=none] (68) at (3.7, -2 * \triht + .2) {$\scriptscriptstyle \bar 1$};
		\node [style=none] (69) at (4.2, -2 * \triht + \upmid + .35) {$\scriptscriptstyle 0$};
		\node [style=none] (70) at (4.2, -1* \triht + .2) {$\scriptscriptstyle 1$};
		\node [style=none] (71) at (3.7, .2) {$\scriptscriptstyle 1^{*}$};
		\node [style=none] (72) at (4.2, \upmid + .35) {$\scriptscriptstyle 2$};
		\node [style=none] (73) at (4.2, 1* \triht + .2) {$\scriptscriptstyle 3$};
		\node [style=none] (74) at (3.7, 2 * \triht + .2) {$\scriptscriptstyle 3^{*}$};
		\node [style=none] (75) at (4.2, 2 * \triht + \upmid + .35) {$\scriptscriptstyle 4$};
		\node [style=none] (76) at (4.2, 3* \triht + .2) {$\scriptscriptstyle 5$};
	\end{pgfonlayer}
	\begin{pgfonlayer}{edgelayer}
		\draw (13) to (14);
		\draw (19) to (13);
		\draw (20) to (19);
		\draw (21) to (20);
		\draw (14) to (22);
		\draw (22) to (21);
		\draw (23) to (24);
		\draw (24) to (25);
		\draw (25) to (26);
		\draw (26) to (27);
		\draw (27) to (28);
		\draw (28) to (23);
		\draw [style=dashed line] (23) to (29);
		\draw [style=dashed line] (29) to (34);
		\draw [style=dashed line] (30) to (24);
		\draw [style=dashed line] (25) to (31);
		\draw [style=dashed line] (30) to (31);
		\draw [style=dashed line] (32) to (26);
		\draw [style=dashed line] (27) to (33);
		\draw [style=dashed line] (32) to (33);
		\draw [style=dashed line] (34) to (28);
	\end{pgfonlayer}
	

    \clip       (-\cols, -\height) rectangle (\cols, \height);
    \draw[gray] (-\cols, -\height) rectangle (\cols, \height);

    \pgfmathsetmacro{\from}{-2 *\cols}
    \pgfmathsetmacro{\to}{2 * \cols}
    \foreach\i in {\from, ..., \to} {
        \draw[xslant=\slant]  (\i, -\height) -- (\i, \height);
        \draw[xslant=-\slant] (\i, -\height) -- (\i, \height);
    }

    \foreach\j in {-\rows, ..., \rows} {
        \pgfmathsetmacro{\y}{0.5 * \j * tan(60)}
        \draw (-\cols, \y) -- (\cols, \y);
    }

\end{tikzpicture}

\caption{Coordinates for Moves 1 and 2, odd chamber.}  \label{fig:moves1and2odd}
\end{figure}

{\bf Odd chambers}: Assume now that $w, w'$ and $w''$ belong to chamber I. We will again take $t$ to be such that $w''$ is outward and to the right from $w'$, when viewed from the origin. Then the NE, W, and SE edges of $\ch_{w'}$ and of $\ch_{w''}$ are collinear, while the E, NW, and SW edges of $\ch_{w''}$ are one root string farther out than the corresponding edges of $\ch_{w'}$.  As in the even chambers, $q^{w''}_{\bullet}$ is constant on each small outer $W(A_{2})$ hexagon comprising its 0-, 1-, and 2-shells, and each of these hexagons contains an alcove on the boundary of $\ch_{w}$. So it suffices to verify the theorem for $x$ on the 0-shell of $\ch_{w}$. 

In the hexagon $\ch_{w}$, the diagonal from the S vertex intersects the NW edge at the alcove labeled 2, and the diagonal from the NW vertex intersects the SE edge at the alcove labeled $\bar 2$. We find that, for an Up alcove in $\ch_{w}$ with coordinates $(i_{0}, i_{1}, i_{2})$,
\begin{equation} \label{e:Move2OddDeltaq}
\gD q_{2} := q_{x}^{w''} - q_{x}^{w'} = \gd_{i_{0},o>0} + \gd_{i_{0},*>0} + \gd_{i_{1},o>\bar 2} + \gd_{i_{1},*>0} + \gd_{i_{2},o>0} + \gd_{i_{2},*>2} - 1.
\end{equation}

Consider $x$ an Up alcove on the 0-shell of $\ch_{w}$. As observed in Move 1, the coordinates of $x$ are some permutation of $(o,o,e)$ and $q_{x}^{w'}=1$ or 0 (depending on whether or not $x$ is in a special segment). In fact, one odd coordinate is positive and the other is negative, so (except when $i_{1}=\bar 1$), $\gD q_{2} = 0$ and $q_{x}^{w''}=q_{x}^{w'}$. The Up alcove $x$ where $i_{1}=\bar 1$, with $\gD q_{2} = 1$, lies on the SE edge of $\ch_{w}$, just outside the special segment for $w$, whence $q_{x}^{w}=q_{x}^{w'}=0$ and $q_{x}^{w''}=1$. Since the $\ga_{1}$ diagonal which determines the upper endpoint of the SE special segment moves up-and-right one string from $w'$ to $w''$, the alcove $x$ lies in the small $W(A_{2})$ hexagon determined by the new alcove $y$ at the upper right of the longer SE special segment for $w''$, confirming that $q_{y}^{w''}=1$, as expected. 

For Down alcoves, the analog of \eqref{e:Move2OddDeltaq} is
\begin{equation} \label{e:Move2OddDeltaqDown}
\gD q_{2} := q_{x}^{w''} - q_{x}^{w'} = \gd_{i_{0},o>0} + \gd_{i_{0},*>0} + \gd_{i_{1},*>\bar 2} + \gd_{i_{1},o>0} + \gd_{i_{2},*>0} + \gd_{i_{2},o>2} - 1.
\end{equation}
We again leave it to the reader to check that this implies the formulas in the theorem statement for Down alcoves on the 0-shell of $\ch_{w}$, as required. In summary, for every $x$ on the 0-, 1-, or 2-shell of $\ch_{w''}$, the value $q_{x}^{w''}$ is as claimed. 
\end{proof}

\begin{Rem} \label{r:qpositive}
One can use the same reflection Move 1 and Move 2 formulas to prove inductively that $q^{w}_{x}>0$ for all $x$ on and inside the 3-shell of $\ch_{w}$ whenever $w$ is not spiral. This immediately implies that the Lookup Conjecture is (trivially) true for non-spiral elements in $\widetilde A_{2}$. In fact, this was our original approach, until we noticed that the Translation Move could be used to give more detailed information on the precise values of~$q$.
\end{Rem}


\section{Base cases} \label{s:basecases}
In this section we analyze the base cases for $w$.  By Remark \ref{r:base-rs}, if $X_w$ is a non-spiral rationally smooth Schubert variety, then $w$ is necessarily a base case, so as a consequence of our analysis, we obtain a description of the rationally Schubert varieties in type $\tilde{A}_2$ (Corollary \ref{c:ratsmoothSchubertvar}).

There are four types of base cases for $w$ in a chamber $\cc$: either $\cc$ is even or odd; the alcove $wA_{\circ}$ lies in the root strip adjacent to a fundamental root strip bounding $\cc$; and $w A_{\circ}$ shares either an edge or a vertex with a spiral (fundamental root strip) alcove. Recall also: either there are two simple reflections $s$ such that $ws < w$,
in which case $w$ is of Type 1, or there is one simple reflection $s$ such that $ws < w$, in which case $w$ is of Type 2.

{\em Case 1}: $\cc$ is even, and the alcove for $w$ shares an edge with a fundamental strip.  In this case, $w$ is a non-spiral element of the form $zs$, where $z$ is an {\em even}-length spiral element, $s$ is a simple reflection, and $\ell(w)=\ell(z)+1$.  Such an element is necessarily Type 1.  Following Billey and Crites \cite[Section 2.5]{BiCr:12}, we call such an element $w$ \emph{twisted spiral}; note that by definition, a twisted spiral element has odd length.  We claim that $q_{x}^{w}=0$ for all $x\le w$.  In fact, this follows from the work of Billey and Crites: they proved that the Schubert varieties corresponding to twisted spiral elements are rationally smooth at every point (see \cite[Theorem 1]{BiCr:12}), so the Carrell-Peterson criterion implies the claim.  We will give a different proof using the methods of this paper.

The proof is by induction on $\ell(w)$ using translation Moves 1 and 2, as in Theorem \ref{t:outershells}. The smallest value of $\ell(w)$ which can occur is $\ell(w) = 3$, and one can check directly that the claim holds in this case.  Now assume, as in the proof of Theorem \ref{t:outershells}, that $w$ is in chamber IV, with $w < ws =: w' < w't =: w''$ as before. Then $w$ and $w''$ are both twisted spiral. We are assuming the result for $w$, and wish to prove it for $w''$. We already know, by Theorem \ref{t:outershells}, that $q_{x}^{w''}=0$ for $x$ on the 0-, 1-, or 2-shell of $\ch_{w''}$. For $x$ belonging to an Up alcove in the interior of $\ch_{w}$, the coordinates of $x$ are some permutation of even $>0$, odd $>0$, and odd* $>0$. Thus by \eqref{e:Move1deltaq} and \eqref{e:Move2EvenDeltaq}, $\gD q_{1} = 1$ and $\gD q_{2} = -1$. Since by induction $q_{x}^{w} = 0$, we have $q_{x}^{w''} = 0$. As usual, we leave it to the reader to check the result for Down alcoves.

\medskip

{\em Case 2}: $\cc$ is even, and the alcove for $w$ shares a vertex with a fundamental strip.  See Figure \ref{fig:fish}.  In this case,
$w$ is of Type 2, $\ell(w)$ is even, and $\ell(w) \geq 4$.  We know by Theorem \ref{t:outershells} that $q^w_x = 0$ for $x$ on the 0- and 1-shells.  We claim that $q^w_x = 1$ for all $x$ on or inside the 2-shell.  

We prove the claim as follows.  The $k$-shells for $k \ge 2$ coincide with the hexagon $\ch_{z}$ where $z := z_1$ shown in Figure \ref{fig:fish}.  Then $z$ is twisted spiral, so for all $x \le z$, $q^z_x = 0$. We will compute $\Delta q := q^{w}_{x} - q^{z}_{x}$. Notice that $l(w)=l(z)+3$. 

Label the alcoves $x$ in $\ch_{z}$ by coordinates $(i_{0},i_{1},i_{2})$ as described after the statement of Theorem \ref{t:outershells} (with $z$ in place of $x$). We need only count how many additional alcoves there are in $\ch_{w} \setminus \ch_{z}$ on the root strings through $x$ in each direction $\ga_{0}, \ga_{1}, \ga_{2}$.
 Fix $j,\  0\le j\le 2$. First, notice that if the coordinate $i_{j}=\bar 1$, meaning $x$ lies along one of the edges of $\ch_{z}$, then there are two additional alcoves, one of each orientation, at \emph{each} end of the $\ga_{j}$ root string through $x$. So this root string contributes $+2$ to $\Delta q$. The same statements are true if $i_{j}=o^{*}$. But if $i_{j} = o$ or $e$, then there is only one additional alcove at each end of the $\ga_{j}$ root string through $x$, one Up and one Down. So these strings contribute $+1$ to $\Delta q$.

One checks immediately that the coordinates of the alcoves along the edges of $\ch_{z}$ are some permutation of $(\bar 1, o, e)$, whereas the interior alcove coordinates are some permutation of $(o^{*}, o, e)$. Therefore in any case,
$$
\Delta q = 2+1+1-3 =1,
$$
which, together with the fact that $q^{z}_{x}=0$, proves the claim.

We remark that in the setting of Case 2, the $2$-shell is empty if and only if $\ell(w) = 4$, and then $q^w_x = 0$ for all $x \leq w$.

\medskip

\begin{figure}[htbp!]
%

\pgfmathsetmacro{\cols}{5}
\pgfmathsetmacro{\rows}{5}
\pgfmathsetmacro{\slant}{cot(60)}
\pgfmathsetmacro{\height}{0.5 * \rows * tan(60)}
\pgfmathsetmacro{\triht}{sin(60)}
\pgfmathsetmacro{\upmid}{0.25 * sec(30)}
\pgfmathsetmacro{\downmid}{\triht - \upmid}

\begin{tikzpicture}
    

\input{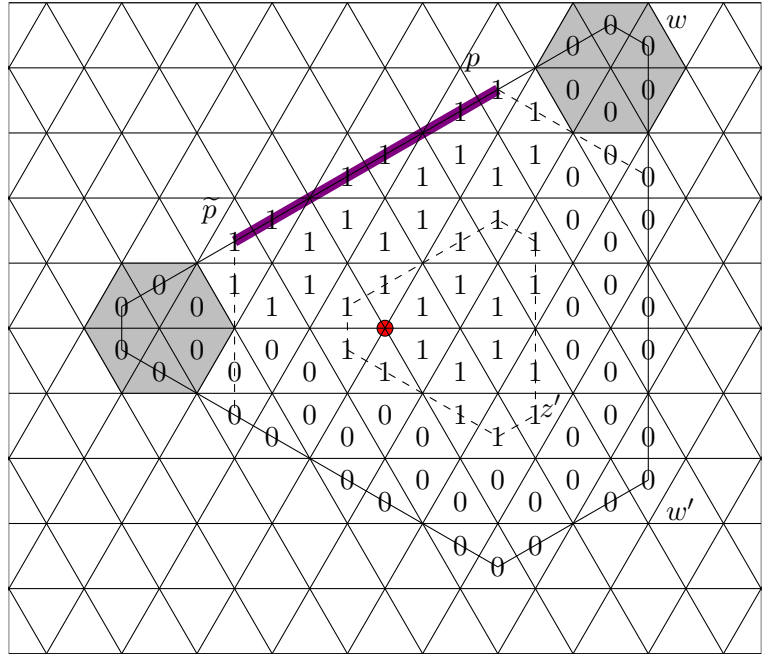}

	\begin{pgfonlayer}{nodelayer}
		\node [style=red dot] (4) at (0, 0) {};
		\node [style=none, label={above left:$p$}] (32) at (1.5, 3 * \triht + \downmid) {1};
		\node [style=none, label={above left:$\tilde p$}] (33) at (-2, 1 * \triht + \upmid) {1};
		\node [style=none] (37) at (3, 3 * \triht + \upmid) {0};
		\node [style=none] (38) at (3, -2 * \triht + \downmid) {0};
		\node [style=none] (39) at (-2.5, -1 * \triht + \downmid) {0};
		\node [style=none] (40) at (1, -3 * \triht + \upmid) {0};
		\node [style=none] (41) at (2.5, 3 * \triht + \downmid) {0};
		\node [style=none] (42) at (2.5, -3 * \triht + \downmid) {0};
		\node [style=none] (43) at (-2.5, 0 * \triht + \upmid) {0};
		\node [style=none] (44) at (2, -3 * \triht + \upmid) {0};
		\node [style=none] (45) at (-2, 0 * \triht + \downmid) {1};
		\node [style=none] (46) at (2, 3 * \triht + \upmid) {1};
		\node [style=none, label={above right:$w$}] (53) at (3.5, 4 * \triht + \upmid) {0};
		\node [style=none] (54) at (3.5, 3 * \triht + \downmid) {0};
		\node [style=none] (55) at (3.5, 2 * \triht + \upmid) {0};
		\node [style=none] (56) at (3.5, 1 * \triht + \downmid) {0};
		\node [style=none] (57) at (3.5, 0 * \triht + \upmid) {0};
		\node [style=none] (58) at (3.5, -1 * \triht + \downmid) {0};
		\node [style=none] (59) at (3.5, -2 * \triht + \upmid) {0};
		\node [style=none, label={below right:$w'$}] (60) at (3.5, -3 * \triht + \downmid) {0};
		\node [style=none] (61) at (3, 2 * \triht + \downmid) {0};
		\node [style=none] (62) at (3, 1 * \triht + \upmid) {0};
		\node [style=none] (63) at (3, 0 * \triht + \downmid) {0};
		\node [style=none] (64) at (3, -1 * \triht + \upmid) {0};
		\node [style=none] (65) at (2.5, 2 * \triht + \upmid) {0};
		\node [style=none] (66) at (2.5, 1 * \triht + \downmid) {0};
		\node [style=none] (67) at (2.5, 0 * \triht + \upmid) {0};
		\node [style=none] (68) at (2.5, -1 * \triht + \downmid) {0};
		\node [style=none] (69) at (2.5, -2 * \triht + \upmid) {0};
		\node [style=none] (70) at (3, 4 * \triht + \downmid) {0};
		\node [style=none] (72) at (2, 2 * \triht + \downmid) {1};
		\node [style=none] (73) at (2, 1 * \triht + \upmid) {1};
		\node [style=none] (74) at (2, 0 * \triht + \downmid) {1};
		\node [style=none] (75) at (2, -1 * \triht + \upmid) {1};
		\node [style=none] (76) at (2, -2 * \triht + \downmid) {1};
		\node [style=none] (77) at (1.5, 2 * \triht + \upmid) {1};
		\node [style=none] (78) at (1.5, 1 * \triht + \downmid) {1};
		\node [style=none] (79) at (1.5, 0 * \triht + \upmid) {1};
		\node [style=none] (80) at (1.5, -1 * \triht + \downmid) {1};
		\node [style=none] (81) at (1, 2 * \triht + \downmid) {1};
		\node [style=none] (82) at (1, 1 * \triht + \upmid) {1};
		\node [style=none] (83) at (1, 0 * \triht + \downmid) {1};
		\node [style=none] (84) at (1, -1 * \triht + \upmid) {1};
		\node [style=none] (85) at (1, 3 * \triht + \upmid) {1};
		\node [style=none] (87) at (0.5, 2 * \triht + \upmid) {1};
		\node [style=none] (88) at (0.5, 1 * \triht + \downmid) {1};
		\node [style=none] (89) at (0.5, 0 * \triht + \upmid) {1};
		\node [style=none] (90) at (0.5, -1 * \triht + \downmid) {1};
		\node [style=none] (92) at (0, 2 * \triht + \downmid) {1};
		\node [style=none] (93) at (0, 1 * \triht + \upmid) {1};
		\node [style=none] (94) at (0, 0 * \triht + \downmid) {1};
		\node [style=none] (95) at (-0.5, 2 * \triht + \upmid) {1};
		\node [style=none] (96) at (-0.5, 1 * \triht + \downmid) {1};
		\node [style=none] (97) at (-0.5, 0 * \triht + \upmid) {1};
		\node [style=none] (99) at (-1, 1 * \triht + \upmid) {1};
		\node [style=none] (100) at (-1, 0 * \triht + \downmid) {1};
		\node [style=none] (102) at (-1.5, 1 * \triht + \downmid) {1};
		\node [style=none] (104) at (-3, 0 * \triht + \downmid) {0};
		\node [style=none] (105) at (-3, -1 * \triht + \upmid) {0};
		\node [style=none] (107) at (-3.5, 0 * \triht + \upmid) {0};
		\node [style=none] (108) at (-3.5, -1 * \triht + \downmid) {0};
		\node [style=none] (109) at (-2, -1 * \triht + \upmid) {0};
		\node [style=none] (110) at (-2, -2 * \triht + \downmid) {0};
		\node [style=none] (112) at (-1.5, 0 * \triht + \upmid) {1};
		\node [style=none] (113) at (-1.5, -1 * \triht + \downmid) {0};
		\node [style=none] (114) at (-1.5, -2 * \triht + \upmid) {0};
		\node [style=none] (116) at (-1, -1 * \triht + \upmid) {0};
		\node [style=none] (117) at (-1, -2 * \triht + \downmid) {0};
		\node [style=none] (119) at (-0.5, -1 * \triht + \downmid) {1};
		\node [style=none] (120) at (-0.5, -2 * \triht + \upmid) {0};
		\node [style=none] (121) at (-0.5, -3 * \triht + \downmid) {0};
		\node [style=none] (122) at (0, -1 * \triht + \upmid) {1};
		\node [style=none] (123) at (0, -2 * \triht + \downmid) {0};
		\node [style=none] (124) at (0, -3 * \triht + \upmid) {0};
		\node [style=none] (126) at (0.5, -2 * \triht + \upmid) {0};
		\node [style=none] (127) at (0.5, -3 * \triht + \downmid) {0};
		\node [style=none] (129) at (1, -2 * \triht + \downmid) {1};
		\node [style=none] (130) at (1, -4 * \triht + \downmid) {0};
		\node [style=none] (131) at (1.5, -2 * \triht + \upmid) {1};
		\node [style=none] (132) at (1.5, -3 * \triht + \downmid) {0};
		\node [style=none] (133) at (1.5, -4 * \triht + \upmid) {0};
		\node [style=none] (134) at (2, -4 * \triht + \downmid) {0};
		\node [style=none] (136) at (3, -3 * \triht + \upmid) {0};
		\node [style=none] (138) at (2.5, 4 * \triht + \upmid) {0};
		\node [style=none] (139) at (2.2, -2 * \triht + 0.7) {$z'$};
		\node [style=none] (139) at (3.5, 5 * \triht) {};
		\node [style=none] (140) at (2.5, 5 * \triht) {};
		\node [style=none] (141) at (2, 4 * \triht) {};
		\node [style=none] (142) at (2.5, 3 * \triht) {};
		\node [style=none] (143) at (3.5, 3 * \triht) {};
		\node [style=none] (144) at (4, 4 * \triht) {};
		\node [style=none] (145) at (-2.5, 1 * \triht) {};
		\node [style=none] (146) at (-3.5, 1 * \triht) {};
		\node [style=none] (147) at (-4, 0) {};
		\node [style=none] (148) at (-3.5, -1 * \triht) {};
		\node [style=none] (149) at (-2.5, -1 * \triht) {};
		\node [style=none] (150) at (-2, 0) {};
	\end{pgfonlayer}
	\begin{pgfonlayer}{edgelayer}
		\draw [style=filled grey] (140.center)
			 to (141.center)
			 to (142.center)
			 to (143.center)
			 to (144.center)
			 to (139.center)
			 to cycle;
		\draw [style=filled grey] (150.center)
			 to (145.center)
			 to (146.center)
			 to (147.center)
			 to (148.center)
			 to (149.center)
			 to cycle;
		\draw [style=thick purple line] (32.center) to (33.center);
		\draw [style=none] (53.center) to (70.center);
		\draw [style=none] (70.center) to (107.center);
		\draw [style=none] (107.center) to (108.center);
		\draw [style=none] (108.center) to (133.center);
		\draw [style=none] (133.center) to (60.center);
		\draw [style=none] (60.center) to (53.center);
		\draw [style=dashed line] (76.center) to (73.center);
		\draw [style=dashed line] (73.center) to (78.center);
		\draw [style=dashed line] (78.center) to (97.center);
		\draw [style=dashed line] (97.center) to (119.center);
		\draw [style=dashed line] (119.center) to (131.center);
		\draw [style=dashed line] (131.center) to (76.center);
		\draw [style=dashed line] (32.center) to (55.center);
		\draw [style=dashed line] (33.center) to (110.center);
	\end{pgfonlayer}
	

    \clip       (-\cols, -\height) rectangle (\cols, \height);
    \draw[gray] (-\cols, -\height) rectangle (\cols, \height);

    \pgfmathsetmacro{\from}{-2 *\cols}
    \pgfmathsetmacro{\to}{2 * \cols}
    \foreach\i in {\from, ..., \to} {
        \draw[xslant=\slant]  (\i, -\height) -- (\i, \height);
        \draw[xslant=-\slant] (\i, -\height) -- (\i, \height);
    }

    \foreach\j in {-\rows, ..., \rows} {
        \pgfmathsetmacro{\y}{0.5 * \j * tan(60)}
        \draw (-\cols, \y) -- (\cols, \y);
    }

\end{tikzpicture}
\caption{Base Case 3}  \label{fig:basecase3}
\end{figure}

{\em Case 3}: $\cc$ is odd, and the alcove for $w$ shares an edge with a fundamental strip.  In this case, $w$ is a non-spiral element of the form $zs$, where $z$ is an {\em odd}-length spiral element, $s$ is a simple reflection, and $\ell(w)=\ell(z)+1$.  See Figure \ref{fig:basecase3}. In this case, $w$ is of Type 1, $\ell(w)$ is even,
and $\ell(w) \geq 4$.
If $\ell(w) = 4$, there are no special segments, the $0$-, $1$-, and $2$-shells exhaust $\ch_w$, and one can check that $q^w_x = 0$ for all $x \leq w$.  

Assume now that 
$\ell(w) \geq 6$.  In this case, there is exactly one special segment along an outer edge of $\ch_w$. The right $R(w)$-orbits are little
$A_2$ hexagons, and $q^w_x$ is the same for all
$x$ in a given $R(w)$-orbit.  Suppose $x \le w$.  We claim that $q_{x}^{w}=1$, except that $q_{x}^{w}=0$ if $x$ is on the 0-, 1-, or 2-shell of $\ch_w$ and
$x R(w)$ does not intersect the special segment of $\ch_{w}$.  Observe that if $x$ is on the 0-, 1-, or 2-shell of $\ch_w$, then the claim follows directly from Theorem \ref{t:outershells}.  Therefore assume $x$ is on a $k$-shell for $k \ge 3$.  We will show (see Theorem \ref{t:qHeredity} case I.1.(ii)) that the $k$-shells for $k \ge 3$ coincide with $\ch_{z'}$ where $z'$ is the element shown in Figures \ref{fig:basecase3} and \ref{fig:edge3shell}.  Since $x$ is in $\ch_{z'}$ and $z'$ is twisted spiral,
we have $q^{z'}_x = 0$. Translation into the chamber takes $z'$ to $w'$. By Theorem \ref{t:translationmove}, $q^{w'}_x = q^{z'}_x + 2 = 0 + 2 = 2$. 
 
To compute $q^w_x$, note that the hexagon $\ch_{w'}$ is $\ch_w$ with the two short corners truncated: cut off the small grey hexagons in the NE and W corners (when $w$ is as drawn in Figure \ref{fig:edge3shell}) to obtain a semiregular hexagon $\ch_{w'}$ all of whose short edges have length 4. The fact that $x$ is inside the
little hexagon $\ch_{z'}$ implies that \emph{none} of the root strings through $x$ contain \emph{any} grey alcoves in the portion of $\ch_w$ that was cut off to form $\ch_{w'}$. By Lemma \ref{l:lengthChange}, $\ell(w')=\ell(w)-1$. 
Thus, $q^w_x = q^{w'}_x - 1 = 2 - 1 = 1$.

\medskip

\begin{figure}[htbp!]
\include{Fig-basecase4}
\caption{Base Case 4}  \label{fig:basecase4}
\end{figure}

{\em Case 4}: $\cc$ is odd, and the alcove for $w$ shares a vertex with a fundamental strip.  In this case, $w$ is of Type 2, so the right $R(w)$-orbits are little diamonds.  
In this case, $\ell(w)$ is odd, and $\ell(w) \geq 5$.  If $\ell(w) = 5$, then there are no special segments, and one can check that if $x$ is on the $0$-shell or the $1$-shell, $q^w_x = 0$;
if $x$ is one of the two elements on the $2$-shell, then $q^w_x = 2$.  

Assume now that $\ell(w) \geq 7$. There is one special segment on the 0-shell.  Suppose $x \leq w$.  By Theorem \ref{t:outershells}, if $x$ is on the 0-shell or 1-shell of $\ch_w$, then $q^w_x = 0$ if and only if $R(w)$-orbit of $x$ does not intersect the special segment.  If $x$ is on the 0- or 1- shell and $x R(w)$ intersects the special segment (shown in purple in Figure \ref{fig:basecase4}), then $q^w_x = 1$.

If $x$ is not in the 0-shell or 1-shell of $\ch_w$, then $q^w_x$ is either $1$ or $2$.  For such $x$, there is an (unfortunately somewhat complicated) explicit description of the values of $q^w_x$.  We begin with some geometry.  Label the vertices of $\ch_w$ as follows.  Set $u_0 = w$, and let $u_1$ denote the vertex of the special edge for which $u_{0}u_{1}$ is an edge.  Moving from $u_0$ to $u_1$ determines a direction---either clockwise or counterclockwise---and continuing in that direction, label the remaining vertices $u_2, u_3, u_4$ and $u_5$, in order.  If the direction is counterclockwise, then for $i \geq 1$, $u_i = w_i$; otherwise, $u_i = w_{6-i}$. See Figure \ref{fig:basecase4}, in which $u_{i}=w_{i}$ for all $i$.

Let $D_i$ be the part of the diagonal through $u_i$ which
is not on the 0- or 1-shells.  Note that $D_0 = D_3$.  Among the set of alcove centers $xq$ in $\ch_w$, the complement of the set of $xq$ on the
0- and 1-shells is the disjoint union of the 
$xq$ in the triangular region with vertices
$D_0 \cap D_4$, $D_0 \cap D_5$, and $D_4 \cap D_5$ (shown in red in Figure \ref{fig:basecase4}); the segments $D_1$ and $D_2$ (colored green); and the
$xq$ on four additional segments: two parallel segments between $D_1$ and $D_4$, 
and two between $D_2$ and $D_5$ (all four colored blue).  Two alcoves $A'$ and $A''$ near the $u_4 u_{5}$ edge each lie on two of these blue segments, and each blue segment has either $A'$ or $A''$ at one end.

Now we can explicitly describe $q^w_x$.  If $x$ is in the red triangular region, then $q^w_x = 2$.  If $x$ is on a green segment $D_1$ or $D_2$, then $q^w_x = 1$.  On each of the remaining four blue segments, $q^w_x$ is as follows.  If $x$ is on either end of the segment, then $q^w_x = 2$, and moving from one end of the segment to the other, $q^w_x$ alternates in the
pattern $2, 1, 2,  \dots, 1,2$.

To prove these assertions, recall that $w$ is of Type 2, so there is a unique simple reflection $s$ such that $w > ws =: w'$. Then $w'$ is of Type 1 and shares an edge with the same fundamental root strip with which $w$ shares a vertex; see Figure \ref{fig:basecase4} (where $\ch_{w'}$ is indicated by the dashed line). In particular, $w'$ is in Base Case 3, and $\ch_{w'}$ is the 1-shell of $\ch_{w}$. For each $x\le w'$ we can determine $q^{w}_{x}$ via Move 1 and the known values of $q^{w'}_{x}$ from Case 3, above. We use the earlier Move 1 formulas for $\Delta q_{1}$ in terms of coordinates $(i_{0},i_{1},i_{2})$ as in Figures \ref{fig:moves1and2even} and \ref{fig:moves1and2odd}, from \eqref{e:Move1deltaq} and the subsequent text (with the roles of $w$ and $w'$ reversed---in particular, the coordinates are based on the alcoves around the $\ch_{w'}$ hexagon edges). 

By the discussion in the second paragraph of Case 4, we may assume that $x$ is on or inside the 2-shell of $w$; i.e., the 1-shell of $w'$. Thus, $x$ is either on a green or blue root string, or in the red triangle. We treat each of these cases in turn. By symmetry, it suffices to analyze the NW-to-SE green and blue root strings where $i_{1}$ is constant, plus the red triangle.

On the green string where $i_{1}=\bar 2$, the alcove coordinates are either $(o, \bar 2, o^{*})$ or $(o^{*}, \bar 2, o)$, according to whether the alcove is Down or Up. This string does not contain any alcoves between two parallel diagonals, so the formula for $\Delta q$ is the same for Up and Down alcoves: $\Delta q = 1+1-1 = 1$. Hence $q^{w}_{x}=q^{w'}_{x} + \Delta q = 0+1=1$.

On the lower blue string, $i_{1}=\bar 1^{*}$. The Up alcove coordinates are $(e,\bar 1^{*},o^{*})$ so $\Delta q=1+1+1-1=2$ and $q^{w}_{x} = 0+2 = 2$. (This includes the bottom right alcove on the string, $A'$.) The Down alcove coordinates are $(o^{*},\bar 1^{*},e)$. This is a string between two parallel diagonals, so the middle coordinate $\bar 1^{*}$ does not contribute to an increase in $q$. In fact, $\Delta q = 1+1-1 = 1$ and $q^{w}_{x} = 0+1 = 1$.

On the upper blue string, $i_{1}=\bar 1$. The Up alcove coordinates are $(o,\bar 1,e)$, $\Delta q = 1-1 = 0$, and $q^{w}_{x} = 1 + 0 = 1$. The Down alcoves other than the bottom right one have coordinates $(e,\bar 1,o)$. This again is a string between two parallel diagonals, so $\Delta q = 1+1-1 = 1$, and $q^{w}_{x}=1+1=2$. The bottom right alcove, $A''$, has coordinates $(e, \bar 1, 1)$. This alcove lies between \emph{both} pairs of parallel diagonals, so $\Delta q = 1+1+1-1 = 2$, whence $q^{w}_{x} = 0+2=2$.

Finally, assume $x$ lies on or inside the red triangle. One checks that the coordinates of $x$ are some permutation of $(o,o^{*},e)$. Since this region avoids the strings between parallel diagonals, the formulas for Up and Down alcoves are the same: $\Delta q = 1+1-1 = 1$ and $q^{w}_{x}=1+1=2$.

\begin{Cor} \label{c:ratsmoothSchubertvar}
The Schubert variety $X_{w}$ is rationally smooth if and only if one of the following holds:

\noindent (a) $\ell(w) \leq 3$.

\noindent (b) $w$ is in an even chamber, the alcove for $w$ shares a vertex with a fundamental strip, and $\ell(w) = 4$.  Equivalently, $w = s_i s_j s_i s_k$ for $\{i,j,k\} = \{0,1,2 \}$.

\noindent (c) $w$ is in an odd chamber, the alcove for $w$ shares an edge with a fundamental strip, and $\ell(w) = 4$.  Equivalently, $w = s_k s_i s_j s_i$ for $\{i,j,k\} = \{0,1,2 \}$.

\noindent (d) $w$ is twisted spiral.
\end{Cor}

\begin{proof}
If $w$ is spiral, then $X_{w}$ is rationally smooth if and only if $\ell(w) \leq 3$ by \cite{GrLi:15}.  Therefore, we may
 assume $w$ is non-spiral.  
 By Remark \ref{r:base-rs}, if $X_{w}$ is rationally smooth, then $X_{w}$ is a base case.  Our analysis of the base cases
 shows that $q^w_x = 0$ for all $x \leq w$ (which is equivalent to rational smoothness of $X_{w}$) in exactly the following cases: Base Case 1; 
 Base Case 2 with $\ell(w) = 4$; Base Case 3 with $\ell(w) = 4$.  
 These are exactly parts (d), (b), and (c) of the Corollary; part (a) of the corollary (for non-spiral $w$)
 corresponds to Base Case 1 with $\ell(w) = 3$.  The reduced expressions in (b) and (c) can be obtained by inspection.
 \end{proof}
 
 Note that conditions of the corollary are mutually exclusive, except that the twisted spiral element $w = s_i s_j s_i$ satisfies both (a) and (d) of the corollary.
We will see below that the Schubert varieties in cases (a)-(c) of the corollary are smooth (see Corollary \ref{c:ratsmoothSchubertvar}), so $X_w$ is rationally smooth if and only if $X_w$ is smooth or $w$ is twisted spiral.  This result was first proved by Billey and Crites \cite{BiCr:12}; see Remark \ref{r:BilleyCrites} for further discussion.


\section{The nrs locus for affine \texorpdfstring{$A_{2}$}{A\_2}} \label{s:nrs}
In this section we show that for a non-spiral Schubert variety $X_w$ of type $\tilde{A}_2$, the point $x \cb$ is nrs in $X_w$ if and only if $q^w_x>0$
(Corollary \ref{c:trivial-rs}).  This result has a number of consequences.  As noted in the introduction, combined with the work in \cite{GrLi:15} for spiral Schubert varieties, it implies the Lookup Conjecture for type $\tilde{A}_2$.  
Moreover, combined with our analysis of $q^w_x$ in previous sections, Corollary \ref{c:trivial-rs} yields
a description of the nrs locus of a non-spiral Schubert variety in terms of the geometry of the Bruhat hexagon $\ch_w$ (see Remark \ref{r:shell-rs})
as well as a description of the maximal nrs $z < w$ (see Corollary \ref{c:Maxnrs}).

The results of this section rely on the following key theorem.

\begin{Thm} \label{t:qHeredity}
Let $w\in W$ be non-spiral, and $y\le x\le w$. If $q_{x}^{w}>0$ then $q_{y}^{w}>0$.
\end{Thm}

\begin{proof}
The proof relies heavily on the characterization of the locus where $q^{w}_{\bullet}>0$ in Corollary \ref{c:qpositive}, and amounts to showing that this region contains the Bruhat hexagons of all its elements. 

First, if $w$ is a twisted spiral element, as in Case 1 in Section \ref{s:basecases}, then $q_{x}^{w}=0$ for all $x\le w$ and there is nothing to prove. So assume henceforth that $w$ is not twisted spiral.

Assume the hypothesis of the theorem. Our strategy is to identify a small number (at most four) of elements $z$ such that $x \le z \le w$, $q_{z}^{w}>0$, and the entire hexagon $\ch_{z}$ is contained in the region where $q^{w}_{\bullet}>0$. In particular, since $y\le z$, we deduce that $q_{y}^{w}>0$, as required. As a corollary, we will have identifed the maximal elements $z<w$ having $q_{z}^{w}>0$.

Let $t(\ga)$ be the translation into the chamber of $w$. 

\noindent I. We will treat first the cases where $q_{x}^{w}>0$ by virtue of part (a) or (b) of Corollary \ref{c:qpositive}.

\begin{figure}[htbp!]
%

\pgfmathsetmacro{\cols}{7}
\pgfmathsetmacro{\rows}{8}
\pgfmathsetmacro{\slant}{cot(60)}
\pgfmathsetmacro{\height}{0.5 * \rows * tan(60)}
\pgfmathsetmacro{\triht}{sin(60)}
\pgfmathsetmacro{\upmid}{0.25 * sec(30)}
\pgfmathsetmacro{\downmid}{\triht - \upmid}

\begin{tikzpicture}
    

\input{sample.tikzstyles}

	\begin{pgfonlayer}{nodelayer}
		\node [style=none] (2) at (-7, 1 * \triht) {};
		\node [style=none] (3) at (7, 1 * \triht) {};
		\node [style=red dot] (4) at (0, 0) {};
		\node [style=none, label={below left:$X$}] (5) at (5, 4 * \triht + \upmid) {};
		\node [style=none] (7) at (-3, -8 * \triht) {};
		\node [style=none] (8) at (5, 8 * \triht) {};
		\node [style=none] (11) at (5, -8 * \triht) {};
		\node [style=none] (12) at (-3, 8 * \triht) {};
		\node [style=small black dot, label={right:$z$}] (13) at (4, 3 * \triht + \upmid) {};
		\node [style=small black dot, label={above:$z_1$}] (14) at (2, 4 * \triht + \downmid) {};
		\node [style=small black dot, label={right:$z_5$}] (19) at (4, -2 * \triht + \downmid) {};
		\node [style=small black dot, label={below:$z_4$}] (20) at (0.5, -4 * \triht + \upmid) {};
		\node [style=small black dot, label={left:$z_3$}] (21) at (-3, -2 * \triht + \downmid) {};
		\node [style=small black dot, label={left:$z_2$}] (22) at (-3, 1 * \triht + \upmid) {};
		\node [style=small black dot, label={right:$w$}] (23) at (5.5, 4 * \triht + \upmid) {};
		\node [style=small black dot, label={above:$w_1$}] (24) at (2, 6 * \triht + \downmid) {};
		\node [style=small black dot, label={left:$w_2$}] (25) at (-4.5, 2 * \triht + \upmid) {};
		\node [style=small black dot, label={left:$w_3$}] (26) at (-4.5, -3 * \triht + \downmid) {};
		\node [style=small black dot, label={below:$w_4$}] (27) at (0.5, -6 * \triht + \upmid) {};
		\node [style=small black dot, label={right:$w_5$}] (28) at (5.5, -3 * \triht + \downmid) {};
		\node [style=none] (29) at (6, 4 * \triht) {};
		\node [style=none] (30) at (5.5, 5 * \triht) {};
		\node [style=none] (31) at (4.5, 5 * \triht) {};
		\node [style=none] (32) at (3.5, 3 * \triht) {};
		\node [style=none] (33) at (5.5, 3 * \triht) {};
		\node [style=none] (34) at (2.5, 7 * \triht) {};
		\node [style=none] (35) at (1.5, 7 * \triht) {};
		\node [style=none] (36) at (1, 6 * \triht) {};
		\node [style=none] (37) at (2, 4 * \triht) {};
		\node [style=none] (38) at (3, 6 * \triht) {};
		\node [style=none] (39) at (-4.5, 3 * \triht) {};
		\node [style=none] (40) at (-5, 2 * \triht) {};
		\node [style=none] (41) at (-4.5, 1 * \triht) {};
		\node [style=none] (42) at (-2.5, 1 * \triht) {};
		\node [style=none] (43) at (-3.5, 3 * \triht) {};
		\node [style=none] (44) at (-5, -2 * \triht) {};
		\node [style=none] (45) at (-4.5, -3 * \triht) {};
		\node [style=none] (46) at (-3.5, -3 * \triht) {};
		\node [style=none] (47) at (-2.5, -1 * \triht) {};
		\node [style=none] (48) at (-4.5, -1 * \triht) {};
		\node [style=none] (49) at (0, -6 * \triht) {};
		\node [style=none] (50) at (1, -6 * \triht) {};
		\node [style=none] (51) at (1.5, -5 * \triht) {};
		\node [style=none] (52) at (0.5, -3 * \triht) {};
		\node [style=none] (53) at (-0.5, -5 * \triht) {};
		\node [style=none] (54) at (5.5, -3 * \triht) {};
		\node [style=none] (55) at (6, -2 * \triht) {};
		\node [style=none] (56) at (5.5, -1 * \triht) {};
		\node [style=none] (57) at (3.5, -1 * \triht) {};
		\node [style=none] (58) at (4.5, -3 * \triht) {};
		
		\node [style=none] (60) at (5, 4 * \triht) {};
		\node [style=none] (61) at (2, 6 * \triht) {};
		\node [style=none] (62) at (-4, 2 * \triht) {};
		\node [style=none] (63) at (-4, -2 * \triht) {};
		\node [style=none] (64) at (0.5, -5 * \triht) {};
		\node [style=none] (65) at (5, -2 * \triht) {};
		
		\node [style=none] (66) at (4.5, 3 * \triht + \downmid) {};
		\node [style=none] (67) at (2, 5 * \triht + \upmid) {};
		\node [style=none] (68) at (-3.5, 1 * \triht + \downmid) {};
		\node [style=none] (69) at (-3.5, -2 * \triht + \upmid) {};
		\node [style=none] (70) at (0.5, -5 * \triht + \downmid) {};
		\node [style=none] (71) at (4.5, -2 * \triht + \upmid) {};
	\end{pgfonlayer}
	\begin{pgfonlayer}{edgelayer}
		\draw [style=red line] (2.center) to (3.center);
		\draw [style=green line] (7.center) to (8.center);
		\draw [style=blue line] (12.center) to (11.center);
		\draw [style=filled grey] (33.center)
			 to (29.center)
			 to (30.center)
			 to (31.center)
			 to (32.center)
			 to cycle;
		\draw [style=filled grey] (36.center)
			 to (37.center)
			 to (38.center)
			 to (34.center)
			 to (35.center)
			 to cycle;
		\draw [style=filled grey] (43.center)
			 to (39.center)
			 to (40.center)
			 to (41.center)
			 to (42.center)
			 to cycle;
		\draw [style=filled grey] (47.center)
			 to (48.center)
			 to (44.center)
			 to (45.center)
			 to (46.center)
			 to cycle;
		\draw [style=filled grey] (52.center)
			 to (53.center)
			 to (49.center)
			 to (50.center)
			 to (51.center)
			 to cycle;
		\draw [style=filled grey] (58.center)
			 to (54.center)
			 to (55.center)
			 to (56.center)
			 to (57.center)
			 to cycle;
		\draw [style=dashed line] (13) to (14);
		\draw [style=dashed line] (19) to (13);
		\draw [style=dashed line] (20) to (19);
		\draw [style=dashed line] (21) to (20);
		\draw [style=dashed line] (14) to (22);
		\draw [style=dashed line] (22) to (21);
		\draw (23) to (24);
		\draw (24) to (25);
		\draw (25) to (26);
		\draw (26) to (27);
		\draw (27) to (28);
		\draw (28) to (23);
		\draw [style=dotted line] (60) to (61);
		\draw [style=dotted line] (61) to (62);
		\draw [style=dotted line] (62) to (63);
		\draw [style=dotted line] (63) to (64);
		\draw [style=dotted line] (64) to (65);
		\draw [style=dotted line] (65) to (60);
		
		\draw [style=dotted line] (66) to (67);
		\draw [style=dotted line] (67) to (68);
		\draw [style=dotted line] (68) to (69);
		\draw [style=dotted line] (69) to (70);
		\draw [style=dotted line] (70) to (71);
		\draw [style=dotted line] (71) to (66);
	\end{pgfonlayer}


    \clip       (-\cols, -\height) rectangle (\cols, \height);
    \draw[gray] (-\cols, -\height) rectangle (\cols, \height);

    \pgfmathsetmacro{\from}{-2 *\cols}
    \pgfmathsetmacro{\to}{2 * \cols}
    \foreach\i in {\from, ..., \to} {
        \draw[xslant=\slant]  (\i, -\height) -- (\i, \height);
        \draw[xslant=-\slant] (\i, -\height) -- (\i, \height);
    }

    \foreach\j in {-\rows, ..., \rows} {
        \pgfmathsetmacro{\y}{0.5 * \j * tan(60)}
        \draw (-\cols, \y) -- (\cols, \y);
    }

\end{tikzpicture}
\caption{3-shell of $\ch_{w}$: $w$ of Type 1}  \label{fig:3shell}
\end{figure}

\noindent
1. Consider first $w$ of Type $\tau=1$. We need to analyze the elements on or inside the $4-\tau = 3$-shell of $\ch_{w}$. Label the simple reflections $s, t, u$ so that $w<ws$, $w>wt$, and $w>wu$. The alcoves associated to the six elements $wr$ for $r \in R(w) = \langle t, u \rangle$ form a small regular hexagon $H$. 

(i) Assume first that $H$ lies entirely within the chamber containing $w$. Adjoin to $H$ the additional alcove corresponding to $z:=wtuts=wutus = t(-\ga) w$, and call the resulting region $X$. (Its shape is the unique convex pentagonal ``heptiamond.'') We use $X$ to help describe the $3$-shell of $\ch_w$, as follows.  Reflect $X$ across the three lines $H_{\ga,i}$,  $H_{\gb,j}$, and $H_{\gg,k}$ in the same fashion as was used to define the vertices of $\ch_{w}$ in \eqref{e:hexagonVertices}, forming five additional copies of $X$; see Figure \ref{fig:3shell}. It is clear that the hexagon $\ch_{z}$ coincides with the 3-shell of $\ch_{w}$; indeed, $w=t(\ga)z$. By Corollary~\ref{c:qpositive}, $q_{y}^{w}>0$ for all $y\le z$. So if $q_{x}^{w}>0$ ``because'' $x$ lies on or inside the 3-shell of $w$, and if $y\le x$, then $y\le x\le z$ implies $q_{y}^{w}>0$.

\begin{figure}[htbp!]
%

\pgfmathsetmacro{\cols}{6}
\pgfmathsetmacro{\rows}{7}
\pgfmathsetmacro{\slant}{cot(60)}
\pgfmathsetmacro{\height}{0.5 * \rows * tan(60)}
\pgfmathsetmacro{\triht}{sin(60)}
\pgfmathsetmacro{\upmid}{0.25 * sec(30)}
\pgfmathsetmacro{\downmid}{\triht - \upmid}

\begin{tikzpicture}
    

\input{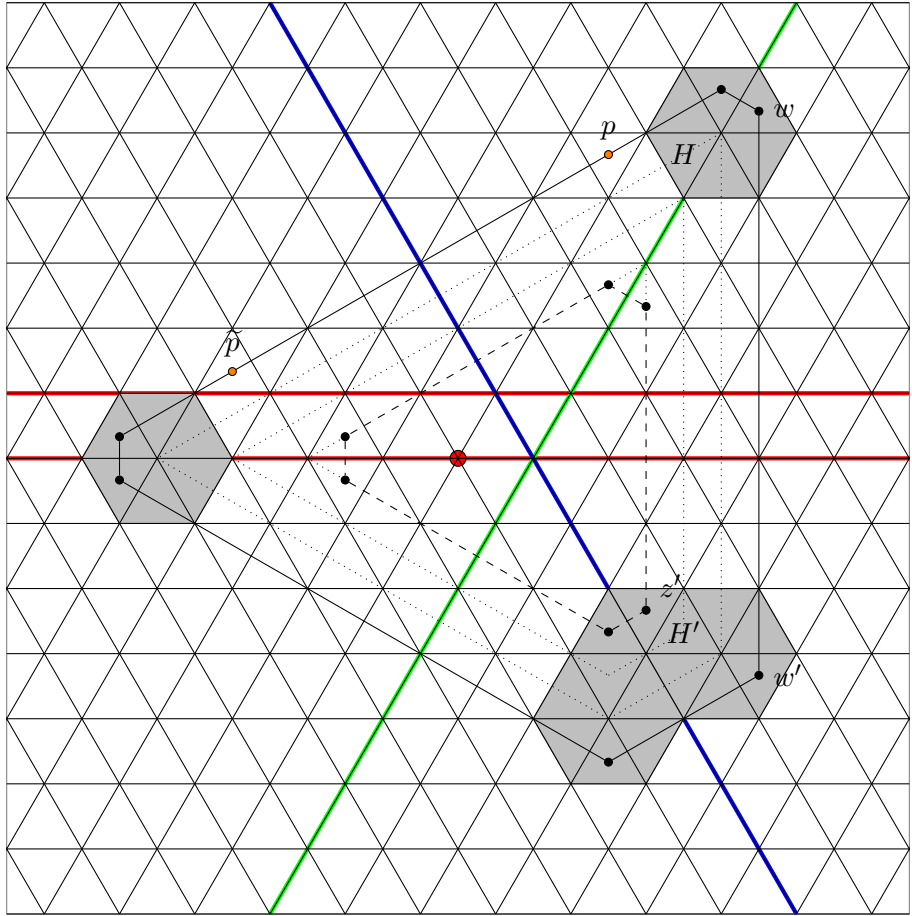}

	\begin{pgfonlayer}{nodelayer}
		\node [style=none] (2) at (-6, 1 * \triht) {};
		\node [style=none] (3) at (6, 1 * \triht) {};
		\node [style=red dot] (4) at (0, 0) {};
		\node [style=none] (5) at (-6, 0) {};
		\node [style=none] (6) at (6, 0) {};
		\node [style=none] (7) at (-2.5, -7 * \triht) {};
		\node [style=none] (8) at (4.5, 7 * \triht) {};
		\node [style=none] (9) at (3, 4 * \triht + \downmid) {$H$};
		\node [style=none] (10) at (3, -3 * \triht + \upmid) {$H'$};
		\node [style=none] (11) at (4.5, -7 * \triht) {};
		\node [style=none] (12) at (-2.5, 7 * \triht) {};
		\node [style=small black dot] (13) at (2.5, 2 * \triht + \upmid) {};
		\node [style=small black dot] (14) at (2, 2 * \triht + \downmid) {};
		\node [style=small black dot, label={above right:$z'$}] (19) at (2.5, -3 * \triht + \downmid) {};
		\node [style=small black dot] (20) at (2, -3 * \triht + \upmid) {};
		\node [style=small black dot] (21) at (-1.5, -1 * \triht + \downmid) {};
		\node [style=small black dot] (22) at (-1.5, \upmid) {};
		\node [style=small black dot, label={right:$w$}] (23) at (4, 5 * \triht + \upmid) {};
		\node [style=small black dot] (24) at (3.5, 5 * \triht + \downmid) {};
		\node [style=small black dot] (25) at (-4.5, \upmid) {};
		\node [style=small black dot] (26) at (-4.5, -1 * \triht + \downmid) {};
		\node [style=small black dot] (27) at (2, -5 * \triht + \upmid) {};
		\node [style=small black dot, label={right:$w'$}] (28) at (4, -4 * \triht + \downmid) {};
		\node [style=none] (29) at (4.5, 5 * \triht) {};
		\node [style=none] (30) at (4, 6 * \triht) {};
		\node [style=none] (31) at (3, 6 * \triht) {};
		\node [style=none] (32) at (2.5, 5 * \triht) {};
		\node [style=none] (34) at (4, 4 * \triht) {};
		\node [style=none] (33) at (3, 4 * \triht) {};
		\node [style=none] (39) at (-4.5, 1 * \triht) {};
		\node [style=none] (40) at (-5, 0) {};
		\node [style=none] (41) at (-4.5, -1 * \triht) {};
		\node [style=none] (42) at (-3.5, -1 * \triht) {};
		\node [style=none] (43) at (-3, 0) {};
		\node [style=none] (44) at (-3.5, 1 * \triht) {};
		\node [style=none] (49) at (1.5, -5 * \triht) {};
		\node [style=none] (50) at (2.5, -5 * \triht) {};
		\node [style=none] (51) at (3, -4 * \triht) {};
		\node [style=none] (52) at (2, -2 * \triht) {};
		\node [style=none] (53) at (1, -4 * \triht) {};
		\node [style=none] (54) at (4, -4 * \triht) {};
		\node [style=none] (55) at (4.5, -3 * \triht) {};
		\node [style=none] (56) at (4, -2 * \triht) {};
		\node [style=none] (57) at (2, -2 * \triht) {};
		\node [style=none] (58) at (3, -4 * \triht) {};
		\node [style=none] (59) at (3.5, 5 * \triht) {};
		\node [style=none] (60) at (3.5, -3 * \triht) {};
		\node [style=none] (61) at (2, -4 * \triht) {};
		\node [style=none] (62) at (-4, 0) {};
		\node [style=none] (63) at (3, -3 * \triht + \upmid) {};
		\node [style=none] (64) at (2, -4 * \triht + \downmid) {};
		\node [style=orange dot, label={above:$p$}] (65) at (2, 4 * \triht + \downmid) {};
		\node [style=orange dot, label={above:$\tilde p$}] (66) at (-3, 1 * \triht + \upmid) {};
		\node [style=none] (70) at (2.5, 3 * \triht) {};
		\node [style=none] (71) at (-2, 0) {};
	\end{pgfonlayer}
	\begin{pgfonlayer}{edgelayer}
		\draw [style=red line] (2.center) to (3.center);
		\draw [style=red line] (5.center) to (6.center);
		\draw [style=green line] (7.center) to (8.center);
		\draw [style=blue line] (12.center) to (11.center);
		\draw [style=filled grey] (34.center)
			 to (29.center)
			 to (30.center)
			 to (31.center)
			 to (32.center)
			 to (33.center)
			 to cycle;
		\draw [style=filled grey] (44.center)
			 to (39.center)
			 to (40.center)
			 to (41.center)
			 to (42.center)
			 to (43.center)
			 to cycle;
		\draw [style=filled grey] (52.center)
			 to (53.center)
			 to (49.center)
			 to (50.center)
			 to (51.center)
			 to cycle;
		\draw [style=filled grey] (58.center)
			 to (54.center)
			 to (55.center)
			 to (56.center)
			 to (57.center)
			 to cycle;
		\draw [style=dashed line] (13) to (14);
		\draw [style=dashed line] (19) to (13);
		\draw [style=dashed line] (20) to (19);
		\draw [style=dashed line] (21) to (20);
		\draw [style=dashed line] (14) to (22);
		\draw [style=dashed line] (22) to (21);
		\draw (23) to (24);
		\draw (24) to (25);
		\draw (25) to (26);
		\draw (26) to (27);
		\draw (27) to (28);
		\draw (28) to (23);
		\draw [style=dotted line] (59.center) to (62.center);
		\draw [style=dotted line] (62.center) to (61.center);
		\draw [style=dotted line] (61.center) to (60.center);
		\draw [style=dotted line] (60.center) to (59.center);
		\draw [style=dotted line] (33.center) to (43.center);
		\draw [style=dotted line] (43.center) to (64.center);
		\draw [style=dotted line] (64.center) to (63.center);
		\draw [style=dotted line] (63.center) to (33.center);
		\draw [style=dotted line] (13.center) to (70.center);
		\draw [style=dotted line] (14.center) to (70.center);
		\draw [style=dotted line] (21.center) to (71.center);
		\draw [style=dotted line] (22.center) to (71.center);
	\end{pgfonlayer}
	

    \clip       (-\cols, -\height) rectangle (\cols, \height);
    \draw[gray] (-\cols, -\height) rectangle (\cols, \height);

    \pgfmathsetmacro{\from}{-2 *\cols}
    \pgfmathsetmacro{\to}{2 * \cols}
    \foreach\i in {\from, ..., \to} {
        \draw[xslant=\slant]  (\i, -\height) -- (\i, \height);
        \draw[xslant=-\slant] (\i, -\height) -- (\i, \height);
    }

    \foreach\j in {-\rows, ..., \rows} {
        \pgfmathsetmacro{\y}{0.5 * \j * tan(60)}
        \draw (-\cols, \y) -- (\cols, \y);
    }

\end{tikzpicture}

\caption{3-shell of $\ch_{w}$: $w$ of Type 1, adjacent to chamber wall}  \label{fig:edge3shell}
\end{figure}

(ii) If $H$ does not lie entirely within the chamber of $w$, then $w$ must be in a root strip adjacent to one of the fundamental root strips (i.e., a base case). If $w$ were to lie in an even chamber, then $w$ would be twisted spiral (Base Case 1), contrary to assumption. So $w$ is in an odd chamber (and in Base Case 3). Three of the alcoves in $H$ lie in the chamber, and the other three lie in the adjacent fundamental root strip. Reflecting $H$ in the \emph{other} wall of this chamber produces a new hexagon $H'$ lying entirely within its chamber. If $w'\in H'$ corresponds to $w$ under that reflection, then $H'$ consists of the alcoves for the six elements $w'r$ for $r \in R(w)$. Letting $z':=w'tuts=w'utus = t(-\ga') w'$ (where $t(\ga')$ is translation into the chamber of $w'$), one checks that $\ch_{z'}$ is the 3-shell of $\ch_{w}$ (provided the latter is non-empty); see Figure \ref{fig:edge3shell}. Now the result of the theorem for $x$ on or inside the 3-shell of $\ch_{w}$ follows as before.

\begin{figure}[htbp!]
%

\pgfmathsetmacro{\cols}{7}
\pgfmathsetmacro{\rows}{8}
\pgfmathsetmacro{\slant}{cot(60)}
\pgfmathsetmacro{\height}{0.5 * \rows * tan(60)}
\pgfmathsetmacro{\triht}{sin(60)}
\pgfmathsetmacro{\upmid}{0.25 * sec(30)}
\pgfmathsetmacro{\downmid}{\triht - \upmid}

\begin{tikzpicture}
    

\input{sample.tikzstyles}

	\begin{pgfonlayer}{nodelayer}
		\node [style=none] (2) at (-7, 1 * \triht) {};
		\node [style=none] (3) at (7, 1 * \triht) {};
		\node [style=red dot] (4) at (0, 0) {};
		\node [style=none] (7) at (-3, -8 * \triht) {};
		\node [style=none] (8) at (5, 8 * \triht) {};
		\node [style=none] (9) at (5.5, 4 * \triht + \upmid) {$X$};
		\node [style=none] (11) at (5, -8 * \triht) {};
		\node [style=none] (12) at (-3, 8 * \triht) {};
		\node [style=small black dot, label={right:$z_2$}] (13) at (4.5, 4 * \triht + \upmid) {};
		\node [style=small black dot] (14) at (2.5, 5 * \triht + \downmid) {};
		\node [style=small black dot] (19) at (4.5, -3 * \triht + \downmid) {};
		\node [style=small black dot] (20) at (1, -5 * \triht + \upmid) {};
		\node [style=small black dot] (21) at (-4, -2 * \triht + \downmid) {};
		\node [style=small black dot] (22) at (-4, 1 * \triht + \upmid) {};
		\node [style=small black dot, label={right:$w$}] (23) at (6, 4 * \triht + \downmid) {};
		\node [style=small black dot] (24) at (2, 7 * \triht + \upmid) {};
		\node [style=small black dot] (25) at (-5, 2 * \triht + \downmid) {};
		\node [style=small black dot] (26) at (-5, -3 * \triht + \upmid) {};
		\node [style=small black dot] (27) at (0.5, -7 * \triht + \downmid) {};
		\node [style=small black dot] (28) at (6, -3 * \triht + \upmid) {};
		\node [style=none] (29) at (5.5, 4 * \triht + \upmid) {};
		\node [style=none] (30) at (2, 6 * \triht + \downmid) {};
		\node [style=none] (31) at (-4.5, 2 * \triht + \upmid) {};
		\node [style=none] (32) at (-4.5, -3 * \triht + \downmid) {};
		\node [style=none] (33) at (0.5, -6 * \triht + \upmid) {};
		\node [style=none] (99) at (5.5, -3 * \triht + \downmid) {};
		\node [style=none] (65) at (5.5, 3 * \triht) {};
		\node [style=none] (66) at (6.5, 5 * \triht) {};
		\node [style=none] (67) at (4.5, 5 * \triht) {};
		\node [style=none] (68) at (4, 4 * \triht) {};
		\node [style=none] (69) at (4.5, 3 * \triht) {};
		\node [style=none] (34) at (3, 6 * \triht) {};
		\node [style=none] (35) at (2, 8 * \triht) {};
		\node [style=none] (36) at (1, 6 * \triht) {};
		\node [style=none] (37) at (1.5, 5 * \triht) {};
		\node [style=none] (38) at (2.5, 5 * \triht) {};
		\node [style=none] (39) at (-3.5, 3 * \triht) {};
		\node [style=none] (40) at (-5.5, 3 * \triht ) {};
		\node [style=none] (41) at (-4.5, 1 * \triht) {};
		\node [style=none] (42) at (-3.5, 1 * \triht) {};
		\node [style=none] (43) at (-3, 2 * \triht) {};
		\node [style=none] (44) at (-4.5, -1 * \triht) {};
		\node [style=none] (45) at (-5.5, -3 * \triht) {};
		\node [style=none] (46) at (-3.5, -3 * \triht) {};
		\node [style=none] (47) at (-3, -2 * \triht) {};
		\node [style=none] (48) at (-3.5, -1 * \triht) {};
		\node [style=none] (49) at (-0.5, -5 * \triht) {};
		\node [style=none] (50) at (0.5, -7 * \triht) {};
		\node [style=none] (51) at (1.5, -5 * \triht) {};
		\node [style=none] (52) at (1, -4 * \triht) {};
		\node [style=none] (53) at (0, -4 * \triht) {};
		\node [style=none] (54) at (4.5, -3 * \triht) {};
		\node [style=none] (55) at (6.5, -3 * \triht) {};
		\node [style=none] (56) at (5.5, -1 * \triht) {};
		\node [style=none] (57) at (4.5, -1 * \triht) {};
		\node [style=none] (58) at (4, -2 * \triht) {};
		\node [style=small black dot, label={right:$z_1$}] (59) at (5, 3 * \triht + \upmid) {};
		\node [style=small black dot] (60) at (5, -2 * \triht + \downmid) {};
		\node [style=small black dot] (61) at (0, -5 * \triht + \upmid) {};
		\node [style=small black dot] (62) at (-3.5, -3 * \triht + \downmid) {};
		\node [style=small black dot] (63) at (-3.5, 2 * \triht + \upmid) {};
		\node [style=small black dot] (64) at (1.5, 5 * \triht + \downmid) {};
		\node [style=none] (70) at (5, 4 * \triht) {};
		\node [style=none] (71) at (2, 6 * \triht) {};
		\node [style=none] (72) at (-4, 2 * \triht) {};
		\node [style=none] (73) at (-4, -2 * \triht) {};
		\node [style=none] (74) at (0.5, -5 * \triht) {};
		\node [style=none] (75) at (5, -2 * \triht) {};
	\end{pgfonlayer}
	\begin{pgfonlayer}{edgelayer}
		\draw [style=red line] (2.center) to (3.center);
		\draw [style=green line] (7.center) to (8.center);
		\draw [style=blue line] (12.center) to (11.center);
		\draw [style=filled grey] (69.center)
			 to (65.center)
			 to (66.center)
			 to (67.center)
			 to (68.center)
			 to cycle;
		\draw [style=filled grey] (36.center)
			 to (37.center)
			 to (38.center)
			 to (34.center)
			 to (35.center)
			 to cycle;
		\draw [style=filled grey] (43.center)
			 to (39.center)
			 to (40.center)
			 to (41.center)
			 to (42.center)
			 to cycle;
		\draw [style=filled grey] (47.center)
			 to (48.center)
			 to (44.center)
			 to (45.center)
			 to (46.center)
			 to cycle;
		\draw [style=filled grey] (52.center)
			 to (53.center)
			 to (49.center)
			 to (50.center)
			 to (51.center)
			 to cycle;
		\draw [style=filled grey] (58.center)
			 to (54.center)
			 to (55.center)
			 to (56.center)
			 to (57.center)
			 to cycle;
		\draw [style=dashed line] (13) to (14);
		\draw [style=dashed line] (19) to (13);
		\draw [style=dashed line] (20) to (19);
		\draw [style=dashed line] (21) to (20);
		\draw [style=dashed line] (14) to (22);
		\draw [style=dashed line] (22) to (21);
		\draw (23) to (24);
		\draw (24) to (25);
		\draw (25) to (26);
		\draw (26) to (27);
		\draw (27) to (28);
		\draw (28) to (23);
		\draw [style=dashed line] (59) to (60);
		\draw [style=dashed line] (60) to (61);
		\draw [style=dashed line] (61) to (62);
		\draw [style=dashed line] (62) to (63);
		\draw [style=dashed line] (63) to (64);
		\draw [style=dashed line] (64) to (59);
		\draw [style=dotted line] (29.center) to (30.center);
		\draw [style=dotted line] (30.center) to (31.center);
		\draw [style=dotted line] (31.center) to (32.center);
		\draw [style=dotted line] (32.center) to (33.center);
		\draw [style=dotted line] (33.center) to (99.center);
		\draw [style=dotted line] (99.center) to (29.center);
		\draw [style=dotted line] (59.center) to (70.center);
		\draw [style=dotted line] (13.center) to (70.center);
		\draw [style=dotted line] (14.center) to (71.center);
		\draw [style=dotted line] (64.center) to (71.center);
		\draw [style=dotted line] (63.center) to (72.center);
		\draw [style=dotted line] (22.center) to (72.center);
		\draw [style=dotted line] (21.center) to (73.center);
		\draw [style=dotted line] (62.center) to (73.center);
		\draw [style=dotted line] (61.center) to (74.center);
		\draw [style=dotted line] (20.center) to (74.center);
		\draw [style=dotted line] (60.center) to (75.center);
		\draw [style=dotted line] (19.center) to (75.center);
	\end{pgfonlayer}


    \clip       (-\cols, -\height) rectangle (\cols, \height);
    \draw[gray] (-\cols, -\height) rectangle (\cols, \height);

    \pgfmathsetmacro{\from}{-2 *\cols}
    \pgfmathsetmacro{\to}{2 * \cols}
    \foreach\i in {\from, ..., \to} {
        \draw[xslant=\slant]  (\i, -\height) -- (\i, \height);
        \draw[xslant=-\slant] (\i, -\height) -- (\i, \height);
    }

    \foreach\j in {-\rows, ..., \rows} {
        \pgfmathsetmacro{\y}{0.5 * \j * tan(60)}
        \draw (-\cols, \y) -- (\cols, \y);
    }

\end{tikzpicture}
\caption{2-shell of $\ch_{w}$: $w$ of Type 2}  \label{fig:2shell}
\end{figure}

\noindent
2. Next, consider $w$ of Type $\tau=2$. We need to analyze the elements on or inside the $4-\tau = 2$-shell of $\ch_{w}$. Label the simple reflections $s, t, u$ so that $w>ws$, $w<wt$, and $w<wu$. Let $X$ be the convex pentagonal heptiamond consisting of the alcoves $w$ and  $wsr$ for $r\in\langle t, u \rangle=R(ws)$. 

(i) If $X$ lies completely within the chamber of $w$, then we can reflect it as in \eqref{e:hexagonVertices} to form copies of $X$ at each of the vertices of $\ch_{w}$; see Figure \ref{fig:2shell}. Let $z_{1}:= wstu$ and $z_{2}:=wsut$, and consider $\ch_{z_{i}}$ for $i=1, 2$. These hexagons lie on or inside the 2-shell of $\ch_{w}$, with the 2-shell alcoves consisting of every second edge of $\ch_{z_{1}}$ going one way, and every second edge of $\ch_{z_{2}}$ going the other. Hence if $x$ lies on or inside the 2-shell of $\ch_{w}$, then $x\le z_{1}$ or $x\le z_{2}$, and we complete the proof of the theorem as before.

(ii) Now suppose $X$ does not lie completely within the chamber of $w$, so $w$ lies in a root strip adjacent to one of the walls of its chamber (i.e., is a base case). Three of the alcoves of the hexagonal part of $X$ lie in the chamber of $w$; relabel if necessary so that $z_{1}$ is one of these. The other three hexagon alcoves lie inside an adjacent fundamental root strip; one of these is $z_{2}$. 

\begin{figure}[htbp!]
%

\pgfmathsetmacro{\cols}{6}
\pgfmathsetmacro{\rows}{7}
\pgfmathsetmacro{\slant}{cot(60)}
\pgfmathsetmacro{\height}{0.5 * \rows * tan(60)}
\pgfmathsetmacro{\triht}{sin(60)}
\pgfmathsetmacro{\upmid}{0.25 * sec(30)}
\pgfmathsetmacro{\downmid}{\triht - \upmid}

\begin{tikzpicture}
    

\input{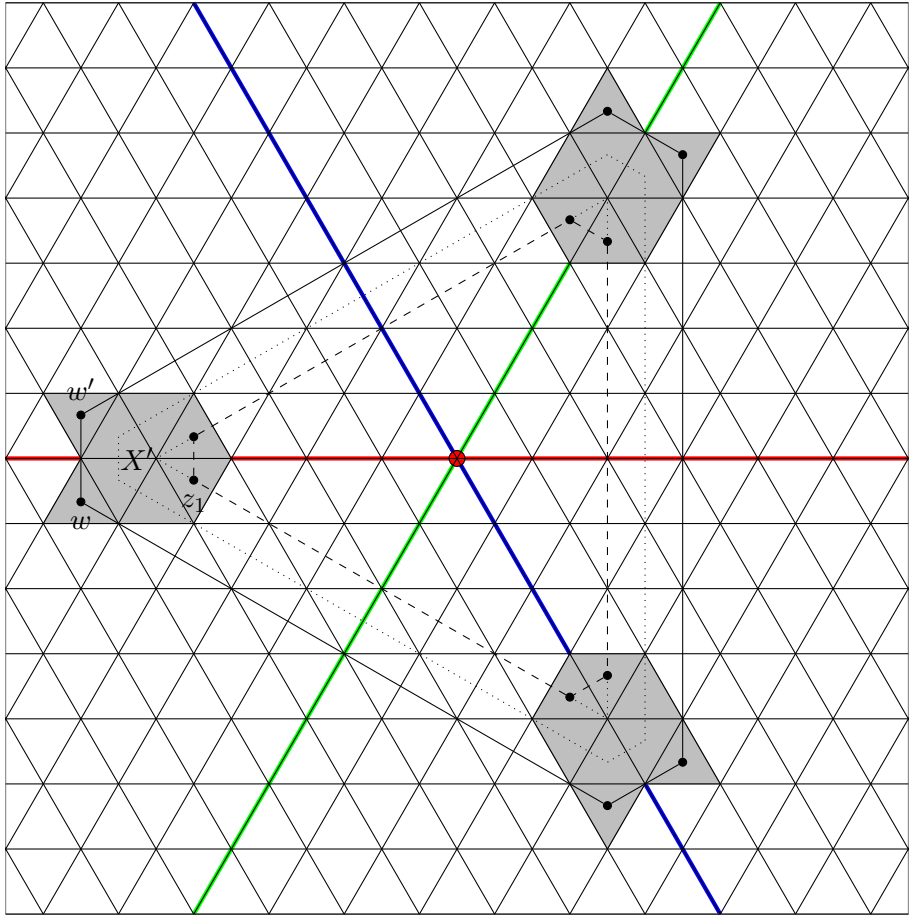}

	\begin{pgfonlayer}{nodelayer}
		\node [style=none] (2) at (-6, 0) {};
		\node [style=none] (3) at (6, 0) {};
		\node [style=red dot] (4) at (0, 0) {};
		\node [style=none] (7) at (-3.5, -7 * \triht) {};
		\node [style=none] (8) at (3.5, 7 * \triht) {};
		\node [style=none] (11) at (3.5, -7 * \triht) {};
		\node [style=none] (12) at (-3.5, 7 * \triht) {};
		\node [style=small black dot] (23) at (3, 4 * \triht + \downmid) {};
		\node [style=small black dot] (24) at (2, 5 * \triht + \upmid) {};
		\node [style=small black dot, label={above:$w'$}] (25) at (-5,  \downmid) {};
		\node [style=small black dot, label={below:$w$}] (26) at (-5, -1 * \triht + \upmid) {};
		\node [style=small black dot] (27) at (2, -6 * \triht + \downmid) {};
		\node [style=small black dot] (28) at (3, -5 * \triht + \upmid) {};
		\node [style=none] (29) at (2.5, 4 * \triht + \upmid) {};
		\node [style=none] (30) at (2, 4 * \triht + \downmid) {};
		\node [style=none] (31) at (-4.5,  \upmid) {};
		\node [style=none] (32) at (-4.5, -1 * \triht + \downmid) {};
		\node [style=none] (33) at (2, -5 * \triht + \upmid) {};
		\node [style=none] (99) at (2.5, -5 * \triht + \downmid) {};
		\node [style=none] (65) at (2.5, 3 * \triht) {};
		\node [style=none] (66) at (3.5, 5 * \triht) {};
		\node [style=none] (67) at (2, 6 * \triht) {};
		\node [style=none] (69) at (1.5, 3 * \triht) {};
		\node [style=none] (38) at (1, 4 * \triht) {};
		\node [style=none] (39) at (-3.5, 1 * \triht) {};
		\node [style=none] (40) at (-5.5, 1 * \triht) {};
		\node [style=none] (41) at (-5.5, -1 * \triht) {};
		\node [style=none] (42) at (-3.5, -1 * \triht) {};
		\node [style=none] (43) at (-3, 0) {};
		\node [style=none] (48) at (-3.5, -1 * \triht) {};
		\node [style=none] (52) at (2.5, -3 * \triht) {};
		\node [style=none] (55) at (3.5, -5 * \triht) {};
		\node [style=none] (56) at (2.5, -3 * \triht) {};
		\node [style=none] (57) at (1.5, -3 * \triht) {};
		\node [style=none] (58) at (1, -4 * \triht) {};
		\node [style=small black dot] (59) at (2, 3 * \triht + \upmid) {};
		\node [style=small black dot] (60) at (2, -4 * \triht + \downmid) {};
		\node [style=small black dot] (61) at (1.5, -4 * \triht + \upmid) {};
		\node [style=small black dot, label={below:$z_1$}] (62) at (-3.5, -1 * \triht + \downmid) {};
		\node [style=small black dot] (63) at (-3.5,  \upmid) {};
		\node [style=small black dot] (64) at (1.5, 3 * \triht + \downmid) {};
		\node [style=none] (100) at (2.5, 5 * \triht) {};
		\node [style=none] (101) at (-5, 0) {};
		\node [style=none] (102) at (-4.25, 0) {$X'$};
		\node [style=none] (104) at (2, -6 * \triht) {};
		\node [style=none] (105) at (2.5, -5 * \triht) {};
		\node [style=none] (110) at (2, 4 * \triht) {};
		\node [style=none] (111) at (-4, -0) {};
		\node [style=none] (112) at (2, -4 * \triht) {};
	\end{pgfonlayer}
	\begin{pgfonlayer}{edgelayer}
		\draw [style=red line] (2.center) to (3.center);
		\draw [style=green line] (7.center) to (8.center);
		\draw [style=blue line] (12.center) to (11.center);
		\draw [style=filled grey] (69.center)
			 to (65.center)
			 to (66.center)
			 to (100.center)
			 to (67.center)
			 to (38.center)
			 to cycle;
		\draw [style=filled grey] (43.center)
			 to (39.center)
			 to (40.center)
			 to (101.center)
			 to (41.center)
			 to (42.center)
			 to cycle;
		\draw [style=filled grey] (104.center)
			 to (105.center)
			 to (55.center)
			 to (56.center)
			 to (57.center)
			 to (58.center)
			 to cycle;
		\draw (23) to (24);
		\draw (24) to (25);
		\draw (25) to (26);
		\draw (26) to (27);
		\draw (27) to (28);
		\draw (28) to (23);
		\draw [style=dashed line] (59) to (60);
		\draw [style=dashed line] (60) to (61);
		\draw [style=dashed line] (61) to (62);
		\draw [style=dashed line] (62) to (63);
		\draw [style=dashed line] (63) to (64);
		\draw [style=dashed line] (64) to (59);
		\draw [style=dotted line] (29.center) to (30.center);
		\draw [style=dotted line] (30.center) to (31.center);
		\draw [style=dotted line] (31.center) to (32.center);
		\draw [style=dotted line] (32.center) to (33.center);
		\draw [style=dotted line] (33.center) to (99.center);
		\draw [style=dotted line] (99.center) to (29.center);
		\draw [style=dotted line] (59.center) to (110.center);
		\draw [style=dotted line] (64.center) to (110.center);
		\draw [style=dotted line] (63.center) to (111.center);
		\draw [style=dotted line] (62.center) to (111.center);
		\draw [style=dotted line] (61.center) to (112.center);
		\draw [style=dotted line] (60.center) to (112.center);
	\end{pgfonlayer}


    \clip       (-\cols, -\height) rectangle (\cols, \height);
    \draw[gray] (-\cols, -\height) rectangle (\cols, \height);

    \pgfmathsetmacro{\from}{-2 *\cols}
    \pgfmathsetmacro{\to}{2 * \cols}
    \foreach\i in {\from, ..., \to} {
        \draw[xslant=\slant]  (\i, -\height) -- (\i, \height);
        \draw[xslant=-\slant] (\i, -\height) -- (\i, \height);
    }

    \foreach\j in {-\rows, ..., \rows} {
        \pgfmathsetmacro{\y}{0.5 * \j * tan(60)}
        \draw (-\cols, \y) -- (\cols, \y);
    }

\end{tikzpicture}
\caption{2-shell of $\ch_{w}$: $w$ of Type 2, adjacent to even chamber wall}  \label{fig:fish}
\end{figure}

(a) If the chamber of $w$ is even (Base Case 2), it is convenient to adjoin one additional spiral alcove to $X$, corresponding to the reflection $w'$ of $w$ across the adjacent chamber wall. Call the alcoves $w$ and $w'$, lying outside the hexagonal part, the ``tail'' of the resulting shape $X'$ (it looks somewhat like a fish; see Figure \ref{fig:fish}). Reflect $X'$ to form two additional ``fish,'' so that the vertices of $\ch_{w}$ lie at the centers of the six tail alcoves of the three fish. Figure \ref{fig:fish} depicts the situation when $w$ lies in Chamber IV. The six ``fish head'' alcove centers are the vertices of $\ch_{z_{1}}$, and this is also the 2-shell of $\ch_{w}$. The theorem follows (for $x$ in the 2-shell) as before. 

\begin{figure}[htbp!]
%

\pgfmathsetmacro{\cols}{7}
\pgfmathsetmacro{\rows}{8}
\pgfmathsetmacro{\slant}{cot(60)}
\pgfmathsetmacro{\height}{0.5 * \rows * tan(60)}
\pgfmathsetmacro{\triht}{sin(60)}
\pgfmathsetmacro{\upmid}{0.25 * sec(30)}
\pgfmathsetmacro{\downmid}{\triht - \upmid}

\begin{tikzpicture}
    

\input{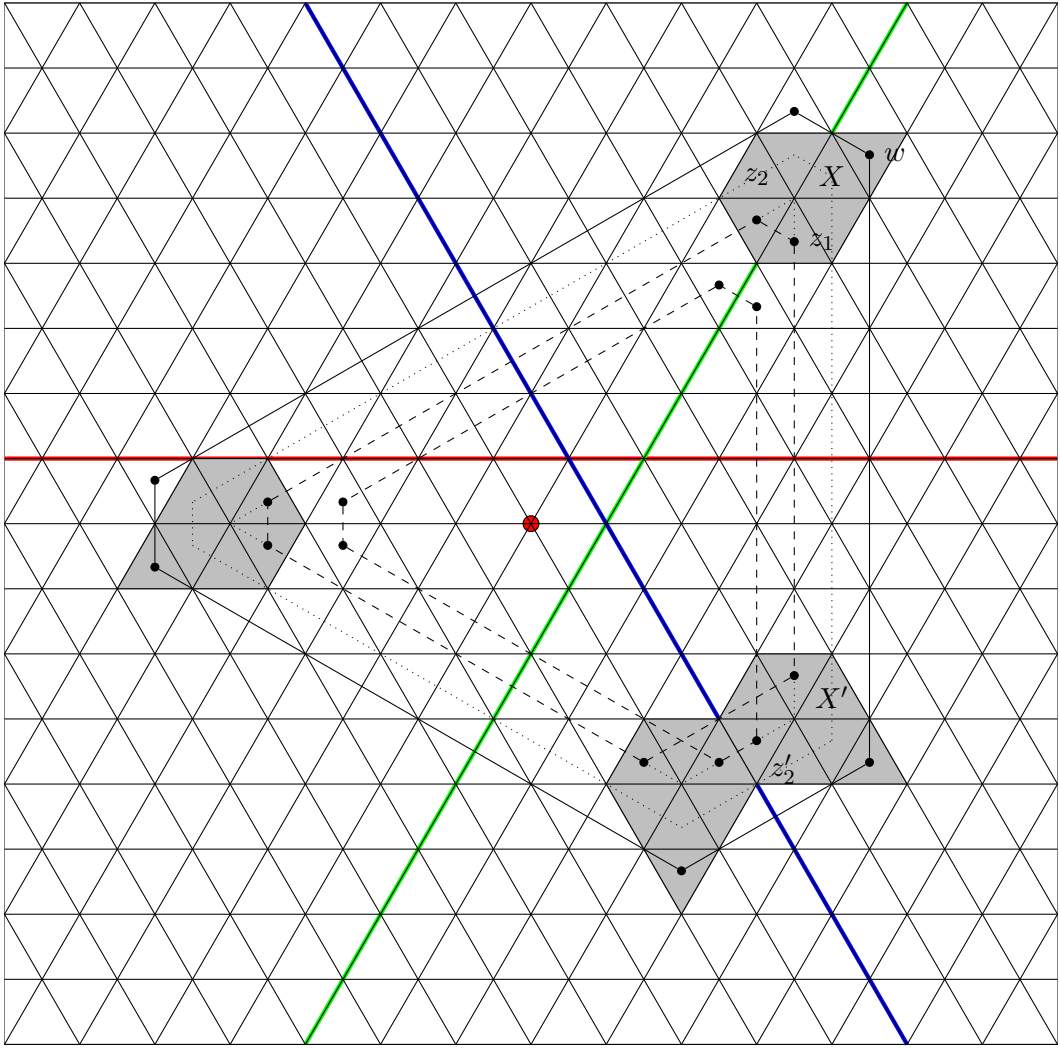}

	\begin{pgfonlayer}{nodelayer}
		\node [style=none] (2) at (-7, 1 * \triht) {};
		\node [style=none] (3) at (7, 1 * \triht) {};
		\node [style=red dot] (4) at (0, 0) {};
		\node [style=none] (7) at (-3, -8 * \triht) {};
		\node [style=none] (8) at (5, 8 * \triht) {};
		\node [style=none] (11) at (5, -8 * \triht) {};
		\node [style=none] (12) at (-3, 8 * \triht) {};
		\node [style=small black dot] (13) at (3, 3 * \triht + \upmid) {};
		\node [style=small black dot] (14) at (2.5, 3 * \triht + \downmid) {};
		\node [style=small black dot, label={below right:$z'_2$}] (19) at (3, -4 * \triht + \downmid) {};
		\node [style=small black dot] (20) at (2.5, -4 * \triht + \upmid) {};
		\node [style=small black dot] (21) at (-2.5, -1 * \triht + \downmid) {};
		\node [style=small black dot] (22) at (-2.5,  \upmid) {};
		\node [style=small black dot, label={right:$w$}] (23) at (4.5, 5 * \triht + \downmid) {};
		\node [style=small black dot] (24) at (3.5, 6 * \triht + \upmid) {};
		\node [style=small black dot] (25) at (-5,  \downmid) {};
		\node [style=small black dot] (26) at (-5, -1 * \triht + \upmid) {};
		\node [style=small black dot] (27) at (2, -6 * \triht + \downmid) {};
		\node [style=small black dot] (28) at (4.5, -4 * \triht + \upmid) {};
		\node [style=none] (29) at (4, 5 * \triht + \upmid) {};
		\node [style=none] (30) at (3.5, 5 * \triht + \downmid) {};
		\node [style=none] (31) at (-4.5,  \upmid) {};
		\node [style=none] (32) at (-4.5, -1 * \triht + \downmid) {};
		\node [style=none] (33) at (2, -5 * \triht + \upmid) {};
		\node [style=none] (99) at (4, -4 * \triht + \downmid) {};
		\node [style=none] (65) at (4, 4 * \triht) {};
		\node [style=none] (66) at (5, 6 * \triht) {};
		\node [style=none] (68) at (2.5, 5 * \triht) {};
		\node [style=none] (69) at (3, 4 * \triht) {};
		\node [style=none] (38) at (2.5, 5 * \triht) {};
		\node [style=none] (39) at (-3.5, 1 * \triht) {};
		\node [style=none] (41) at (-5.5, -1 * \triht) {};
		\node [style=none] (42) at (-3.5, -1 * \triht) {};
		\node [style=none] (43) at (-3, 0) {};
		\node [style=none] (48) at (-3.5, -1 * \triht) {};
		\node [style=none] (49) at (1, -4 * \triht) {};
		\node [style=none] (50) at (2, -6 * \triht) {};
		\node [style=none] (51) at (3, -4 * \triht) {};
		\node [style=none] (52) at (2.5, -3 * \triht) {};
		\node [style=none] (53) at (1.5, -3 * \triht) {};
		\node [style=none] (54) at (3, -4 * \triht) {};
		\node [style=none] (55) at (5, -4 * \triht) {};
		\node [style=none] (56) at (4, -2 * \triht) {};
		\node [style=none] (57) at (3, -2 * \triht) {};
		\node [style=none] (58) at (2.5, -3 * \triht) {};
		\node [style=small black dot, label={right:$z_1$}] (59) at (3.5, 4 * \triht + \upmid) {};
		\node [style=small black dot] (60) at (3.5, -3 * \triht + \downmid) {};
		\node [style=small black dot] (61) at (1.5, -4 * \triht + \upmid) {};
		\node [style=small black dot] (62) at (-3.5, -1 * \triht + \downmid) {};
		\node [style=small black dot] (63) at (-3.5,  \upmid) {};
		\node [style=small black dot] (64) at (3, 4 * \triht + \downmid) {};
		\node [style=none] (70) at (3, 5 * \triht + \upmid) {$z_2$};
		\node [style=none] (100) at (3, 6 * \triht) {};
		\node [style=none] (101) at (-4.5, 1 * \triht) {};
		\node [style=none] (102) at (4, 5 * \triht + \upmid) {$X$};
		\node [style=none] (103) at (4, -3 * \triht + \upmid) {$X'$};
		\node [style=none] (121) at (3.5, 5 * \triht) {};
		\node [style=none] (122) at (-4, 0) {};
		\node [style=none] (123) at (3.5, -3 * \triht) {};
		\node [style=none] (124) at (2, -4 * \triht) {};
	\end{pgfonlayer}
	\begin{pgfonlayer}{edgelayer}
		\draw [style=red line] (2.center) to (3.center);
		\draw [style=green line] (7.center) to (8.center);
		\draw [style=blue line] (12.center) to (11.center);
		\draw [style=filled grey] (68.center)
			 to (69.center)
			 to (65.center)
			 to (66.center)
			 to (100.center)
			 to (38.center);
		\draw [style=filled grey] (43.center)
			 to (39.center)
			 to (101.center)
			 to (41.center)
			 to (42.center)
			 to cycle;
		\draw [style=filled grey] (52.center)
			 to (53.center)
			 to (49.center)
			 to (50.center)
			 to (51.center)
			 to cycle;
		\draw [style=filled grey] (58.center)
			 to (54.center)
			 to (55.center)
			 to (56.center)
			 to (57.center)
			 to cycle;
		\draw [style=dashed line] (13) to (14);
		\draw [style=dashed line] (19) to (13);
		\draw [style=dashed line] (20) to (19);
		\draw [style=dashed line] (21) to (20);
		\draw [style=dashed line] (14) to (22);
		\draw [style=dashed line] (22) to (21);
		\draw (23) to (24);
		\draw (24) to (25);
		\draw (25) to (26);
		\draw (26) to (27);
		\draw (27) to (28);
		\draw (28) to (23);
		\draw [style=dashed line] (59) to (60);
		\draw [style=dashed line] (60) to (61);
		\draw [style=dashed line] (61) to (62);
		\draw [style=dashed line] (62) to (63);
		\draw [style=dashed line] (63) to (64);
		\draw [style=dashed line] (64) to (59);
		\draw [style=dotted line] (29.center) to (30.center);
		\draw [style=dotted line] (30.center) to (31.center);
		\draw [style=dotted line] (31.center) to (32.center);
		\draw [style=dotted line] (32.center) to (33.center);
		\draw [style=dotted line] (33.center) to (99.center);
		\draw [style=dotted line] (99.center) to (29.center);
		\draw [style=dotted line] (59.center) to (121.center);
		\draw [style=dotted line] (64.center) to (121.center);
		\draw [style=dotted line] (62.center) to (122.center);
		\draw [style=dotted line] (63.center) to (122.center);
		\draw [style=dotted line] (60.center) to (123.center);
		\draw [style=dotted line] (19.center) to (123.center);
		\draw [style=dotted line] (20.center) to (124.center);
		\draw [style=dotted line] (61.center) to (124.center);
	\end{pgfonlayer}


    \clip       (-\cols, -\height) rectangle (\cols, \height);
    \draw[gray] (-\cols, -\height) rectangle (\cols, \height);

    \pgfmathsetmacro{\from}{-2 *\cols}
    \pgfmathsetmacro{\to}{2 * \cols}
    \foreach\i in {\from, ..., \to} {
        \draw[xslant=\slant]  (\i, -\height) -- (\i, \height);
        \draw[xslant=-\slant] (\i, -\height) -- (\i, \height);
    }

    \foreach\j in {-\rows, ..., \rows} {
        \pgfmathsetmacro{\y}{0.5 * \j * tan(60)}
        \draw (-\cols, \y) -- (\cols, \y);
    }

\end{tikzpicture}
\caption{2-shell of $\ch_{w}$: $w$ of Type 2, adjacent to odd chamber wall}  \label{fig:edge2shell}
\end{figure}

(b) Now assume the chamber of $w$ is odd (Base Case 4). Reflect $X$ across the other wall of the chamber to obtain $X'$ (see Figure \ref{fig:edge2shell}). Let $z'_{2}$ be the reflection of $z_{2}$. Observe that five of the six edges of $\ch_{z_{1}}$ coincide with five of the six edges of the 2-shell of $\ch_{w}$, while the sixth edge of $\ch_{z_{1}}$ is inside the 2-shell. The remaining edge of the 2-shell of $\ch_{w}$ coincides with one of the edges of $\ch_{z'_{2}}$, whose remaining five edges lie inside the 2-shell. These observations imply that the set of alcoves on or inside the 2-shell of $\ch_{w}$ equals the union of the alcoves on or inside $\ch_{z_{1}}$ and $\ch_{z'_{2}}$. Now the claim of the theorem follows (for $x$ in the 2-shell) as before. 

\medskip\noindent
II. Now consider the situation $y\le x\le w$ with $q_{x}^{w}>0$, $w$ of Type $\tau=1$ or $2$, and $x$ is on the $k$-shell of $\ch_{w}$ for $0\le k < 4-\tau$. Then $x$ must be of the form $x=x_{0}v$ for $v\in R(w)$ and $x_{0}$ is on a special segment $E$ of $\ch_{w}$, as in Theorem \ref{t:outershells}.

\begin{figure}[htbp!]
%

\pgfmathsetmacro{\cols}{7}
\pgfmathsetmacro{\rows}{8}
\pgfmathsetmacro{\slant}{cot(60)}
\pgfmathsetmacro{\height}{0.5 * \rows * tan(60)}
\pgfmathsetmacro{\triht}{sin(60)}
\pgfmathsetmacro{\upmid}{0.25 * sec(30)}
\pgfmathsetmacro{\downmid}{\triht - \upmid}

\begin{tikzpicture}
    

\input{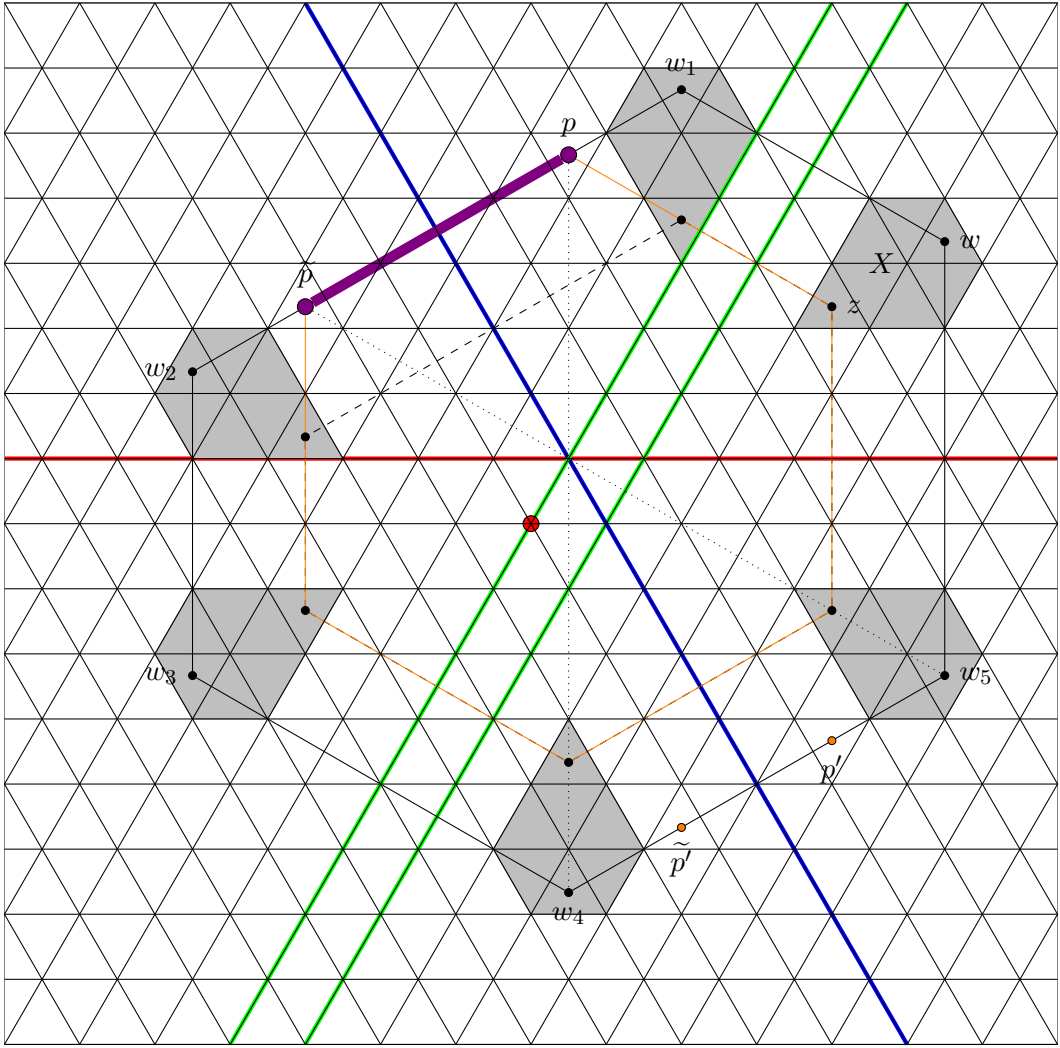}

	\begin{pgfonlayer}{nodelayer}
		\node [style=none] (2) at (-7, 1 * \triht) {};
		\node [style=none] (3) at (7, 1 * \triht) {};
		\node [style=red dot] (4) at (0, 0) {};
		\node [style=none, label={below left:$X$}] (5) at (5, 4 * \triht + \upmid) {};
		\node [style=none] (7) at (-3, -8 * \triht) {};
		\node [style=none] (8) at (5, 8 * \triht) {};
		\node [style=none] (11) at (5, -8 * \triht) {};
		\node [style=none] (12) at (-3, 8 * \triht) {};
		\node [style=small black dot, label={right:$z$}] (13) at (4, 3 * \triht + \upmid) {};
		\node [style=small black dot] (14) at (2, 4 * \triht + \downmid) {};
		\node [style=small black dot] (19) at (4, -2 * \triht + \downmid) {};
		\node [style=small black dot] (20) at (0.5, -4 * \triht + \upmid) {};
		\node [style=small black dot] (21) at (-3, -2 * \triht + \downmid) {};
		\node [style=small black dot] (22) at (-3, 1 * \triht + \upmid) {};
		\node [style=small black dot, label={right:$w$}] (23) at (5.5, 4 * \triht + \upmid) {};
		\node [style=small black dot, label={above:$w_1$}] (24) at (2, 6 * \triht + \downmid) {};
		\node [style=small black dot, label={left:$w_2$}] (25) at (-4.5, 2 * \triht + \upmid) {};
		\node [style=small black dot, label={left:$w_3$}] (26) at (-4.5, -3 * \triht + \downmid) {};
		\node [style=small black dot, label={below:$w_4$}] (27) at (0.5, -6 * \triht + \upmid) {};
		\node [style=small black dot, label={right:$w_5$}] (28) at (5.5, -3 * \triht + \downmid) {};
		\node [style=none] (29) at (6, 4 * \triht) {};
		\node [style=none] (30) at (5.5, 5 * \triht) {};
		\node [style=none] (31) at (4.5, 5 * \triht) {};
		\node [style=none] (32) at (3.5, 3 * \triht) {};
		\node [style=none] (33) at (5.5, 3 * \triht) {};
		\node [style=none] (34) at (2.5, 7 * \triht) {};
		\node [style=none] (35) at (1.5, 7 * \triht) {};
		\node [style=none] (36) at (1, 6 * \triht) {};
		\node [style=none] (37) at (2, 4 * \triht) {};
		\node [style=none] (38) at (3, 6 * \triht) {};
		\node [style=none] (39) at (-4.5, 3 * \triht) {};
		\node [style=none] (40) at (-5, 2 * \triht) {};
		\node [style=none] (41) at (-4.5, 1 * \triht) {};
		\node [style=none] (42) at (-2.5, 1 * \triht) {};
		\node [style=none] (43) at (-3.5, 3 * \triht) {};
		\node [style=none] (44) at (-5, -2 * \triht) {};
		\node [style=none] (45) at (-4.5, -3 * \triht) {};
		\node [style=none] (46) at (-3.5, -3 * \triht) {};
		\node [style=none] (47) at (-2.5, -1 * \triht) {};
		\node [style=none] (48) at (-4.5, -1 * \triht) {};
		\node [style=none] (49) at (0, -6 * \triht) {};
		\node [style=none] (50) at (1, -6 * \triht) {};
		\node [style=none] (51) at (1.5, -5 * \triht) {};
		\node [style=none] (52) at (0.5, -3 * \triht) {};
		\node [style=none] (53) at (-0.5, -5 * \triht) {};
		\node [style=none] (54) at (5.5, -3 * \triht) {};
		\node [style=none] (55) at (6, -2 * \triht) {};
		\node [style=none] (56) at (5.5, -1 * \triht) {};
		\node [style=none] (57) at (3.5, -1 * \triht) {};
		\node [style=none] (58) at (4.5, -3 * \triht) {};
		\node [style=purple dot, label={above:$p$}] (59) at (0.5, 5 * \triht + \downmid) {};
		\node [style=purple dot, label={above:$\tilde p$}] (60) at (-3, 3 * \triht + \upmid) {};
		\node [style=none] (61) at (4, 8 * \triht) {};
		\node [style=none] (62) at (-4, -8 * \triht) {};
		\node [style=orange dot, label={below:$p'$}] (63) at (4, -4 * \triht + \downmid) {};
		\node [style=orange dot, label={below:$\widetilde{p'}$}] (64) at (2, -5 * \triht + \upmid) {};
	\end{pgfonlayer}
	\begin{pgfonlayer}{edgelayer}
		\draw [style=red line] (2.center) to (3.center);
		\draw [style=green line] (7.center) to (8.center);
		\draw [style=blue line] (12.center) to (11.center);
		\draw [style=filled grey] (33.center)
			 to (29.center)
			 to (30.center)
			 to (31.center)
			 to (32.center)
			 to cycle;
		\draw [style=filled grey] (36.center)
			 to (37.center)
			 to (38.center)
			 to (34.center)
			 to (35.center)
			 to cycle;
		\draw [style=filled grey] (43.center)
			 to (39.center)
			 to (40.center)
			 to (41.center)
			 to (42.center)
			 to cycle;
		\draw [style=filled grey] (47.center)
			 to (48.center)
			 to (44.center)
			 to (45.center)
			 to (46.center)
			 to cycle;
		\draw [style=filled grey] (52.center)
			 to (53.center)
			 to (49.center)
			 to (50.center)
			 to (51.center)
			 to cycle;
		\draw [style=filled grey] (58.center)
			 to (54.center)
			 to (55.center)
			 to (56.center)
			 to (57.center)
			 to cycle;
		\draw [style=dashed line] (13) to (14);
		\draw [style=dashed line] (19) to (13);
		\draw [style=dashed line] (20) to (19);
		\draw [style=dashed line] (21) to (20);
		\draw [style=dashed line] (14) to (22);
		\draw [style=dashed line] (22) to (21);
		\draw (23) to (24);
		\draw (24) to (25);
		\draw (25) to (26);
		\draw (26) to (27);
		\draw (27) to (28);
		\draw (28) to (23);
		\draw [style=green line] (61.center) to (62.center);
		\draw [style=thick purple line] (59) to (60);
		\draw [style=orange line] (60) to (21);
		\draw [style=orange line] (21) to (20);
		\draw [style=orange line] (20) to (19);
		\draw [style=orange line] (19) to (13);
		\draw [style=orange line] (13) to (59);
		\draw [style=dotted line] (60) to (28);
		\draw [style=dotted line] (59) to (27);
	\end{pgfonlayer}


    \clip       (-\cols, -\height) rectangle (\cols, \height);
    \draw[gray] (-\cols, -\height) rectangle (\cols, \height);

    \pgfmathsetmacro{\from}{-2 *\cols}
    \pgfmathsetmacro{\to}{2 * \cols}
    \foreach\i in {\from, ..., \to} {
        \draw[xslant=\slant]  (\i, -\height) -- (\i, \height);
        \draw[xslant=-\slant] (\i, -\height) -- (\i, \height);
    }

    \foreach\j in {-\rows, ..., \rows} {
        \pgfmathsetmacro{\y}{0.5 * \j * tan(60)}
        \draw (-\cols, \y) -- (\cols, \y);
    }

\end{tikzpicture}
\caption{Special edge, Type 1}  \label{fig:specialEdgeType1}
\end{figure}

Recall that $E$ is determined by the intersection of a hexagon edge with the (crossing) diagonals from the vertices, say $A$ and $B$, of the opposite edge. (See the dotted lines in Figure \ref{fig:specialEdgeType1}.) Now $A$ and $B$ are reflections of each other across the root hyperplane $H=H_{\gg,k}$ used in the construction of $\ch_{w}$ (see the last two equations in \eqref{e:hexagonVertices}). Then $H$ is the perpendicular bisector of the segment $AB$. Since the diagonals from $A$ and $B$ both make $60^{\circ}$ angles with the segment $AB$, and intersect at a point $D$ in the interior of $\ch_{w}$, it is clear from similar triangles that the points $A'$ and $B'$ where the diagonals intersect the opposite (parallel!) edge of $\ch_{w}$ are also reflections of each other across the hyperplane $H$. In particular, one of them is in the Bruhat hexagon determined by the other. This implies that for one endpoint $p$ of $E$,  $x_{0}\le p$ and $q_{p}^{w}=1$. Likewise $x=x_{0}v\le p$. We will be done if we show all of $\ch_{p}$ is contained in the region where $q^{w}_{\bullet}>0$: namely, the $(4-\tau)$-shell and its interior when $w$ is of Type $\tau$, together with the alcoves ``near'' the special edge $E$.

Since a special segment exists, we know $w$ belongs to an odd chamber, and without loss of generality we may assume $w$ belongs to Chamber I. We will treat the special segment along the NW edge of $\ch_{w}$; the argument for the SE special segment is similar.  We will give the arguments for the cases where $w$ is not in a root strip adjacent to a chamber wall (as in Figures \ref{fig:3shell} and \ref{fig:2shell}), and leave the arguments for the cases where $w$ is in a root strip (as in Figures \ref{fig:edge3shell} and \ref{fig:fish}) to the reader.

\medskip\noindent
1. Suppose first that $w$ is of Type $\tau=1$, so $4-\tau=3$. Continue the notation and setup surrounding Figure \ref{fig:3shell}. In particular, we have the region $X$ and its five reflections across the hyperplanes which define the vertices of both $\ch_{w}$ and, as we observed earlier, its 3-shell $\ch_{z}$ (or $\ch_{z'}$ when $w$ is adjacent to $H_{\ga_{1},1}$). A key observation is that the South copy of X (containing the vertex $w_{4}$ of $\ch_{w}$) is offset from the North copy of X by exactly three $\ga_{2}$ (vertical) root strings to the left. The point $p$ is the intersection of the vertical string through $w_{4}$ with the NW edge of $\ch_{w}$; see Figure \ref{fig:specialEdgeType1}. Since $p$ belongs to Chamber II, it is closer to $H_{\ga_{1},0}$ than to $H_{\ga_{1},1}$ (whereas the reverse is true for $w$ in Chamber I). These observations imply that four of the vertices of $\ch_{p}$ (shown in orange in the figure) coincide with four of the vertices of the 3-shell (shown as dashed lines) of $\ch_{w}$. The remaining two vertices of $\ch_{p}$ are the endpoints $p, \tilde p$ of the special segment $E$ along the NE edge of $\ch_{w}$. Hence the alcoves on or inside $\ch_{p}$ consist of (i) the alcoves on or inside the 3-shell of $\ch_{w}$, and (ii) the alcoves $x_{0}v$ for $x_{0}$ on $E$ and $v\in R(w)$. All of these alcoves have $q^{w}_{\bullet}>0$, as desired.

\begin{figure}[htbp!]
%

\pgfmathsetmacro{\cols}{7}
\pgfmathsetmacro{\rows}{8}
\pgfmathsetmacro{\slant}{cot(60)}
\pgfmathsetmacro{\height}{0.5 * \rows * tan(60)}
\pgfmathsetmacro{\triht}{sin(60)}
\pgfmathsetmacro{\upmid}{0.25 * sec(30)}
\pgfmathsetmacro{\downmid}{\triht - \upmid}

\begin{tikzpicture}
    

\input{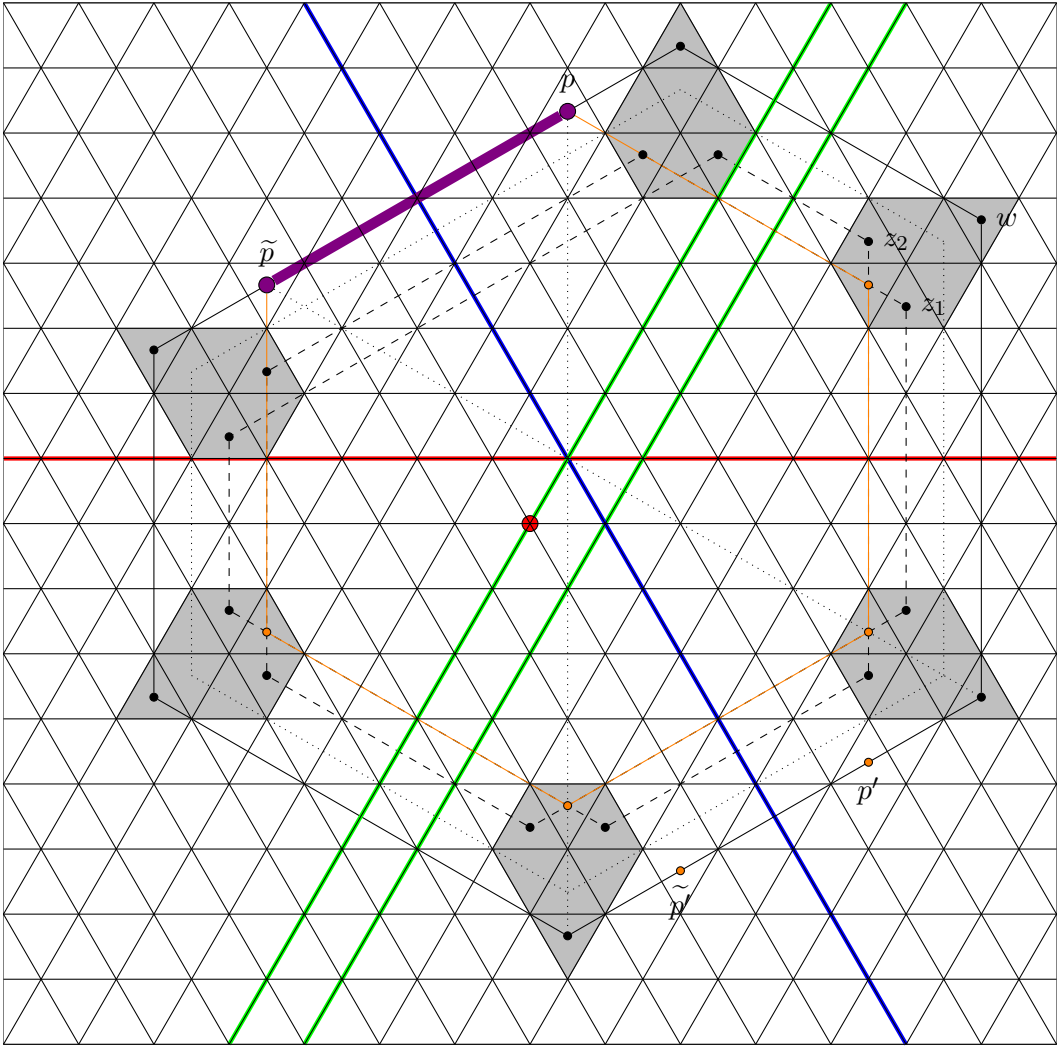}

	\begin{pgfonlayer}{nodelayer}
		\node [style=none] (2) at (-7, 1 * \triht) {};
		\node [style=none] (3) at (7, 1 * \triht) {};
		\node [style=red dot] (4) at (0, 0) {};
		\node [style=none] (7) at (-3, -8 * \triht) {};
		\node [style=none] (8) at (5, 8 * \triht) {};
		\node [style=none] (11) at (5, -8 * \triht) {};
		\node [style=none] (12) at (-3, 8 * \triht) {};
		\node [style=small black dot, label={right:$z_2$}] (13) at (4.5, 4 * \triht + \upmid) {};
		\node [style=small black dot] (14) at (2.5, 5 * \triht + \downmid) {};
		\node [style=small black dot] (19) at (4.5, -3 * \triht + \downmid) {};
		\node [style=small black dot] (20) at (1, -5 * \triht + \upmid) {};
		\node [style=small black dot] (21) at (-4, -2 * \triht + \downmid) {};
		\node [style=small black dot] (22) at (-4, 1 * \triht + \upmid) {};
		\node [style=small black dot, label={right:$w$}] (23) at (6, 4 * \triht + \downmid) {};
		\node [style=small black dot] (24) at (2, 7 * \triht + \upmid) {};
		\node [style=small black dot] (25) at (-5, 2 * \triht + \downmid) {};
		\node [style=small black dot] (26) at (-5, -3 * \triht + \upmid) {};
		\node [style=small black dot] (27) at (0.5, -7 * \triht + \downmid) {};
		\node [style=small black dot] (28) at (6, -3 * \triht + \upmid) {};
		\node [style=none] (29) at (5.5, 4 * \triht + \upmid) {};
		\node [style=none] (30) at (2, 6 * \triht + \downmid) {};
		\node [style=none] (31) at (-4.5, 2 * \triht + \upmid) {};
		\node [style=none] (32) at (-4.5, -3 * \triht + \downmid) {};
		\node [style=none] (33) at (0.5, -6 * \triht + \upmid) {};
		\node [style=none] (99) at (5.5, -3 * \triht + \downmid) {};
		\node [style=none] (65) at (5.5, 3 * \triht) {};
		\node [style=none] (66) at (6.5, 5 * \triht) {};
		\node [style=none] (67) at (4.5, 5 * \triht) {};
		\node [style=none] (68) at (4, 4 * \triht) {};
		\node [style=none] (69) at (4.5, 3 * \triht) {};
		\node [style=none] (34) at (3, 6 * \triht) {};
		\node [style=none] (35) at (2, 8 * \triht) {};
		\node [style=none] (36) at (1, 6 * \triht) {};
		\node [style=none] (37) at (1.5, 5 * \triht) {};
		\node [style=none] (38) at (2.5, 5 * \triht) {};
		\node [style=none] (39) at (-3.5, 3 * \triht) {};
		\node [style=none] (40) at (-5.5, 3 * \triht) {};
		\node [style=none] (41) at (-4.5, 1 * \triht) {};
		\node [style=none] (42) at (-3.5, 1 * \triht) {};
		\node [style=none] (43) at (-3, 2 * \triht) {};
		\node [style=none] (44) at (-4.5, -1 * \triht) {};
		\node [style=none] (45) at (-5.5, -3 * \triht) {};
		\node [style=none] (46) at (-3.5, -3 * \triht) {};
		\node [style=none] (47) at (-3, -2 * \triht) {};
		\node [style=none] (48) at (-3.5, -1 * \triht) {};
		\node [style=none] (49) at (-0.5, -5 * \triht) {};
		\node [style=none] (50) at (0.5, -7 * \triht) {};
		\node [style=none] (51) at (1.5, -5 * \triht) {};
		\node [style=none] (52) at (1, -4 * \triht) {};
		\node [style=none] (53) at (0, -4 * \triht) {};
		\node [style=none] (54) at (4.5, -3 * \triht) {};
		\node [style=none] (55) at (6.5, -3 * \triht) {};
		\node [style=none] (56) at (5.5, -1 * \triht) {};
		\node [style=none] (57) at (4.5, -1 * \triht) {};
		\node [style=none] (58) at (4, -2 * \triht) {};
		\node [style=small black dot, label={right:$z_1$}] (59) at (5, 3 * \triht + \upmid) {};
		\node [style=small black dot] (60) at (5, -2 * \triht + \downmid) {};
		\node [style=small black dot] (61) at (0, -5 * \triht + \upmid) {};
		\node [style=small black dot] (62) at (-3.5, -3 * \triht + \downmid) {};
		\node [style=small black dot] (63) at (-3.5, 2 * \triht + \upmid) {};
		\node [style=small black dot] (64) at (1.5, 5 * \triht + \downmid) {};
		\node [style=none] (100) at (4, 8 * \triht) {};
		\node [style=none] (101) at (-4, -8 * \triht) {};
		\node [style=purple dot, label={above:$p$}] (102) at (0.5, 6 * \triht + \upmid) {};
		\node [style=purple dot, label={above:$\tilde p$}] (103) at (-3.5, 3 * \triht + \downmid) {};
		\node [style=orange dot] (104) at (4.5, 3 * \triht + \downmid) {};
		\node [style=orange dot] (105) at (4.5, -2 * \triht + \upmid) {};
		\node [style=orange dot] (106) at (0.5, -5 * \triht + \downmid5) {};
		\node [style=orange dot] (107) at (-3.5, -2 * \triht + \upmid) {};
		\node [style=orange dot, label={below:$p'$}] (108) at (4.5, -4 * \triht + \upmid) {};
		\node [style=orange dot, label={below:$\widetilde{p'}$}] (109) at (2, -6 * \triht + \downmid) {};

	\end{pgfonlayer}
	\begin{pgfonlayer}{edgelayer}
		\draw [style=red line] (2.center) to (3.center);
		\draw [style=green line] (7.center) to (8.center);
		\draw [style=blue line] (12.center) to (11.center);
		\draw [style=filled grey] (69.center)
			 to (65.center)
			 to (66.center)
			 to (67.center)
			 to (68.center)
			 to cycle;
		\draw [style=filled grey] (36.center)
			 to (37.center)
			 to (38.center)
			 to (34.center)
			 to (35.center)
			 to cycle;
		\draw [style=filled grey] (43.center)
			 to (39.center)
			 to (40.center)
			 to (41.center)
			 to (42.center)
			 to cycle;
		\draw [style=filled grey] (47.center)
			 to (48.center)
			 to (44.center)
			 to (45.center)
			 to (46.center)
			 to cycle;
		\draw [style=filled grey] (52.center)
			 to (53.center)
			 to (49.center)
			 to (50.center)
			 to (51.center)
			 to cycle;
		\draw [style=filled grey] (58.center)
			 to (54.center)
			 to (55.center)
			 to (56.center)
			 to (57.center)
			 to cycle;
		\draw [style=dashed line] (13) to (14);
		\draw [style=dashed line] (19) to (13);
		\draw [style=dashed line] (20) to (19);
		\draw [style=dashed line] (21) to (20);
		\draw [style=dashed line] (14) to (22);
		\draw [style=dashed line] (22) to (21);
		\draw (23) to (24);
		\draw (24) to (25);
		\draw (25) to (26);
		\draw (26) to (27);
		\draw (27) to (28);
		\draw (28) to (23);
		\draw [style=dashed line] (59) to (60);
		\draw [style=dashed line] (60) to (61);
		\draw [style=dashed line] (61) to (62);
		\draw [style=dashed line] (62) to (63);
		\draw [style=dashed line] (63) to (64);
		\draw [style=dashed line] (64) to (59);
		\draw [style=dotted line] (29.center) to (30.center);
		\draw [style=dotted line] (30.center) to (31.center);
		\draw [style=dotted line] (31.center) to (32.center);
		\draw [style=dotted line] (32.center) to (33.center);
		\draw [style=dotted line] (33.center) to (99.center);
		\draw [style=dotted line] (99.center) to (29.center);
		\draw [style=green line] (100.center) to (101.center);
		\draw [style=thick purple line] (102) to (103);
		\draw [style=orange line] (103) to (107.center);
		\draw [style=orange line] (107.center) to (106.center);
		\draw [style=orange line] (106.center) to (105.center);
		\draw [style=orange line] (105.center) to (104.center);
		\draw [style=orange line] (104.center) to (102);
		\draw [style=dotted line] (102) to (27);
		\draw [style=dotted line] (103) to (28);
	\end{pgfonlayer}


    \clip       (-\cols, -\height) rectangle (\cols, \height);
    \draw[gray] (-\cols, -\height) rectangle (\cols, \height);

    \pgfmathsetmacro{\from}{-2 *\cols}
    \pgfmathsetmacro{\to}{2 * \cols}
    \foreach\i in {\from, ..., \to} {
        \draw[xslant=\slant]  (\i, -\height) -- (\i, \height);
        \draw[xslant=-\slant] (\i, -\height) -- (\i, \height);
    }

    \foreach\j in {-\rows, ..., \rows} {
        \pgfmathsetmacro{\y}{0.5 * \j * tan(60)}
        \draw (-\cols, \y) -- (\cols, \y);
    }

\end{tikzpicture}
\caption{Special edge, Type 2}  \label{fig:specialEdgeType2}
\end{figure}

\noindent
2. Suppose finally that $w$ is of Type $\tau=2$, so $4-\tau=2$. Continue the notation and setup surrounding Figure \ref{fig:2shell}. This time each of the four vertices of $\ch_{p}$ not at the endpoints of $E$ is either in the alcove between $z_{1}$ and $z_{2}$, or one of its reflections; see Figure \ref{fig:specialEdgeType2}. So the alcoves on or inside $\ch_{p}$ are all either (i) on or inside the 2-shell of $\ch_{w}$, or (ii) of the form $x_{0}v$ for $x_{0}$ on $E$ and $v\in R(w)$.  Again these all have $q^{w}_{\bullet}> 0$. \end{proof}

\begin{Cor} \label{c:trivial-rs}
If $w$ is not spiral, then $X_{w}$ is not rationally smooth at $xB$ if and only if $q_{x}^{w}>0$.
\end{Cor}

\begin{proof}
Theorem \ref{t:qHeredity} implies that $q^w_x = 0$ if and only if $q^w_y = 0$ for all $y$ with $x \leq y \leq w$, which by the Carrell-Peterson criterion is equivalent to
the statement that $X_{w}$ is rationally smooth at $x \cb$.
\end{proof}

\begin{Rem} \label{r:shell-rs}
This corollary, combined with the results of previous sections concerning the values of $q^w_x$, leads to a description of
 the nrs locus of $X_w$ in terms of the geometry of the Bruhat hexagon $\ch_w$.   If $w$ is not a base case, then Corollary \ref{c:qpositive} implies
 that the only rationally smooth $x$ are certain $x$ on the 0-, 1-, or 2- shell of $\ch_w$, and these $x$ are explicitly described in Theorem \ref{t:outershells}.   If $w$ is a base case, the nrs locus can be determined from Section \ref{s:basecases}.
\end{Rem}

\begin{Cor} \label{c:LookupConjecture}
The Lookup Conjecture is true for $\tilde A_{2}$.
\end{Cor}

\begin{proof}
If $w$ is a spiral element, then the Lookup Conjecture is true for $X_{w}$ by \cite[Theorem 10.6]{GrLi:21}. If $w$ is not spiral, then $X_{w}$ is not rationally smooth at $xB$ if and only if $q_{x}^{w}>0$, by Corollary \ref{c:trivial-rs}. So the Lookup Conjecture is trivially true for $X_{w}$.
\end{proof}

For a spiral element $w$, the maximal nrs $z<w$ are described in \cite[Corollary 10.5]{GrLi:21}.  For non-spiral elements, the maximal nrs $z<w$ are described in the next corollary (note that $X_w$ is rationally smooth if $w$ is twisted spiral).

\begin{Cor} \label{c:Maxnrs}
Assume that $w$ is not spiral or twisted spiral, in chamber $\cc$, of Type $\tau=1$ or $2$. Suppose $\ell(w) \geq 6$.  Let $\ga$ be the translation into $\cc$. When $\cc$ is odd, with (potential) special edges $w_{1}w_{2}$ and $w_{4}w_{5}$, let $\gb$ be the root such that $w_{2}=w_{1}+\gl\gb$ for some $\gl>0$. Set $p=w_{1}+\gb,\ p'=w_{5}+\gb$. The maximal nrs $z<w$ are as follows:
\begin{enumerate}
\item{$\cc$ is even, $\tau=1$}: $z=t(-\ga)w$ (Figure \ref{fig:3shell}).
\item{$\cc$ is even, $\tau=2$}: Generically, $z=z_{1}, z_{2}$ (Figure \ref{fig:2shell}). In Base Case 2, $z=z_{1}$ (Figure \ref{fig:fish}).
\item{$\cc$ is odd, $\tau=1$}: Generically, $z=p, p'$ (Figure \ref{fig:specialEdgeType1}). Omit $p$ or $p'$ if it is not a special segment endpoint (e.g.\ Figure \ref{fig:edge3shell}).   
\item{$\cc$ is odd, $\tau=2$}: Generically, $z=z_{1}, z_{2}, p, p'$ (Figure \ref{fig:specialEdgeType2}); omit $p$ or $p'$ if it is not a special segment endpoint. In Base Case 4, $z=z_{1}, z'_{2}, p$ (Figures \ref{fig:basecase4} and \ref{fig:edge2shell}).
\end{enumerate}
\end{Cor}

\begin{proof}
Combining Corollary \ref{c:qpositive}, Theorem \ref{t:outershells}, and the descriptions of $q^{w}_{\bullet}$ for the Base Cases (Section \ref{s:basecases}), it follows that $q^{w}_{x} > 0$ iff either: I. $x$ is on or inside the $(4-\tau)$-shell of $\ch_{w}$, or II. $x R(w) \cap E \ne \varnothing$ for some special segment $E$ of $\ch_{w}$ (and then necessarily $\cc$ is odd and $x$ belongs to the $j$-shell for some $0 \le j < 4-\tau$). 

There are at most two maximal elements $z$ in case I, denoted by $z$ (case I.1.(i)), $z'$ (case I.1.(ii)), $z_{1}, z_{2}$ (case I.2.(i)), $z_{1}$ (case I.2.(ii)(a)), $z_{1}, z'_{2}$ (case I.2.(ii)(b)) of the proof (and shown in the accompanying Figures \ref{fig:3shell}--\ref{fig:edge2shell}). 

Likewise there are at most two maximal elements $z$ in case II, namely the special segment endpoints $p$ (resp.\ $p'$) closer to $w_{1}$ (resp.\ $w_{5}$) along $w_{1}w_{2}$ (resp.\ $w_{4}w_{5}$) (if they exist). When $w$ is of Type 1 (in an odd chamber), $p>z$ (Case I.1.(i)) (resp.\ $p>z'$ (Case I.1.(ii)), so we can omit $z$ (resp.\ $z'$) from the list of maximal nrs elements.
\end{proof}

\begin{Rem} \label{r:Maxnrs}
The restriction $\ell(w) \geq 6$ in the preceding corollary is imposed to make the statement more uniform.
Removing this restriction yields a few additional elements $w$, which we can analyze as follows.  The other hypotheses of the corollary force $\ell(w) \geq 4$
(and in fact if $\ell(w) \leq 3$, then $X_w$ is smooth).  In cases (ii) and (iii), $\ell(w)$ is even, so the additional elements satisfy $\ell(w) = 4$; these elements are in Base Case 2 (for case (ii)) or Base Case 3 (for case (iii)), and the corresponding $X_w$ are smooth.
In cases (i) and (iv), $\ell(w)$ is odd.  The requirement that $w$ is not twisted spiral forces $\ell(w) \geq 7$ in case (i), so there are no additional elements to consider.
In case (iv), in each odd chamber, there is one additional element $w$; this element satisfies $\ell(w) = 5$.  Such an element is in Base Case 4; there are two nrs $z < w$: the maximal one
is $z=t(-\ga)w$, which has length $1$; the other is $z = e$.
\end{Rem}

\begin{Cor} \label{c:Codimnrs}
Maintain the setup and hypotheses of Corollary \ref{c:Maxnrs}. 

\noindent (a) The codimension of the nrs locus of $X_{w}$ is 4 in case (i), and 3 in cases (ii)--(iv). More precisely, the maximal nrs point $z$ has length $\ell(w)-4$ in case (i), and they all have length $\ell(w)-3$ in the remaining three cases. 

\noindent (b) The maximal nrs points $z$ have $q^{w}_{z}=1$, except for $z=t(-\ga)w$ in case (i) and $z=z'_{2}$ in Base Case 4, where $q^{w}_{z}=2$.
\end{Cor}

\begin{proof}
One checks immediately that the length of the maximal nrs $z<w$ is $\ell(w)-4$ in case (i), and they all have length $\ell(w)-3$ in the other three cases. The values of $q$ are given by Theorem \ref{t:outershells}, except in case (i), where $z$ lies on the 3-shell of $\ch_{w}$ and has $q^{w}_{z}=2$, as described in the last paragraph of Section \ref{ss:translation}.
\end{proof}

\begin{Rem} \label{r:Codimnrs}
Removing the restriction $\ell(w) \geq 6$ yields 3 additional elements $w$ for which $X_w$ is nrs, by Remark \ref{r:Maxnrs}---one in each odd chamber, of length 5.  For such a $w$, the maximal nrs point $z$
satisfies $\ell(z) = 1$, so the codimension of the nrs locus is $4$.  In this case $q^w_z = 2$.
\end{Rem}

\begin{Rem} \label{r:Maxnrs-spiral}
If $w$ is spiral of length at least $4$, then the maximal nrs points of $X_w$ are described in \cite[Cor.~10.5]{GrLi:21}.  Using the length formulas in
\cite[Prop.~6.9]{GrLi:21}, we see that the
codimension of the nrs locus of $X_{w}$ is 3.
\end{Rem}

\section{The locus of smooth points} \label{s:smooth}

In this section, we fix a non-spiral $w\in W$ and identify explicitly those points $x\le w$ for which $x\cb$ is a smooth point in the Schubert variety $X_{w}$. 

\begin{figure}[htbp!]
%

\pgfmathsetmacro{\cols}{6}
\pgfmathsetmacro{\rows}{7}
\pgfmathsetmacro{\slant}{cot(60)}
\pgfmathsetmacro{\height}{0.5 * \rows * tan(60)}
\pgfmathsetmacro{\triht}{sin(60)}
\pgfmathsetmacro{\upmid}{0.25 * sec(30)}
\pgfmathsetmacro{\downmid}{\triht - \upmid}

\begin{tikzpicture}
    

\input{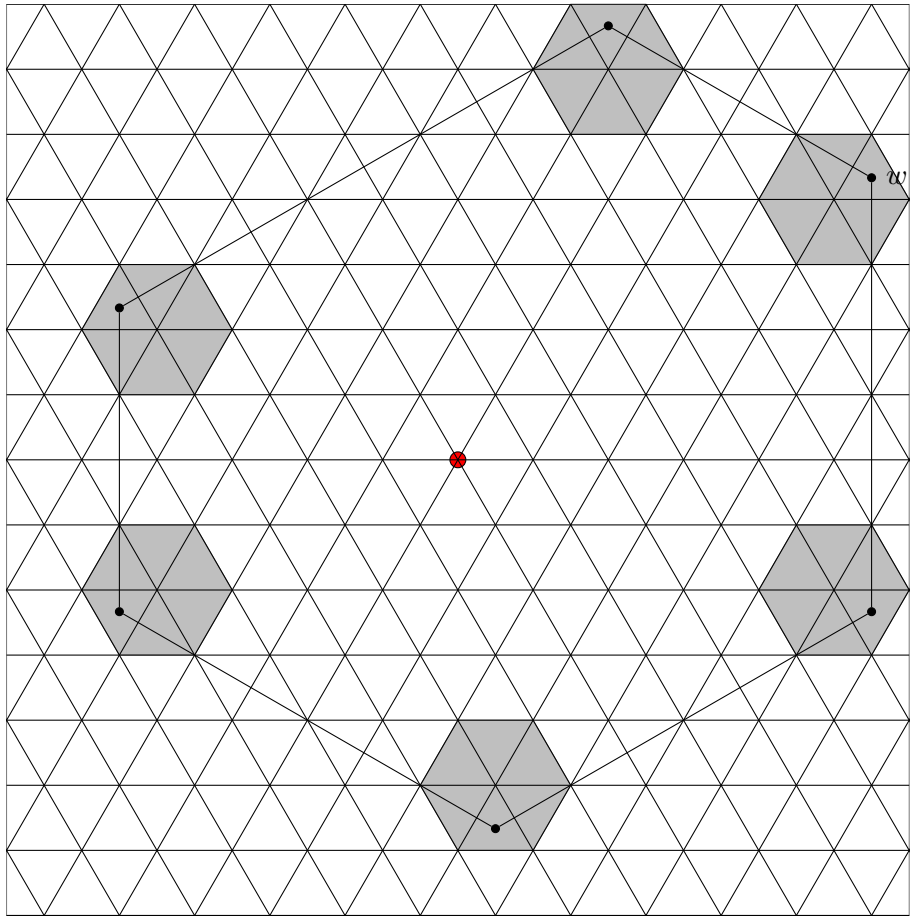}

	\begin{pgfonlayer}{nodelayer}
		\node [style=red dot] (4) at (0, 0) {};
		\node [style=small black dot, label={right:$w$}] (23) at (5.5, 4 * \triht + \upmid) {};
		\node [style=small black dot] (24) at (2, 6 * \triht + \downmid) {};
		\node [style=small black dot] (25) at (-4.5, 2 * \triht + \upmid) {};
		\node [style=small black dot] (26) at (-4.5, -3 * \triht + \downmid) {};
		\node [style=small black dot] (27) at (0.5, -6 * \triht + \upmid) {};
		\node [style=small black dot] (28) at (5.5, -3 * \triht + \downmid) {};
		\node [style=none] (29) at (6, 4 * \triht) {};
		\node [style=none] (30) at (5.5, 5 * \triht) {};
		\node [style=none] (31) at (4.5, 5 * \triht) {};
		\node [style=none] (32) at (4, 4 * \triht) {};
		\node [style=none] (33) at (5.5, 3 * \triht) {};
		\node [style=none] (34) at (2.5, 7 * \triht) {};
		\node [style=none] (35) at (1.5, 7 * \triht) {};
		\node [style=none] (36) at (1, 6 * \triht) {};
		\node [style=none] (37) at (1.5, 5 * \triht) {};
		\node [style=none] (38) at (3, 6 * \triht) {};
		\node [style=none] (39) at (-4.5, 3 * \triht) {};
		\node [style=none] (40) at (-5, 2 * \triht) {};
		\node [style=none] (41) at (-4.5, 1 * \triht) {};
		\node [style=none] (42) at (-3.5, 1 * \triht) {};
		\node [style=none] (43) at (-3.5, 3 * \triht) {};
		\node [style=none] (44) at (-5, -2 * \triht) {};
		\node [style=none] (45) at (-4.5, -3 * \triht) {};
		\node [style=none] (46) at (-3.5, -3 * \triht) {};
		\node [style=none] (47) at (-3, -2 * \triht) {};
		\node [style=none] (48) at (-4.5, -1 * \triht) {};
		\node [style=none] (49) at (0, -6 * \triht) {};
		\node [style=none] (50) at (1, -6 * \triht) {};
		\node [style=none] (51) at (1.5, -5 * \triht) {};
		\node [style=none] (52) at (1, -4 * \triht) {};
		\node [style=none] (53) at (-0.5, -5 * \triht) {};
		\node [style=none] (54) at (5.5, -3 * \triht) {};
		\node [style=none] (55) at (6, -2 * \triht) {};
		\node [style=none] (56) at (5.5, -1 * \triht) {};
		\node [style=none] (57) at (4.5, -1 * \triht) {};
		\node [style=none] (58) at (4.5, -3 * \triht) {};
		\node [style=none] (59) at (4.5, 3 * \triht) {};
		\node [style=none] (60) at (2.5, 5 * \triht) {};
		\node [style=none] (61) at (-3, 2 * \triht) {};
		\node [style=none] (62) at (-3.5, -1 * \triht) {};
		\node [style=none] (63) at (0, -4 * \triht) {};
		\node [style=none] (64) at (4, -2 * \triht) {};
	\end{pgfonlayer}
	\begin{pgfonlayer}{edgelayer}
		\draw [style=filled grey] (33.center)
			 to (29.center)
			 to (30.center)
			 to (31.center)
			 to (32.center)
			 to (59.center)
			 to cycle;
		\draw [style=filled grey] (36.center)
			 to (37.center)
			 to (60.center)
			 to (38.center)
			 to (34.center)
			 to (35.center)
			 to cycle;
		\draw [style=filled grey] (43.center)
			 to (39.center)
			 to (40.center)
			 to (41.center)
			 to (42.center)
			 to (61.center)
			 to cycle;
		\draw [style=filled grey] (47.center)
			 to (62.center)
			 to (48.center)
			 to (44.center)
			 to (45.center)
			 to (46.center)
			 to cycle;
		\draw [style=filled grey] (52.center)
			 to (63.center)
			 to (53.center)
			 to (49.center)
			 to (50.center)
			 to (51.center)
			 to cycle;
		\draw [style=filled grey] (58.center)
			 to (54.center)
			 to (55.center)
			 to (56.center)
			 to (57.center)
			 to (64.center)
			 to cycle;
		\draw (23) to (24);
		\draw (24) to (25);
		\draw (25) to (26);
		\draw (26) to (27);
		\draw (27) to (28);
		\draw (28) to (23);
	\end{pgfonlayer}


    \clip       (-\cols, -\height) rectangle (\cols, \height);
    \draw[gray] (-\cols, -\height) rectangle (\cols, \height);

    \pgfmathsetmacro{\from}{-2 *\cols}
    \pgfmathsetmacro{\to}{2 * \cols}
    \foreach\i in {\from, ..., \to} {
        \draw[xslant=\slant]  (\i, -\height) -- (\i, \height);
        \draw[xslant=-\slant] (\i, -\height) -- (\i, \height);
    }

    \foreach\j in {-\rows, ..., \rows} {
        \pgfmathsetmacro{\y}{0.5 * \j * tan(60)}
        \draw (-\cols, \y) -- (\cols, \y);
    }

\end{tikzpicture}
\caption{Smooth locus for $w$ of Type 1}  \label{fig:smoothtype1}
\end{figure}

\begin{figure}[htbp!]
%

\pgfmathsetmacro{\cols}{7}
\pgfmathsetmacro{\rows}{8}
\pgfmathsetmacro{\slant}{cot(60)}
\pgfmathsetmacro{\height}{0.5 * \rows * tan(60)}
\pgfmathsetmacro{\triht}{sin(60)}
\pgfmathsetmacro{\upmid}{0.25 * sec(30)}
\pgfmathsetmacro{\downmid}{\triht - \upmid}

\begin{tikzpicture}
    

\input{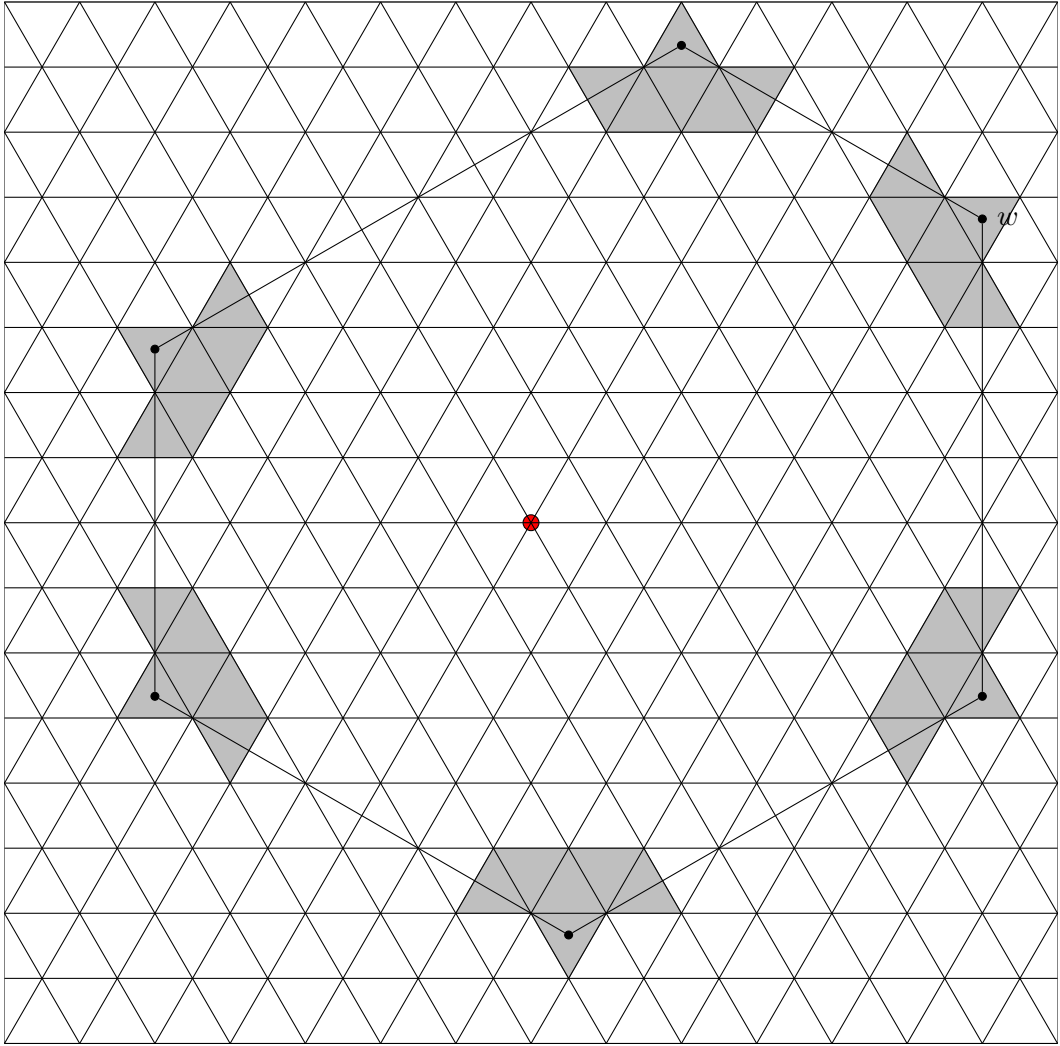}

	\begin{pgfonlayer}{nodelayer}
		\node [style=red dot] (4) at (0, 0) {};
		\node [style=small black dot, label={right:$w$}] (23) at (6, 4 * \triht + \downmid) {};
		\node [style=small black dot] (24) at (2, 7 * \triht + \upmid) {};
		\node [style=small black dot] (25) at (-5, 2 * \triht + \downmid) {};
		\node [style=small black dot] (26) at (-5, -3 * \triht + \upmid) {};
		\node [style=small black dot] (27) at (0.5, -7 * \triht + \downmid) {};
		\node [style=small black dot] (28) at (6, -3 * \triht + \upmid) {};
		\node [style=none] (29) at (6, 4 * \triht) {};
		\node [style=none] (30) at (2.5, 7 * \triht) {};
		\node [style=none] (31) at (-4, 4 * \triht) {};
		\node [style=none] (32) at (-4, -4 * \triht) {};
		\node [style=none] (33) at (2, -6 * \triht) {};
		\node [style=none] (99) at (6.5, -1 * \triht) {};
		\node [style=none] (65) at (5.5, 3 * \triht) {};
		\node [style=none] (66) at (6.5, 5 * \triht) {};
		\node [style=none] (67) at (5.5, 5 * \triht) {};
		\node [style=none] (68) at (5, 6 * \triht) {};
		\node [style=none] (69) at (4.5, 5 * \triht) {};
		\node [style=none] (34) at (3, 6 * \triht) {};
		\node [style=none] (35) at (2, 8 * \triht) {};
		\node [style=none] (36) at (1.5, 7 * \triht) {};
		\node [style=none] (37) at (0.5, 7 * \triht) {};
		\node [style=none] (38) at (1, 6 * \triht) {};
		\node [style=none] (39) at (-4.5, 3 * \triht) {};
		\node [style=none] (40) at (-5.5, 3 * \triht) {};
		\node [style=none] (41) at (-5, 2 * \triht) {};
		\node [style=none] (42) at (-5.5, 1 * \triht) {};
		\node [style=none] (43) at (-4.5, 1 * \triht) {};
		\node [style=none] (44) at (-5, -2 * \triht) {};
		\node [style=none] (45) at (-5.5, -3 * \triht) {};
		\node [style=none] (46) at (-4.5, -3 * \triht) {};
		\node [style=none] (47) at (-4.5, -1 * \triht) {};
		\node [style=none] (48) at (-5.5, -1 * \triht) {};
		\node [style=none] (49) at (0, -6 * \triht) {};
		\node [style=none] (50) at (0.5, -7 * \triht) {};
		\node [style=none] (51) at (1, -6 * \triht) {};
		\node [style=none] (52) at (-0.5, -5 * \triht) {};
		\node [style=none] (53) at (-1, -6 * \triht) {};
		\node [style=none] (54) at (5.5, -3 * \triht) {};
		\node [style=none] (55) at (6.5, -3 * \triht) {};
		\node [style=none] (56) at (6, -2 * \triht) {};
		\node [style=none] (57) at (4.5, -3 * \triht) {};
		\node [style=none] (58) at (5, -4 * \triht) {};
		\node [style=none] (100) at (6.5, 3 * \triht) {};
		\node [style=none] (101) at (3.5, 7 * \triht) {};
		\node [style=none] (102) at (-3.5, 3 * \triht) {};
		\node [style=none] (104) at (-3.5, -3 * \triht) {};
		\node [style=none] (105) at (1.5, -5 * \triht) {};
		\node [style=none] (106) at (5.5, -1 * \triht) {};
	\end{pgfonlayer}
	\begin{pgfonlayer}{edgelayer}
		\draw [style=filled grey] (66.center)
			 to (67.center)
			 to (68.center)
			 to (69.center)
			 to (65.center)
			 to (100.center)
			 to (29.center)
			 to cycle;
		\draw [style=filled grey] (30.center)
			 to (35.center)
			 to (36.center)
			 to (37.center)
			 to (38.center)
			 to (34.center)
			 to (101.center)
			 to cycle;
		\draw [style=filled grey] (31.center)
			 to (39.center)
			 to (40.center)
			 to (41.center)
			 to (42.center)
			 to (43.center)
			 to (102.center)
			 to cycle;
		\draw [style=filled grey] (44.center)
			 to (45.center)
			 to (46.center)
			 to (32.center)
			 to (104.center)
			 to (47.center)
			 to (48.center)
			 to cycle;
		\draw [style=filled grey] (33.center)
			 to (105.center)
			 to (52.center)
			 to (53.center)
			 to (49.center)
			 to (50.center)
			 to (51.center)
			 to cycle;
		\draw [style=filled grey] (99.center)
			 to (106.center)
			 to (57.center)
			 to (58.center)
			 to (54.center)
			 to (55.center)
			 to (56.center)
			 to cycle;
		\draw (23) to (24);
		\draw (24) to (25);
		\draw (25) to (26);
		\draw (26) to (27);
		\draw (27) to (28);
		\draw (28) to (23);
	\end{pgfonlayer}

    \clip       (-\cols, -\height) rectangle (\cols, \height);
    \draw[gray] (-\cols, -\height) rectangle (\cols, \height);

    \pgfmathsetmacro{\from}{-2 *\cols}
    \pgfmathsetmacro{\to}{2 * \cols}
    \foreach\i in {\from, ..., \to} {
        \draw[xslant=\slant]  (\i, -\height) -- (\i, \height);
        \draw[xslant=-\slant] (\i, -\height) -- (\i, \height);
    }

    \foreach\j in {-\rows, ..., \rows} {
        \pgfmathsetmacro{\y}{0.5 * \j * tan(60)}
        \draw (-\cols, \y) -- (\cols, \y);
    }

\end{tikzpicture}
\caption{Smooth locus for $w$ of Type 2}  \label{fig:smoothtype2}
\end{figure}

Recall that $R(w)$ is the subgroup of $W$ generated by the simple reflections which go down from $w$ on the right.  If $w$ is of Type $1$, then $R(w)$ is a subgroup of type $A_2$ generated by the two simple reflections it contains; if $w$ is of Type 2, then $R(w)$ has $2$ elements.

\begin{Thm} \label{t:smoothLocus}
Let $w\in W$ be non-spiral, and $x\le w$. Then $xB$ is a smooth point in the Schubert variety $X_{w}$ if and only if $x$ belongs to one of six alcoves ``near'' one of the six vertices of the hexagon $\ch_{w}$. More precisely, we have the following characterization of the smooth points. (Recall the notation $w_{i},\ 0\le i \le 5$, introduced in \ref{e:hexagonVertices}, for the six vertices of $\ch_{w}$.)
\begin{enumerate}
\item[(1)] Suppose $w$ is of Type 1. 
Then $x$ is smooth if and only if $x\in w_{i}R(w)$ for some $0\le i\le 5$. In other words, $x$ lies in one of the six alcoves of the small $A_{2}$ hexagon inside $\ch_{w}$ determined by some hexagon vertex $w_{i}$. See Figure \ref{fig:smoothtype1}.
\item[(2)] Suppose $w$ is of Type 2. Then $x$ is smooth if and only if, for some $0\le i\le 5$, either: (a) $x=w_{i}$; or (b) $x$ lies in one of the two edge alcoves of $\ch_{w}$ adjacent to $w_{i}$; or (c) $x$ lies in one of the three alcoves between the two alcoves described in (b). This region looks like a folded paper hat; see Figure \ref{fig:smoothtype2}.
\end{enumerate}
\end{Thm}

\begin{Rem} (1) There can be some overlap between the small hexagons or ``paper hats'' (i.e., sets of six smooth alcoves) near two hexagon vertices which are close together. 

\noindent (2) The smooth locus of a spiral Schubert variety has a similar description (see \cite{GrLi:15}), except that the ``hexagon'' $\ch_{w}$ degenerates to a quadrilateral, and some of the alcoves described lie outside the quadrilateral.
\end{Rem}

The proof, which will occupy the next four subsections, consists in showing that the points mentioned in the theorem are smooth, and that the remaining points are singular. We begin with several reductions.

\subsection{Reductions: Smooth Points} \label{ss:reductionsSmooth} 
Here we describe some easy facts which reduce the number of points which need to be proved smooth.

{\bf Reduction 1:} Suppose $w$ is of Type 1. To prove that all the points mentioned in the theorem are smooth, it suffices to prove that the vertices $w_{i}$ are smooth: the others, being in $w_{i}R(w)$, are then smooth by the Simple Move (on the right).

{\bf Reduction 2:} Suppose $w$ is of Type 2, with $w>ws$ for (unique) $s$ simple. To prove that all the points mentioned in the theorem are smooth, it suffices to prove that the vertices $w_{i}$, and the adjacent edge alcoves, say $y_{i}$ and $y'_{i}$, are smooth: the others, $w_{i}s, y_{i}s$, and $y'_{i}s$, are then smooth by the Simple Move.

{\bf Reduction 3:} Suppose $w$ lies in an even chamber. Then there are two simple reflections $s$ and $t$ such that $sw<w$ and $tw<w$, and these generate $L(w)$. Moreover, by Proposition \ref{p.diagonals}, the vertices of $\ch_{w}$ are the elements of the set $L(w) w$. Since $w$ is trivially a smooth point, it follows by Simple Moves on the left that all six vertices of $\ch_{w}$ are smooth. 

So if $w$ is of Type 1 (in an even chamber), all the points of Theorem \ref{t:smoothLocus} (1) are smooth, by Reduction 1. 

Assume $w$ is of Type 2 (in an even chamber), and write $y, y'$ for the edge alcoves of $\ch_{w}$ adjacent to $w$. It is easy to see that (independent of any assumptions about the type or the chamber of $w$), 
\begin{equation} \label{e:lengthAdjacentTow}
\ell(y) = \ell(y') = \ell(w)-1.
\end{equation}
By the codimension 1 condition (see Section 2.2), $y$ and $y'$ are smooth points. Applying the elements of $L(w)$ on the left of $y$ and $y'$ gives all the points $y_{i}$ and $y'_{i}$ ($0\le i\le 5$), so these are also smooth, by the Simple Move. Again we are done, by Reduction 2.

{\bf Reduction 4:} Suppose $w$ lies in an odd chamber. Then there is one simple reflection $s$ for which $s \ch_{w} = \ch_{w}$; it takes $w_{i}$ to $w_{j}$ for $(i,j)=(0,3), (1,2), (5,4)$. (See Figure \ref{fig:hexagon}, where this symmetry is reflection in the blue line, corresponding to $s_{0}$.)  Likewise $s$ takes the hexagon edge alcoves $y_{i}, y'_{i}$ adjacent to $w_{i}$ to the corresponding edge alcoves $y_{j}, y'_{j}$. 

If $w$ is of Type 1, using the Simple Move for $s$ on the left along with Reduction 1, it suffices to prove that $w_{1}$ and $w_{5}$ are smooth. But by Lemma \ref{l:lengthChange},
$$
\ell(w_{1}) = \ell(w_{5}) = \ell(w)-1,
$$
so these are smooth by the codimension 1 condition. Thus we are done verifying smooth points in this case. 

If $w$ is of Type 2, similar arguments to those in the previous paragraph, along with \eqref{e:lengthAdjacentTow}, show that we are reduced to proving $y_{i}$ and $y'_{i}$ are smooth when $i=1, 5$. We will do this in section \ref{ss:verificationSmooth}.

\subsection{Reductions: Singular Points} \label{ss:reductionsSingular} 
This subsection contains reductions for the number of points which need to be proved singular.  

{\bf Reduction $\mathbf{1'}$:} It suffices to prove that the edge alcoves $y$ of $\ch_{w}$ which are two alcoves away from the nearest vertex are singular. (Note that on a ``short edge''---containing only two or four alcoves---there is no such alcove $y$.)  To see this, it is enough to prove the following claim: if $z \leq w$ is one of the elements which needs to be proved singular, 
then either $q^w_z >0$ (in which case $z$ is nrs, whence certainly singular), or $z \leq y$ for one of these edge alcoves $y$, in which case the fact that $y$ is singular at $X_w$ implies
that $z$ is as well.  If $w$ is twisted spiral, then
there is one long edge attached to $w$, and if $x$ is the element two alcoves away from $w$ on this edge, then the Bruhat hexagon $\ch_x$
contains $z$.   Therefore assume $w$ is not twisted spiral.  Corollary \ref{c:qpositive} and the analysis of Base Cases 2, 3, and 4 imply that if $z$ is not on an outer shell of $\ch_w$ (by which we mean the $0$-shell and the $1$-shell, and if $w$ is Type 1, the $2$-shell
as well),  then $q^w_z>0$.  On the other hand, if $z$ lies on an outer shell of $\ch_w$, then the orbit $z R(w)$ intersects an edge of $\ch_w$ in a point at least two alcoves
away from the nearest vertex.  By the endpoint theorem, $z \leq y$ for some $y$ which is two alcoves away from the nearest vertex.  This proves the claim.

{\bf Reduction $\mathbf{2'}$:} Suppose $w$ lies in an even chamber. As in Reduction 3 above, using Simple Moves on the left, it suffices to prove that the edge alcoves $x, x'$ two away from $w$ are singular (assuming these are at least two alcoves away from the other endpoint of their edge).   
Note that if one of the edges attached to $w$ is short, so that (say) $x'$ does not exist, then by symmetry, alternate edges of $\ch_{w}$ are all short, and the reflections of $x$ give all the long-edge alcoves two away from the vertices. So in this situation it suffices to prove that $x$ is singular. 

For later use, note that this argument proves that for $w$ in an even chamber,
all the edge alcoves of $\ch_{w}$ which are two alcoves away from the nearest vertex are less than or equal to $x$ or (if it exists) $x'$.

{\bf Reduction $\mathbf{3'}$:} Suppose $w$ lies in an odd chamber. We claim that, just as in Reduction $2'$, it suffices to prove that the edge alcoves $x, x'$ two away from $w$ (assuming they are at least two away from the other edge endpoint) are singular. 

\begin{figure}[htbp!]
%

\pgfmathsetmacro{\cols}{6}
\pgfmathsetmacro{\rows}{7}
\pgfmathsetmacro{\slant}{cot(60)}
\pgfmathsetmacro{\height}{0.5 * \rows * tan(60)}
\pgfmathsetmacro{\triht}{sin(60)}
\pgfmathsetmacro{\upmid}{0.25 * sec(30)}
\pgfmathsetmacro{\downmid}{\triht - \upmid}

\begin{tikzpicture}
    

\input{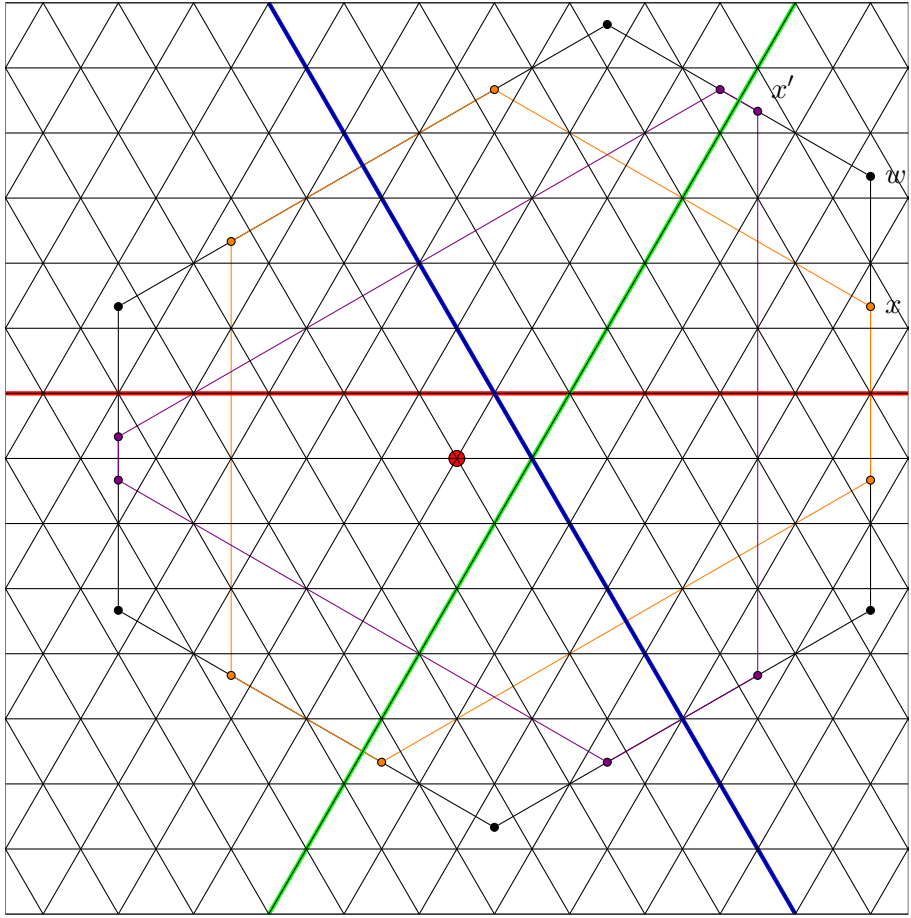}

	\begin{pgfonlayer}{nodelayer}
		\node [style=red dot] (4) at (0, 0) {};
		\node [style=small black dot, label={right:$w$}] (23) at (5.5, 4 * \triht + \upmid) {};
		\node [style=small black dot] (24) at (2, 6 * \triht + \downmid) {};
		\node [style=small black dot] (25) at (-4.5, 2 * \triht + \upmid) {};
		\node [style=small black dot] (26) at (-4.5, -3 * \triht + \downmid) {};
		\node [style=small black dot] (27) at (0.5, -6 * \triht + \upmid) {};
		\node [style=small black dot] (28) at (5.5, -3 * \triht + \downmid) {};
		\node [style=orange dot, label={right:$x$}] (65) at (5.5, 2 * \triht + \upmid) {};
		\node [style=orange dot] (66) at (5.5, -1 * \triht + \downmid) {};
		\node [style=orange dot] (67) at (-1, -5 * \triht + \upmid) {};
		\node [style=orange dot] (68) at (-3, -4 * \triht + \downmid) {};
		\node [style=orange dot] (69) at (-3, 3 * \triht + \upmid) {};
		\node [style=orange dot] (70) at (0.5, 5 * \triht + \downmid) {};
		\node [style=small purple dot, label={above right:$x'$}] (71) at (4, 5 * \triht + \upmid) {};
		\node [style=small purple dot] (72) at (3.5, 5 * \triht + \downmid) {};
		\node [style=small purple dot] (73) at (4, -4 * \triht + \downmid) {};
		\node [style=small purple dot] (74) at (-4.5,  \upmid) {};
		\node [style=small purple dot] (75) at (-4.5, -1 * \triht + \downmid) {};
		\node [style=small purple dot] (76) at (2, -5 * \triht + \upmid) {};
		\node [style=none] (77) at (-6, 1 * \triht) {};
		\node [style=none] (78) at (6, 1 * \triht) {};
		\node [style=none] (79) at (4.5, 7 * \triht) {};
		\node [style=none] (80) at (-2.5, -7 * \triht) {};
		\node [style=none] (81) at (-2.5, 7 * \triht) {};
		\node [style=none] (82) at (4.5, -7 * \triht) {};
	\end{pgfonlayer}
	\begin{pgfonlayer}{edgelayer}
		\draw (23) to (24);
		\draw (24) to (25);
		\draw (25) to (26);
		\draw (26) to (27);
		\draw (27) to (28);
		\draw (28) to (23);
		\draw [style=orange line] (65) to (66);
		\draw [style=orange line] (66) to (67);
		\draw [style=orange line] (67) to (68);
		\draw [style=orange line] (68) to (69);
		\draw [style=orange line] (69) to (70);
		\draw [style=orange line] (70) to (65);
		\draw [style=purple line] (71) to (72);
		\draw [style=purple line] (72) to (74);
		\draw [style=purple line] (74) to (75);
		\draw [style=purple line] (75) to (76);
		\draw [style=purple line] (76) to (73);
		\draw [style=purple line] (73) to (71);
		\draw [style=red line] (77.center) to (78.center);
		\draw [style=green line] (79.center) to (80.center);
		\draw [style=blue line] (81.center) to (82.center);
	\end{pgfonlayer}


    \clip       (-\cols, -\height) rectangle (\cols, \height);
    \draw[gray] (-\cols, -\height) rectangle (\cols, \height);

    \pgfmathsetmacro{\from}{-2 *\cols}
    \pgfmathsetmacro{\to}{2 * \cols}
    \foreach\i in {\from, ..., \to} {
        \draw[xslant=\slant]  (\i, -\height) -- (\i, \height);
        \draw[xslant=-\slant] (\i, -\height) -- (\i, \height);
    }

    \foreach\j in {-\rows, ..., \rows} {
        \pgfmathsetmacro{\y}{0.5 * \j * tan(60)}
        \draw (-\cols, \y) -- (\cols, \y);
    }

\end{tikzpicture}
\caption{Hexagons for the two maximal singular points $x, x'$.}  \label{fig:singularxx'}
\end{figure}

First, suppose both these alcoves exist. Then they are both in the same chamber as $w$. Thus the vertices of their hexagons $\ch_{x}, \ch_{x'}$ are defined using reflection in the same hyperplanes as the vertices of $\ch_{w}$ (see Figure \ref{fig:singularxx'}). These 12 vertices  thus lie at the centers of the alcoves two away from the 6 vertices of $\ch_{w}$. But since each vertex $z$ of $\ch_{x}$ satisfies $z\le x$, it follows that if $x$ is singular so is $z$ (and similarly for $x'$).

\begin{figure}[htbp!]
%

\pgfmathsetmacro{\cols}{6}
\pgfmathsetmacro{\rows}{6}
\pgfmathsetmacro{\slant}{cot(60)}
\pgfmathsetmacro{\height}{0.5 * \rows * tan(60)}
\pgfmathsetmacro{\triht}{sin(60)}
\pgfmathsetmacro{\upmid}{0.25 * sec(30)}
\pgfmathsetmacro{\downmid}{\triht - \upmid}

\begin{tikzpicture}
    

\input{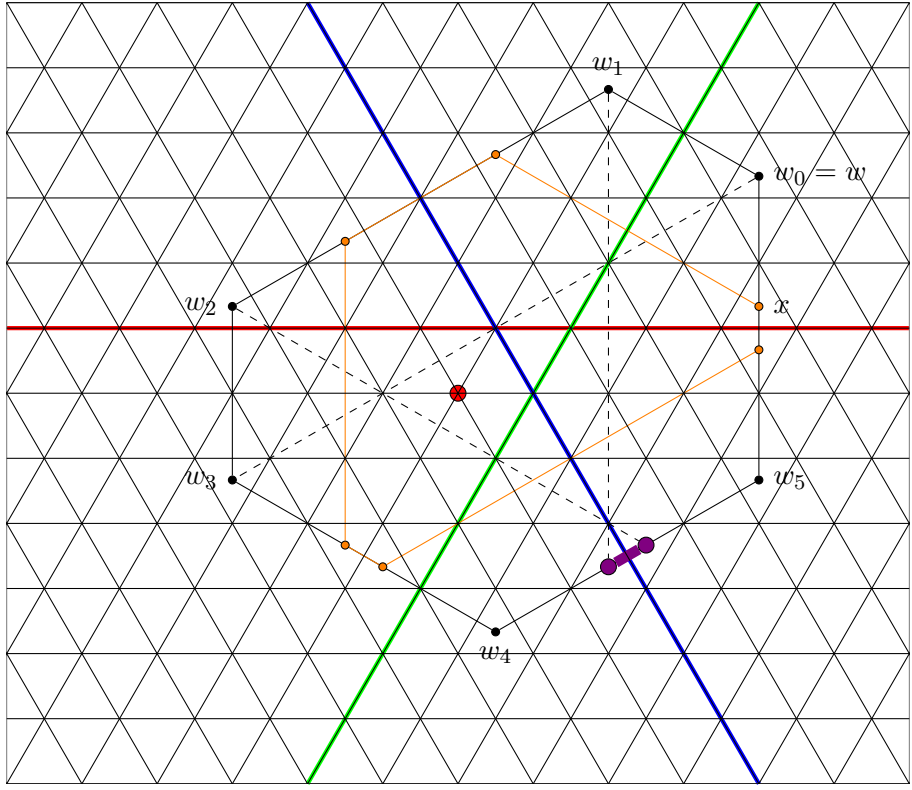}

	\begin{pgfonlayer}{nodelayer}
		\node [style=none] (2) at (-6, 1 * \triht) {};
		\node [style=none] (3) at (6, 1 * \triht) {};
		\node [style=red dot] (4) at (0, 0) {};
		\node [style=none] (7) at (-2, -6 * \triht) {};
		\node [style=none] (8) at (4, 6 * \triht) {};
		\node [style=none] (11) at (4, -6 * \triht) {};
		\node [style=none] (12) at (-2, 6 * \triht) {};
		\node [style=small black dot, label={right:$w_0=w$}] (13) at (4, 3 * \triht + \upmid) {};
		\node [style=small black dot, label={above:$w_1$}] (14) at (2, 4 * \triht + \downmid) {};
		\node [style=small black dot, label={right:$w_5$}] (19) at (4, -2 * \triht + \downmid) {};
		\node [style=small black dot, label={below:$w_4$}] (20) at (0.5, -4 * \triht + \upmid) {};
		\node [style=small black dot, label={left:$w_3$}] (21) at (-3, -2 * \triht + \downmid) {};
		\node [style=small black dot, label={left:$w_2$}] (22) at (-3, 1 * \triht + \upmid) {};
		\node [style=orange dot, label={right:$x$}] (23) at (4, 1 * \triht + \upmid) {};
		\node [style=orange dot] (24) at (0.5, 3 * \triht + \downmid) {};
		\node [style=orange dot] (25) at (-1.5, 2 * \triht + \upmid) {};
		\node [style=orange dot] (26) at (4, \downmid) {};
		\node [style=orange dot] (27) at (-1.5, -3 * \triht + \downmid) {};
		\node [style=orange dot] (28) at (-1, -3 * \triht + \upmid) {};
		\node [style=purple dot] (29) at (2, -3 * \triht + \upmid) {};
		\node [style=purple dot] (30) at (2.5, -3 * \triht + \downmid) {};
	\end{pgfonlayer}
	\begin{pgfonlayer}{edgelayer}
		\draw [style=red line] (2.center) to (3.center);
		\draw [style=green line] (7.center) to (8.center);
		\draw [style=blue line] (12.center) to (11.center);
		\draw (13) to (14);
		\draw (19) to (13);
		\draw (20) to (19);
		\draw (21) to (20);
		\draw (14) to (22);
		\draw (22) to (21);
		\draw [style=dashed line] (13) to (21);
		\draw [style=orange line] (23) to (24);
		\draw [style=orange line] (24) to (25);
		\draw [style=orange line] (25) to (27);
		\draw [style=orange line] (27) to (28);
		\draw [style=thick purple line] (29) to (30);
		\draw [style=dashed line] (22) to (30);
		\draw [style=dashed line] (14) to (29);
		\draw [style=orange line] (28) to (26);
	\end{pgfonlayer}
	

    \clip       (-\cols, -\height) rectangle (\cols, \height);
    \draw[gray] (-\cols, -\height) rectangle (\cols, \height);

    \pgfmathsetmacro{\from}{-2 *\cols}
    \pgfmathsetmacro{\to}{2 * \cols}
    \foreach\i in {\from, ..., \to} {
        \draw[xslant=\slant]  (\i, -\height) -- (\i, \height);
        \draw[xslant=-\slant] (\i, -\height) -- (\i, \height);
    }

    \foreach\j in {-\rows, ..., \rows} {
        \pgfmathsetmacro{\y}{0.5 * \j * tan(60)}
        \draw (-\cols, \y) -- (\cols, \y);
    }

\end{tikzpicture}

\caption{Hexagon for the maximal singular point $x$.}  \label{fig:singularx}
\end{figure}

If one of the edges attached to $w$ is short, so that (say) only $x$ exists, then the six vertices of $\ch_{x}$ (the orange dots in Figure \ref{fig:singularx}) are each two alcoves away from corners of $\ch_{w}$ on three long edges. If $x$ is not smooth, then neither are these other five points (because they are less than $x$ in the Bruhat order). By the hexagon symmetry mentioned in Reduction 4, there is a short edge attached to $w_{3}$ corresponding to the short edge attached to $w$. The remaining edge, say it is the one joining $w_{4}$ and $w_{5}$, may have length 6 or more. But if so, it is necessarily a special edge, containing a special segment which connects the alcoves two away from the vertices $w_{4}$ and $w_{5}$ (the thick purple segment in Figure \ref{fig:singularx}). By Theorem \ref{t:outershells} (b) and (d), $q=1$ at those points, so they are not even rationally smooth, hence certainly not smooth. Thus it suffices to show that $x$ is not smooth. 

\subsection{Final verification: Singular Points} \label{ss:rverificationSingular} 
As shown in Reductions $2'$ and $3'$, it remains to show that the $\ch_{w}$ edge alcoves $x, x'$ two away from $w$ are singular (assuming these alcoves are at least two away from the other endpoint of their edge). Assume for definiteness that $w$ lies in chamber III or IV, and factor $w=uv$ as a product of spiral elements (as in Remark \ref{r:spiralFactorizations}), where $u$ lies in the west fundamental half-strip between chambers III and IV (and is thus one of the spiral elements considered in \cite{GrLi:15}). 

Suppose first that $w$ is in chamber IV. The ``hexagon'' $\ch_{u}$ (which is degenerate: in fact, a quadrilateral) has south vertex $ru$ for a reflection $r$, where $rw=w_{1}$ is the south vertex of $\ch_{w}$. This implies that the southwest edges $E$ of $\ch_{u}$ and $E'$ of $\ch_{w}$ lie on parallel root strings and are the same length (measured by alcoves whose centers lie on each edge). According to Lemma \ref{l:edgeTranslation}, there is a fixed spiral element of $W$ which right-multiplies each alcove on $E$ to give the corresponding alcove on $E'$. Of course, this fixed spiral element is $v$, since $uv=w$ and $u, w$ are the northwest endpoints of $E, E'$ (respectively). Write $v=v_{1}\dots v_{m}$ where each $v_{j}$ is simple and the expression is reduced.

Assume that $E$ and $E'$ have length at least six (since otherwise there will be no alcove on $E'$ which is two alcoves away from a vertex). Consider the element $y$ on $E$ two alcoves away from $u$. (This is the point denoted $y_{2}$ in \cite[Figure 6]{GrLi:15}.) Then $y v =: x$ lies on $E'$ two alcoves away from $w$. Notice that $y v_{1}$ lies just outside the boundary of $\ch_{u}$ (and on the boundary of $\ch_{u v_{1}}$). This means that the pair $(y v_{1}, u v_{1})$ is obtained from the pair $(y,u)$ by the Setup Move.
This argument can be repeated in turn for each of $v_{2}. \dots, v_{m}$, showing that the pair $(x,w) = (yv, uv)$ is obtained from the pair $(y,u)$ by a sequence of
Setup Moves.  By \cite{GrLi:15} $y \cb$ is singular (although rationally smooth) in $X_{u}$.  Repeatedly applying Proposition \ref{p:setup} shows that
$x \cb$ is singular in $X_w$, as required.

We leave it to the reader to check that, using the other canonical spiral factorization of $w$ (recall Remark \ref{r:spiralFactorizations}),  analogous reasoning shows that the point $x'$ two alcoves away from $w$ on the west edge of $\ch_{w}$ is singular in $X_{w}$.

Now suppose $w$ is in chamber III. Almost identical reasoning to the above shows that right multiplication by $v$ carries the northwest edge $E$ of $\ch_{u}$, alcove by alcove, to the northwest edge $E'$ of $\ch_{w}$. Now $y$ is the point denoted $x_{2}$ in \cite[Figure 7]{GrLi:15}, and the argument proceeds as before to prove that the point $x$ two alcoves away from $w$ on $E'$ is singular in $X_{w}$. The point $x'$ two alcoves away from $w$ on the west edge of $\ch_{w}$ is handled by using the other canonical spiral factorization of $w$.

\subsection{Final verification: Smooth Points} \label{ss:verificationSmooth} By Reductions 1--4, it remains to show, for $w$ of Type 2 in an odd chamber, that the edge alcoves one away from the vertices $w_{1}$ and $w_{5}$ correspond to smooth points. Assume for definiteness that $w$ belongs to chamber III.  First consider
the element $y_{5}$, which is one alcove southwest of the north vertex $w_{5}$ of $\ch_{w}$.
By using the first canonical spiral factorization $w=uv$ as in the previous subsection, we have that right multiplication by $v$ is an alcove-by-alcove bijection from the northwest edge $E$ of $\ch_{u}$ to the northwest edge $E'$ of $\ch_{w}$.  Let $y$ denote the (smooth!) alcove $s_{\tilde\ga}x_{1}$ on the northwest edge of $\ch_{u}$ in \cite[Figure 7]{GrLi:15}. 
Reasoning as in the previous subsection shows that the pair $(y_{5},w)$ is obtained from the pair $(y,u)$ by a sequence of Setup Moves, so by Proposition \ref{p:setup}, the point $y_{5} \cb$ is smooth in $X_{w}$.

Using the other canonical spiral factorization of $w$, we deduce similarly that the point $y_{1}$ (one alcove north of the southwest vertex $w_{1}$ of $\ch_{w}$) is smooth in $X_{w}$.

It remains to check the smoothness at $y'_{1}$ and $y'_{5}$, which are just around the corners from the edges of $\ch_{w}$ attached to $w$. For these, we resort to a more laborious method involving subexpressions. First consider $y'_{5}$. Factor $w=uv$ as usual, with $u=u_{1}\dots u_{n}$ and $v=v_{1}\dots v_{m}$ reduced expressions. Our assumptions on $w$ imply that $n$ is odd and $m$ is even. Write $u_{\hat \imath}$ for $u$ with $u_{i}$ removed, $u_{\hat \imath \hat \jmath}$ for $u$ with $u_{i}$ and $u_{j}$ removed, and similarly for $v$. 

The subexpression $u_{\widehat 2} v_{\widehat{m-1}}$ of $w=uv$ gives a reduced expression for $y'_{5}$.  Indeed, the path for $u_{\widehat 2}$ heads NW and then NE, and then $v_{\widehat{m-1}}$ heads NW and finally NE to $y'_{5}$; Proposition \ref{p:reducedPath} implies that such an expression 
for an element like $y'_{5}$ in chamber II is reduced. We claim that this is the \emph{only} subexpression of $w$ that multiplies to $y'_{5}$.

First, it follows from the preceding paragraph that $\ell(y'_{5})=\ell(w)-2$. So $y'_{5}$ cannot be obtained by deleting exactly one simple factor from the reduced expression for $w$: the length of such an element would have  opposite parity to $\ell(w)$. Since the length of a product in any Coxeter group is at most the sum of the lengths of the factors, the only possible forms for a reduced expression of $w$ having length $\ell(w)-2$ are $u v_{\hat \imath \hat \jmath}$, $u_{\hat \imath \hat \jmath}v$, or $u_{\hat \imath}v_{\hat \jmath}$, where each subexpression of $u$ or $v$ has colength equal to the number of deleted factors.

Now, we cannot have $y'_{5} = u v_{\hat \imath \hat \jmath}$, since the path for $u$ goes west, whereas by Proposition \ref{p:reducedPath}, the path for any reduced expression of an element (such as $y'_{5}$) in chamber II can only go NW or NE.

Suppose $y'_{5}=u_{\hat \imath \hat \jmath}v$, and put $z:=u_{\hat \imath \hat \jmath} = y'_{5} v^{-1}$. It is not hard to check that $z$ is (at least two alcoves) due west of $y'_{5}$ (recall that $\ell(v)$ is even, hence is at least 2). But this means that $y'_{5}=z v = u_{\hat \imath \hat \jmath} v$ is a reduced expression whose path ends with a segment heading east. This is impossible for an element in chamber II. 

The only remaining possibility is $y'_{5} = u_{\hat \imath}v_{\hat \jmath}$, where $u_{\hat \imath}$ has length $n-1$ and $v_{\hat \jmath}$ has length $m-1$. According to \cite[Proposition 5.2]{GrLi:15}, this implies that $i\in\{1,2,n-1,n\}$ (and moreover there is a unique subexpression of $u$ multiplying to $u_{\hat \imath}$). We can immediately rule out $i=1,\ n-1$, and $n$ because these subexpressions of $u$ correspond to paths heading southeast, west, and west, respectively, which are impossible for an element in chamber II. That leaves $y'_{5}=u_{\widehat {2}} v_{\hat\jmath}$. But we already know $y'_{5}=u_{\widehat 2} v_{\widehat{m-1}}$, whence $v_{\hat\jmath} = v_{\widehat{m-1}}$. Using \cite[Proposition 5.2]{GrLi:15} again, we conclude $j=m-1$. 

Knowing that there is a unique subexpression of $w=uv$ multiplying to $y'_{5}$, along with the fact (Theorem \ref{t:outershells}) that $y'_{5}$ is rationally smooth, we deduce from \cite[Theorem 2.3]{GrLi:15} that in fact $y'_{5}$ is smooth in $X_{w}$. 

The verification for $y'_{1}$ is completely analogous, using the other canonical spiral factorization of $w$, and is left to the reader.

This completes the proof of Theorem \ref{t:smoothLocus}.  \hfill $\Box$

The next corollary describes the maximal singular points of $X_w$.  In this corollary, we will use the following notation.  If both of the edges of $\ch_w$ attached to $w$ are at least 6 alcoves long, then $x, x'$ will denote the points on these edges on these edges which are two alcoves away from $x$.  If exactly one of these edges
is at least 6 alcoves long, let $x$ be the element on this edge two alcoves away from $w$ (and we say $x'$ does not exist).  We also use the notation of Corollary
\ref{c:Maxnrs} and Remarks \ref{r:Maxnrs} and \ref{r:Maxnrs-spiral}, where the maximal nrs points are described.  

\begin{Cor} \label{c:Maxsing}
The maximal singular points of $X_w$ are as follows.
If both of the edges of $\ch_w$ attached to $w$ are at least 6 alcoves long, then the maximal singular points are $x$ and $x'$.
If neither of these edges
is at least 6 alcoves long, then the maximal singular points coincide with the maximal nrs points, described in Corollary \ref{c:Maxnrs} for non-spiral elements.
Finally, suppose exactly one of these edges is at least 6 elements long.  If the other edge attached to $w$ has length $2$ (which means $w$ is a base case), then in Base Case $4$, the maximal singular points are
$x$ and $z_2'$; otherwise, $x$ is the maximal singular point.  If the other edge attached to $w$ has length $4$, the maximal singular points are as follows (cf.~Corollary \ref{c:Maxnrs}):
\begin{enumerate}
\item{$\cc$ is even, $\tau=1$}: The maximal singular point is $x$, since $t(-\ga)w < x$ (Figure \ref{fig:3shell}).
\item{$\cc$ is even, $\tau=2$}: The maximal singular points are $x$ and $z_1$, where
$z_1$ is the maximal nrs point which is not less than $x$ (Figure \ref{fig:2shell}). 
\item{$\cc$ is odd, $\tau=1$}: The maximal singular points are $x$ and $p$, where $p$ is the maximal nrs point which is not less than $x$
 (Figure \ref{fig:specialEdgeType1}).
\item{$\cc$ is odd, $\tau=2$}: The maximal singular points are $x, z_1$, and $p$, where $z_1$ and $p$ are the maximal nrs points not less than $x$
 (Figure \ref{fig:edge2shell}).
\end{enumerate}
\end{Cor}

\begin{proof}
If $w$ is spiral of length less than $6$, the result can be checked directly.  For spiral $w$ of length at least 6,
the singular locus is described in \cite{GrLi:15}, and one can verify that the set of singular points is contained in
the union of $\ch_x$ and $\ch_{x'}$.  Therefore assume that $w$ is not spiral.
It follows from Reductions $1'$, $2'$ and $3'$ in the proof of Theorem \ref{t:smoothLocus} that if $z$ is a singular point of $X_w$, then
either $z \leq x$ or $z \leq x'$ (if $x'$ exists), or $z$ is nrs.  Therefore, the maximal singular points of $X_w$ are $x$ and (if it exists) $x'$,
together with the maximal nrs points which are not
less than $x$ or $x'$.   The maximal nrs points are described in Corollary
\ref{c:Maxnrs} and Remarks \ref{r:Maxnrs} and \ref{r:Maxnrs-spiral}, and the corollary follows by comparing these elements with $x$ and $x'$.
\end{proof}

\begin{Rem} \label{r:Max-Billey-Crites}
In the proof of their main result identifying the Schubert varieties of type $\tilde{A}_n$ which are rationally smooth at every point, Billey and Crites identify some nrs points in 
Schubert varieties $X_w$ if $w$, when viewed as an infinite permutation, contains the pattern $3412$.  They conjecture that these points are maximal singular points.  We have not checked whether the points they identify are maximal singular.
\end{Rem}

\begin{Rem} \label{r:finitesmall}
There are only 64 elements $w$ such that no edge of $\ch_w$ attached to $w$ has 6 or more alcoves.  Exactly 31 of these $w$ correspond to smooth Schubert varieties (see below).
\end{Rem}

\begin{Cor} \label{c:SingularLocuscodim}
Let $X_{w}$ be a Schubert variety which is not smooth. Then the singular locus of $X_{w}$ has codimension 2, unless both of the edges of $\ch_w$ attached
to $w$ have length at most 4, in which case the singular locus coincides with the nrs locus (whose codimension is given in Corollary \ref{c:Codimnrs} and
Remarks \ref{r:Codimnrs} and \ref{r:Maxnrs-spiral}).
\end{Cor}

\begin{proof}
First suppose $w$ is not spiral.
When $X_{w}$ is singular, at least one of the edge alcoves $x$ or $x'$ two away from $w$ is at least two alcoves away from the other vertex of its edge in the hexagon $\ch_{w}$. As shown in the proof of the theorem, $x$ and/or $x'$ are the singular points in $X_{w}$ having largest length, and they have length $\ell(w)-2$. 
If $w$ is spiral, then the results of the preceding two sentences continue to hold, by the description of the smooth locus in \cite{GrLi:15}.
\end{proof}

\begin{Cor} \label{c:dimsmooth}
Suppose $x \cb$ is a smooth point of $X_w$.  We have $\ell(x) \geq \ell(w) - 6$.  More precisely, for non-spiral $w$, we have:

\noindent (a) If $w$ is of Type 1 in an even chamber, then $\ell(x) \geq \ell(w) - 6$.

\noindent (b) If $w$ is of Type 1 in an odd chamber, then $\ell(x) \geq \ell(w) - 5$.

\noindent (c) If $w$ is of Type 2 in an even chamber, then $\ell(x) \geq \ell(w) - 5$.

\noindent (d) If $w$ is of Type 2 in an odd chamber, then $\ell(x) \geq \ell(w) - 4$.
\end{Cor}

\begin{proof}
If $w$ is spiral, the inequality $\ell(x) \geq \ell(w) - 6$ can be deduced from the description of smooth points given in \cite{GrLi:15}.  In fact, if $w$ is an even (resp.~odd) length
spiral element, then one can show that $\ell(x) \geq \ell(w) -3$ (resp.~ $\ell(x) \geq \ell(w)-4$); we omit the proof.  Therefore  
assume that $w$ is not spiral.  If $x$ is a smooth point near the vertex $w_i$, and $w$ is of Type $\tau$,
then $\ell(x) \geq \ell(w_i) - 3$ if $\tau=1$, and $\ell(x) \geq \ell(w_i)-2$ if $\tau=2$.
Using Lemma \ref{l:lengthChange} one can see that if $w$ is in an odd chamber, then $\ell(w_i) \geq \ell(w) - 2$; if $w$ is in an even chamber, then $\ell(w_i) \geq \ell(w) - 3$.
Considering the possible combinations of type and chamber gives the result.
\end{proof}

The inequalities in the previous corollary are sharp.  For example, the inequality in (a) is sharp, since if $w$ is a Type 1 element of length $7$ or more in an even chamber $\cc$, then $X_w$ will have a smooth point $x \cb$ where $\ell(x) = \ell(w) - 6$.  The point $x$ can be obtained as follows: Let $w_{3}$  as usual be the hexagon vertex opposite $w$, and then let $x$ be the shortest length element in the small $A_2$ hexagon containing $w_{3}$.
Similar examples show that the other inequalities are sharp as well.

In Corollary \ref{c:ratsmoothSchubertvar}, we determined the rationally smooth Schubert varieties.  We can now determine the smooth Schubert varieties.
The previous corollary guarantees that there is no smooth Schubert variety of dimension greater than $6$; in fact, we see that the largest dimension of a smooth Schubert variety is 5.

\begin{Cor} \label{c:smoothSchubertvar}
The Schubert variety $X_{w}$ is smooth if and only if it is rationally smooth and $\ell(w) \leq 5$.
\end{Cor}

\begin{proof}
Since any smooth variety is rationally smooth, we may assume $X_{w}$ is rationally smooth, so $w$ is as in Corollary \ref{c:ratsmoothSchubertvar}.  
First observe that for $w$ as in cases (a), (b), or (c) of that corollary, $X_{w}$ is smooth.  Indeed, $X_{w}$ is smooth if and only if
$e \cb$ is a smooth point of $X_{w}$, and using Kumar's criterion, this can be verified by direct calculation for $w$ in these cases.  Alternatively, one can show that for non-spiral $w$, the smooth points of $X_{w}$ described in Theorem \ref{t:smoothLocus} include all the $x \leq w$.  We may therefore assume that $w$ is twisted spiral (case (d)).
If $\ell(w) \leq 5$, then the smooth points described in Theorem \ref{t:smoothLocus} include all the $x \leq w$, but this fails if $\ell(w) \geq 7$.
This completes the proof.
\end{proof}

Below we list explicitly the $w$ for which $X_{w}$ is smooth. (In the second column, $i, j$, and $k$ are distinct in $\{0,1,2\}$.) We see that there are exactly 31 such $w$.

\begin{table}[h] \label{t:smooth}
\centering
\caption{Smooth Schubert varieties $X_{w}$ in $\tilde A_{2}$.}
\renewcommand{\arraystretch}{1.5}
\begin{tabular}{| c | c | c |}
 \hline 
Length & Reduced Expr.\ for $w$ & Number of Such  \\  \hline
0 &  $e$ & 1 \\ \hline
1 &  $s_i$ & 3 \\ \hline
2 &  $s_i s_j$  & 6 \\ \hline
3 &  $s_i s_j s_i$  &  3 \\ \hline
3 &  $s_i s_j s_k$ &  6  \\ \hline
4 &  $s_i s_j s_i s_k$ &  3    \\ \hline
4 &  $s_k s_i s_j s_i$   &  3  \\ \hline
5 &  $s_i s_j s_k s_i s_k$ & 6 \\ \hline
 \end{tabular} \\[.4in]
\end{table}

 \begin{Rem} \label{r:BilleyCrites}
 Results analogous to Corollaries \ref{c:ratsmoothSchubertvar} and \ref{c:smoothSchubertvar} hold in type $\tilde{A}_n$ for general $n$, due to work
 of Billey and Crites \cite{BiCr:12} and Richmond and Slofstra \cite{RiSl:16}.   Their results use the language of pattern avoidance, using the
 realization of the Weyl group of type $\tilde{A}_n$ as a group of infinite permutations, so the statements are somewhat different.
The results are as follows.  Billey and Crites proved that $X_w$ is rationally smooth if and only if either (the infinite permutation corresponding to)
$w$ avoids the patterns $3412$ and $4231$ or $w$ is twisted spiral.  They conjectured, and verified for $n \leq 5$, that $X_w$ is smooth if and only if  $w$ avoids the patterns $3412$ and $4231$; this conjecture was proved for all $n$ in \cite{RiSl:16}.  Together, these results imply that in type $\tilde{A}_n$, $X_w$ is rationally smooth if and only if either $X_w$ is smooth or $w$ is twisted spiral.  For type $\tilde{A}_2$, this also follows from Corollaries \ref{c:ratsmoothSchubertvar} and \ref{c:smoothSchubertvar} of this paper.  For
small values of $n$, Billey and Crites state the number of $w$ which avoid the patterns $3412$ and $4231$; for $n=2$, there are 31 such permutations.  Thus, their work implies that there are
exactly 31 smooth Schubert varieties of type $\tilde{A}_2$.  This agrees with the count obtained in Table \ref{t:smooth}.  We have not attempted to verify combinatorially that the list from Table \ref{t:smooth} matches the list of permutations that avoid these patterns.
\end{Rem}

\bibliographystyle{alpha}	
\bibliography{library}	

\vfill

\end{document}